\documentclass[10pt]{article}

\usepackage{amssymb,amsmath,amsthm,amsfonts}
\usepackage[pdftex,pagebackref,colorlinks,linktocpage,citecolor=black,linkcolor=black]{hyperref}
\usepackage{graphicx}
\usepackage{epstopdf}
\usepackage{graphics}
\usepackage{subfigure}
\usepackage{bm}
\usepackage{multirow,multicol}
\usepackage{comment}
\usepackage{booktabs}
\usepackage{threeparttable}
\usepackage[ruled]{algorithm2e}
\graphicspath{{./pics/}}
\usepackage{verbatim}
\usepackage{indentfirst}
\usepackage{mathtools}
\usepackage[usenames,dvipsnames,svgnames,table]{xcolor}
\usepackage{cite}
\usepackage{soul}

\allowdisplaybreaks

\numberwithin{equation}{section}
\newtheorem{theorem}{Theorem}[section]
\newtheorem{lemma}{Lemma}[section]
\theoremstyle{remark}
\newtheorem{remark}{Remark}[section]

\def\figwid{0.25\textwidth}
\def\cT {\mathcal T}
\def\bold {\boldsymbol}
\def\Om {\Omega}
\def\J {\mathcal{J}}

\def\to{\rightarrow}

\def\Atilde {\widetilde{\mathcal A}}

\def\strain {\mathcal{\varepsilon}}
\def\stress {\mathcal{\sigma}}
\def\C {\mathcal{C}}

\newcommand{\dx}{\,{\rm d}x}

\newcommand{\dd}{\,{\rm d}}

\begin{document}

\title{An Adaptive Phase-Field Method for Structural Topology Optimization\thanks{The work of B. Jin is supported by Hong Kong RGC General Research
Fund (Project 14306423), and a start-up fund and Direct Grant for Research 2022/2023, both from The Chinese University of Hong Kong. The work of Y. Xu is supported in part by the National Natural Science Foundation of China (12250013, 12261160361 and 12271367), the Science and Technology Commission of Shanghai Municipality (20JC1413800 and 22ZR1445400) and a General Research Fund (KF202318) from Shanghai Normal University. The work of S. Zhu is supported in part by the National Key Basic Research Program (2022YFA1004402), the National Natural Science Foundation of China (12071149) and the Science and Technology Commission of Shanghai Municipality (22ZR1421900 and 22DZ2229014).}}
\author{Bangti Jin\thanks{Department of Mathematics, The Chinese University of Hong Kong, Shatin, New Territories, Hong Kong, P.R. China (b.jin@cuhk.edu.hk, bangti.jin@gmail.com)}
\and Jing Li\thanks{School of Mathematical Sciences, East China Normal University, Shanghai 200241, P.R. China. (betterljing@163.com)} \and Yifeng Xu\thanks{Department of Mathematics and Scientific Computing Key Laboratory of Shanghai Universities, Shanghai Normal University, Shanghai 200234, P.R. China. (yfxu@shnu.edu.cn, mayfxu@gmail.com); Corresponding author.}
\and Shengfeng Zhu\thanks{School of Mathematical Sciences, Key Laboratory of MEA (Ministry of Education) \& Shanghai Key Laboratory of PMMP, East China Normal University, Shanghai 200241, P.R. China. (sfzhu@math.ecnu.edu.cn)}}
\date{\today}
\maketitle
	
\begin{abstract}
In this work, we develop an adaptive algorithm for the efficient numerical solution of the minimum compliance problem in topology optimization. The algorithm employs the phase field
approximation and continuous density field. The adaptive procedure is driven by two
residual type a posteriori error estimators, one for the state variable and the other for the first-order optimality condition of the
objective functional. The adaptive algorithm is provably convergent in the sense that the
sequence of numerical approximations generated by the adaptive algorithm contains a
subsequence convergent to a solution of the continuous first-order optimality system. We provide
several numerical simulations to show the distinct features of the algorithm.\\
\textbf{Keywords}: minimum compliance, adaptive algorithm, topology optimization, convergence, a posteriori error estimator
	\end{abstract}
	
\section{Introduction}
	
This work is concerned with effective numerical algorithms for the classical minimum compliance problem in structural topology optimization, which is a highly active research area in engineering. It mainly deals with the optimal distribution of materials over a given design domain. Over the last few decades, this area has witnessed significant advances, and	found a broad range of practical applications, e.g., in mechanics, heat transfer, acoustics and electromagnetics; see, e.g., the reviews \cite{EschenauerOlhoff:2001,Razvany:2009} and the monographs \cite{BendsoeSigmund:2003,NovotnySokołowski:2013} on various theoretical and practical issues of the discipline.
	
The minimum compliance problem, first introduced by the pioneering work \cite{BendsoeKikuchi:1988}, is one
fundamental formulation in topology optimization. The goal is to find the optimal distribution of isotropic material
that minimizes the work of the external forces at equilibrium, with a prescribed amount of volume. Often
an optimization problem is formulated by penalizing elastic properties of the medium that depends on the
local values of the density field. There are several popular material interpolation schemes, including SIMP (Solid Isotropic Material
with Penalization) \cite{Bendsoe:1989,ZhouRozvany:1991,BendsoeSigmund:1999} and RAMP (Rational Approximation
of Material Properties) \cite{StolpeSvanberg:2001} etc.

Among the different ways to formulate the minimum compliance problem, we focus on the
phase-field approximation, which employs a smooth function, i.e., material density
function, to describe the distribution of two phases, representing the material and void inside the
design domain, and encodes the geometrical information associated with the optimal topology in the
sharp interfaces between the two phases. It penalizes the objective functional by a
functional approximating the total variation of the material density function, which, in the sharp-interface
limit, measures the perimeter of the interfaces between the material and void \cite{BlankGarckeHecht:2016}.
The relaxed functional consists of two terms controlling the width of the interfaces and the decomposition
of the pure phases. Distinct features of the approach include that it greatly facilitates relevant
theoretical analysis, e.g., well-posedness of the formulation and the convergence analysis of discrete
approximation, and it eventually provides optimal solutions at the discrete level free from undesirable mesh dependency
and checkerboard  phenomena. The use of the phase-field approach in topology optimization was first
introduced for design dependent loads \cite{BourdinChambolle:2003,BourdinChambolle:2006}, problems with
stress constraints \cite{BurgerStainko:2006} and minimum compliance case \cite{WangZhou:2004,ZhouWang:2007};
see the works \cite{TakezawaNishiwakiKitamura:2010,
DedeBordenHughes:2012,HQZ,WallinRistinmaa:2013,WallinRistinmaa:2014,Tavakoli:2014,QianHuZhu:2012,LiYang:2022,GLNS:2023} for several contributions.

Numerically, one popular strategy for topology optimization is discretize-and-then-optimize: first discretizes the density field and state variable and then optimize the resulting finite-dimensional optimization problem. The discretization often employs finite element methods (FEMs) to approximate the
state variable and the density, and the same mesh is generally used for both physical quantities. For the minimum compliance problem, the standard approach implements low-order displacement-based finite elements and element-wise density discretization (i.e., piecewise constant). The degree of freedom of the resulting discretization(s) directly determines the computational cost of the overall procedure. Finer meshes improve the accuracy of the geometric description and the approximation of mechanical behavior, but involve a larger number of unknowns for the density and displacement fields and thus increase the overall computational complexity for the optimizer and finite element solver. The approximation of the state variable is nontrivial due to the nearly piecewise constant nature of the elastic properties and geometrical variations, which tend to induce weak singularities in the solution along the
interfaces. Thus, it is highly nontrivial to specify meshes \textit{a priori} for computational efficiency.

One promising idea to tackle the computational challenge is adaptivity, which has
been explored in several recent works \cite{CostaAlves:2003,Stainko:2006,SturlerPaulino:2008,BruggiVerani:2011,
	LambeCzekanski:2018,TroyaTortorelli:2018}. Costa and Alves \cite{CostaAlves:2003} described an adaptive
algorithm combining topology optimization and a Zienkiewicz-Zhu based $h$-refinement technique, which
however, does not include coarsening void regions and filtering procedures, which hence exhibits
undesirable mesh-dependence convergence behaviour. Stainko \cite{Stainko:2006} proposed an adaptive
multilevel refinement driven by the detection of the evolving solid-void interface based on RAMP-based filtering
scheme and linear material interpolation. de Sturler et al. \cite{SturlerPaulino:2008} showed the
advantage provided by coarsening operations within the cavities of the optimal design and commented on the
need for uniform refinement in and around the full material zones, in order to avoid convergence
to sub-optimal results that may arise in early coarse grids. Bruggi and Verani \cite{BruggiVerani:2011}
presented a fully adaptive topology optimization algorithm using goal-oriented error control in order to
simultaneously ensure accurate geometrical description and compliance approximation. Lambe and Czekanski
\cite{LambeCzekanski:2018} studied adaptive mesh refinement using different analysis and design meshes,
based on Helmholtz type density filter. de Troya and Tortorelli \cite{TroyaTortorelli:2018} presented
a structural optimization framework with adaptive mesh refinement and stress constraints. All these
proposals have demonstrated impressive empirical performance under various problem settings in that the
obtained optimal designs are comparable with that by uniform meshes but at a much reduced computational
expense. However, the rigorous convergence analysis is still unavailable for all these algorithms. This is
largely due to the incorporation of \textit{ad hoc} filtering procedures in the algorithms for overcoming numerical
instabilities (e.g., checkboard), and also since the error estimators are derived in a heuristic manner, especially for geometric description
\cite{Stainko:2006,BruggiVerani:2011}. Hence, there is still a demand in developing provably convergent
adaptive topology optimization algorithms.
	
Motivated by these observations, we propose in this work a novel adaptive algorithm for
the volume-constrained minimum compliance problem using the phase field approximation, and call the resulting algorithm an adaptive phase-field method.
The proposed algorithm consists of an iterative procedure of the following form
\begin{equation*}
  \text{OPTIMIZE} \rightarrow \text{ESTIMATE} \rightarrow \text{MARK} \rightarrow \text{REFINE}.
\end{equation*}
The algorithm employs the same triangulation of simplicial elements for approximating the density	and displacement field. The four modules in the loop are defined as follows. The module OPTIMIZE solves a well-posed minimization problem based on the phase-field approximation on the underlying mesh using techniques from mathematical programming. (In contrast, in several existing adaptive techniques, this step employs \textit{ad hoc} filtering techniques to tackle numerical instabilities.) The module ESTIMATE computes two error estimators of residual type \cite{AinsworthOden:2000,Verfurth:2013}, one for the error due to the finite element approximation of the displacement field, and the other for the variational inequalities in the first-order necessary optimality system. These two estimators drive the whole adaptive procedure and are derived via a careful analysis of the necessary optimality system. Thus, the proposed algorithm provides an adaptive discretization that iteratively improves the geometrical description and the approximation of mechanical behavior. The module REFINE refines the marked elements within the triangulation (by the
module MARK) to produce a locally refined mesh. Typical meshes show concentrated refinement layer near the boundaries of the optimal design with a variable size mesh within the bulk, depending on the strain energy density.
	
Numerical simulations on six benchmark problems show that the proposed adaptive algorithm can find the layout of the optimal design in the first few iterations, and improve the accuracy of the final results at later iterations. The obtained optimal designs present topological and compliance errors that	are quantitatively comparable with that based on uniform meshes of larger scale, but at a reduced computational cost. Further, compared with existing adaptive techniques \cite{CostaAlves:2003,Stainko:2006,
SturlerPaulino:2008,BruggiVerani:2011,LambeCzekanski:2018,TroyaTortorelli:2018}, we derive the	error estimators rigorously, and prove the convergence of the proposed adaptive algorithm in the following sense: the sequence of discrete approximations generated by the algorithm contains a subsequence
convergent to a solution of the necessary optimality system; See Section \ref{sec:conv} for the	detailed proof. This result is viable only for the phase field approximation, and not for the filter-based approach, for which there seems no known convergence result, even for the case of uniform mesh refinement. The adaptive algorithm and its convergence analysis represent the main contributions of this work. To the best of our knowledge, this work presents a first  provably convergent adaptive algorithm for a topology optimization problem.
	
The rest of the paper is organized as follows. In Section \ref{sec:formulation},
we describe fundamentals of the volume-constrained minimum compliance problem, recalling both continuous and discrete formulations. Then in Section \ref{sec:adapt}, we present
the proposed adaptive algorithm, and derive suitable error estimators.
In Section \ref{sec:numer}, we present extensive numerical simulations to demonstrate the
excellent performances of the adaptive algorithm, and compare it with relevant uniform
refinement. In Section \ref{sec:conv}, we present the technical and lengthy convergence
analysis of the algorithm.

\section{Mathematical formulation}\label{sec:formulation}

\subsection{Minimum compliance problem}
Let $\Omega\subset \mathbb{R}^2$ be an open bounded domain. The boundary $\partial\Omega$ consists of
two disjoint open sets $\Gamma_D$ and $\Gamma_N$ and $\partial\Omega: = \overline{\Gamma}_D \cup \overline{\Gamma}_N$.
Loaded with a body force $\bold{f}\in\bold{L}^2(\Omega):=[L^2(\Omega)]^2$ and a traction $\bold{g}\in\bold{L}^2(\Gamma_N):=[L^2(\Gamma_N)]^2$
and clamped on the boundary $\Gamma_D$, the elastic structure occupying the domain $\Omega$ is described
by the linear elasticity theory in the stress-displacement formulation:
\begin{equation}\label{eqn:pde}
	\left\{
	\begin{array}{llll}
		\stress = \mathcal{C} \strain(\bold{u}), & \quad \mbox{in}~\Omega,\\
		-\bold{\nabla} \cdot \stress = \bold{f}, & \quad \mbox{in}~\Omega,\\
		\bold{u} = \bold{0}, & \quad \mbox{on}~\Gamma_D, \\
		\stress \cdot \bold{n} = \bold{g}, & \quad \mbox{on}~\Gamma_N,
	\end{array}
	\right.
\end{equation}
where $\bm{u}=\bold{u}(x)$ denotes the displacement field at equilibrium, the related strain field $\mathcal{\varepsilon}(\bold{u})$ is given by
\begin{equation*}
	\mathcal{\varepsilon}(\bold{u})=\tfrac{1}{2}(\bold{\nabla}\bold{u}+(\bold{\nabla}\bold{u})^{\rm T}),
\end{equation*}
$\sigma = (\sigma_{ij})\in\mathbb{R}^{2\times 2}$ is the stress field and $\bold{n}$ is the unit outward normal on the boundary
$\Gamma_N$, and the superscript $\mathrm{T}$ denotes the matrix transpose. Different from the standard linear elasticity theory,
in topology optimization, the elastic tensor $\mathcal{C}$ depends on the density field. There are several different
parameterizations, and we adopt the popular SIMP model \cite{Bendsoe:1989}, which takes the following form
\begin{equation*}
	\C = \rho^p \C_0,
\end{equation*}
where $\rho=\rho(x) $ is the density field taking two values in $\{\underline{\rho}, 1\}$ almost everywhere
in $\Omega$ with $0<\underline{\rho}\ll 1$, $p>1$ is a penalization parameter that is usually assumed to be equal to
3 \cite{BendsoeSigmund:2003}, and $\C_0$ is a symmetric fourth-order elastic tensor of an isotropic medium, given by
\begin{equation*}
	\C_0 \varepsilon = 2\mu \strain  + \lambda \mathrm{tr}(\strain) \mathcal{I}
\end{equation*}
with $\mu,\lambda>0$ being Lam\'{e} constants, $\mathrm{tr}$ being the trace operator of a second-order tensor field and
$\mathcal{I}\in\mathbb{R}^{2\times 2}$ being the identity matrix. The weak formulation of problem \eqref{eqn:pde} reads: find
$\bold{u}\in\bold{V}^D: = \{\bold{v}\in \bold{H}^1(\Omega)~|~\bold{v}|_{\Gamma_D}=\bold{0}\}$ such that
\begin{equation}\label{eqn:var}
	\int_\Omega \C \strain(\bold{u}) : \strain(\bold{v}) \dd x = \int_\Omega \bold{f} \cdot \bold{v} \dd x + \int_{\Gamma_N} \bold{g}\cdot
	\bold{v} \dd s,\quad  \forall\bold{v} \in  \bold{V}^{D}.
\end{equation}
The unique solvability of \eqref{eqn:var} is guaranteed by Korn's inequality (e.g., \cite{Braess:2007,BrennerScott:2008,Ciarlet:2002,Ciarlet:2013}) and Lax-Milgram theorem. In topology optimization, it is common to impose constraint on the volume of the design domain, i.e., the total volume of the design material is fixed. Thus the admissible set $\mathcal{A}$ for the density $\rho$ is given by
\[
\mathcal{A}:=\Big\{\rho \in L^\infty(\Omega):\rho\in\{\underline{\rho},1\}~\mbox{a.e. in}~\Omega,\int_\Om \rho \dd x = V^0\Big\},
\]
with $V^0$ being the total amount of material available for design.

The minimum compliance problem is to minimize the elastic potential energy
\begin{equation*}
	\int_{\Omega} \C \strain(\bold{u}) : \strain(\bold{u}) \dd x,
\end{equation*}
over the admissible set $\mathcal{A}$ subject to the variational formulation \eqref{eqn:var}. However, it is well-known that
this problem is ill-posed in the sense that generally there is no minimizer within the admissible set $\mathcal{A}$ \cite{BendsoeSigmund:2003}.
Numerically, the corresponding numerical scheme often suffers from instabilities like checkerboard
and mesh dependence, reminiscent of the ill-posedness of the continuous optimization problem. In practice, these
drawbacks may be overcome via suitable filtering procedures (applied to the finite-dimensional discrete
problem) \cite{DiazSigmund:1995,SigmundPetersson:1998,LazarovSigmund:2011} or penalization techniques (on the
perimeters) \cite{BourdinChambolle:2003,BourdinChambolle:2006}.

In this work, we follow the penalization strategy pioneered by \cite{BourdinChambolle:2003,BourdinChambolle:2006}.
It is well known within the community of ill-posed problems that one powerful approach to remedy ill-posedness
is regularization \cite{ItoJin:2015}. Below we regularize the objective functional with a total variation penalty,
following the pioneering works \cite{BourdinChambolle:2003,BourdinChambolle:2006}, which is consistent with the
piecewise constant nature of the design (density field). This naturally leads to the following minimization problem:
\begin{equation}\label{min_TV}
	\inf_{\rho\in \mathcal{A}\cap\mathrm{BV}(\Omega)}\left\{\J(\rho)= \int_{\Omega} \C \strain(\bold{u}) : \strain(\bold{u}) \dd x +
	\beta |\rho|_{\mathrm{TV}(\Omega)}\right\},
\end{equation}
where $\beta>0$ is the regularizing parameter, and $|\cdot|_{\rm TV}$ is the total variation seminorm (or Radon
measure of the distributional derivative), defined by
\begin{equation*}
	|\rho|_{\rm TV} = \sup_{\bm\phi\in C_c^1(\overline\Omega;\mathbb{R}^d): |\bm\phi(x)|\leq 1\mbox{ a.e. in }\Omega}\int_\Omega \rho(x) \nabla\cdot \bm\phi(x){\rm d} x.
\end{equation*}
Then the space $\mathrm{BV}(\Omega)$ is defined by
\begin{equation*}
    {\rm BV}(\Omega) = \{v\in L^1(\Omega): \ \ |v|_{\rm TV}<\infty\},
\end{equation*}
with the associated ${\rm BV}(\Omega)$ norm defined by $\|v\|_{\rm BV(\Omega)}=\|v\|_{L^1(\Omega)}+|v|_{\rm TV}$.

The next result gives the continuity of the parameter-to-state map $\rho\mapsto \bold{u}(\rho)$. Then
by a standard argument in calculus of variation, the compact embedding of the space $\mathrm{BV}(\Omega)$
into the space $L^1(\Omega)$ and the fact that $L^1(\Omega)$-convergence
implies almost everywhere convergence up to a subsequence \cite{EvansGariepy:2015}, one can show the existence
of a minimizer to the penalized problem \eqref{min_TV}. Thus, the penalized problem \eqref{min_TV} is well-posed.

\begin{lemma}\label{lem:fm_cont_TV}
If the sequence $\{\rho_n\}_{n\geq 1}\subset \mathcal{A}$ converges to $\rho\in \mathcal{A}$ in $L^1(\Omega)$ and almost everywhere,
then the sequence $\{\bold{u}(\rho_{n})\}_{n\geq1}$ converges to $\bold{u}(\rho)$ in $\bold{H}^1(\Omega)$.
\end{lemma}
\begin{proof}
Let $\bold{u}_n=\bold{u}(\rho_n)$ and $\bold{u}=\bold{u}(\rho)$. By \eqref{eqn:var}, we deduce
\[
\int_\Omega \rho_n^p \C_0 \strain( \bold{u}-\bold{u}_n) : \strain ( \bold{v} ) \dx= \int_\Omega (\rho_n^p - \rho^p)
\C_0 \strain (\bold{u}) : \strain (\bold{v})  \dx, \quad \forall\bold{v} \in \bold{V}^D.
\]
Taking $\bold{v} = \bold{u} - \bold{u}_n\in \bold{V}^D$ gives
\[
2\mu\underline{\rho}^p \|\strain(\bold{u}-\bold{u}_n)\|^2_{\bold{L}^2(\Omega)}
\leq \|(\rho_n^p - \rho^p) \C_0 \strain (\bold{u})\|_{\bold{L}^2(\Omega)} \|\strain(\bold{u}-\bold{u}_n)\|_{\bold{L}^2(\Omega)}.
\]
Since $|\rho_n^p-\rho^p|^2\leq 1$, Lebesgue's dominated convergence theorem \cite[Theorem 1.19]{EvansGariepy:2015} implies
$\|(\rho_n^p - \rho^p) \C_0 \strain (\bold{u})\|_{\bold{L}^2(\Omega)}\to 0.$
By Korn's inequality on the space $\bold{V}^D$ \cite{Braess:2007,Ciarlet:2013}, $\|\bold{u}_n - \bold{u} \|_{\bold{H}^1(\Omega)}
\leq c(\Omega,\Gamma_D) \|\strain(\bold{u}-\bold{u}_n)\|_{\bold{L}^2(\Omega)}$,  the desired convergence holds for the sequence $\{\bold{u}(\rho_n)\}_{n\geq 1}$.
\end{proof}

\subsection{Phase field approximation}
The numerical treatment of the total variation penalty is highly nontrivial, since the density
$\rho$ is piecewise constant almost everywhere in the domain $\Omega$. In practice, it is often more convenient
to relax the formulation via the phase field approximation. Specifically we further take
a phase-field approximation to the functional $\J$, where the total variation seminorm $|\rho|_{\rm TV}$ in \eqref{min_TV}
is approximated by a Modica--Mortola type functional \cite{ModicaMortola:1977,ModicaMortola:1977b}
\[
\mathcal{F}_\gamma(\rho):=
\left\{\begin{array}{ll}
	\dfrac{\gamma}{2} \|\nabla \rho\|^2_{L^2(\Om)}+\dfrac{1}{\gamma}\int_{\Omega}W(\rho)\dx & \mbox{if}~\rho\in H^1(\Om),\\
	+\infty & \mbox{otherwise},
\end{array}
\right.
\]
with $\gamma>0$ being a small number, controlling the width of the interface between the void and the material,
and $W: \mathbb{R}\rightarrow\mathbb{R}$ given by a double-well potential
\begin{equation*}
	W(s)=\tfrac{1}{4}(s-\underline{\rho})^2(s-1)^2.
\end{equation*}
The functional $\mathcal{F}_\gamma$ was first used by Cahn and Hilliard \cite{CahnHilliard:1958} to study the mixture of two immiscible
fluids in the phase field theory. It is well-known that the sequence $\{\mathcal{F}_{\gamma}\}_{\gamma>0}$ $\Gamma$-converges to $c_W|\cdot
|_{\rm TV}$ in $L^1(\Omega)$ with the constant $c_W=\int_{\underline{\rho}}^{1}\sqrt{2W(s)}\dd s$ \cite{Modica:1987}. The relaxed
optimization problem is then given by
\begin{subequations}
	\begin{align}
		&  \inf_{\rho\in \Atilde} \left\{\J_{\gamma}(\rho) = \int_{\Omega} \rho^p\C_0 \strain(\bold{u}) : \strain(\bold{u}) \dd x + \widetilde{\beta}\mathcal{F}_\gamma(\rho)\right\},\label{eqn:min-phase}\\
		&\mbox{subject to } {\bold{u}\in \bold{V}^{D}:}~\int_\Omega \C \strain(\bold{u}) : \strain(\bold{v}) \dd x = \int_\Omega \bold{f} \cdot \bold{v} \dd x + \int_{\Gamma_N} \bold{g}\cdot \bold{v} \dd s\quad  \forall\bold{v} \in  \bold{V}^{D},\label{eqn:var-form}
	\end{align}
\end{subequations}
where 
$\widetilde{\beta}=\beta/c_W$ and the relaxed admissible set
$\widetilde{\mathcal{A}}$ is defined as
\begin{equation*}
	\widetilde{\mathcal{A}}:=\left\{{\rho \in H^1(\Omega)}: \rho \in [\underline{\rho},1]~\mbox{ a.e. in}~\Omega,\int_\Om \rho \dd x = V^0\right\}.
\end{equation*}

By repeating the argument for Lemma \ref{lem:fm_cont_TV}, we have the following continuity result with $\Atilde$
in place of $\mathcal{A}$.
\begin{lemma}\label{lem:fm_cont_phase}
	If $\{\rho_n\}_{n\geq 1}\subset \Atilde$ converges to $\rho\in \Atilde$ in $L^1(\Omega)$ and almost everywhere,
	then the sequence $\{\bold{u}(\rho_{n})\}_{n\geq1}$ converges to $\bold{u}(\rho)$ in $\bold{H}^1(\Omega)$.
\end{lemma}

For the relaxed optimization problem \eqref{eqn:min-phase}--\eqref{eqn:var-form}, we have the following existence result.
\begin{theorem}\label{thm:min_phase}
	For each fixed $\gamma>0$, there exists at least one minimizer to problem \eqref{eqn:min-phase}--\eqref{eqn:var-form}.
\end{theorem}
\begin{proof}
Since the functional $\J_\gamma$ is nonnegative, there exists a minimizing sequence $\{\rho_n\}_{n\geq1}\subset\widetilde{\mathcal{A}}$
such that $\J_\gamma(\rho_n)\to m_\gamma:=\inf_{\rho\in\widetilde{\mathcal{A}}}\J_\gamma(\rho)$. If $m_\gamma=\infty$,
the statement holds trivially. Now we assume $m_\gamma<\infty$. This implies $\sup_n\|\nabla\rho_n\|_{L^2(\Om)}<\infty$,
which, along with the bound $\underline{\rho}\leq \rho_n\leq 1$ and the fact the domain $\Om$ is bounded, yields $\|\rho_n\|_{H^1(\Om)}\leq c$. {The uniform boundedness in $H^1(\Omega)$, Banach-Alaoglu theorem, Sobolev (compact) embedding $H^1(\Omega)\hookrightarrow L^1(\Omega)$ and closedness and convexity of the set $\widetilde{\mathcal{A}}$ imply that} there exists a subsequence, again labeled by $\{\rho_n\}_{n\geq 1}$,
and some $\rho^\ast\in\widetilde{\mathcal{A}}$ such that
\[ \rho_{n}\rightharpoonup\rho^\ast\quad\mbox{in}~H^1(\Om),\quad \rho_n\to\rho^\ast\quad\mbox{in}~L^1(\Om),\quad\rho_n\to\rho^\ast\quad\mbox{a.e. in}~\Om,
\]
{where the last relation follows from \cite[Theorem 1.21, p. 29]{EvansGariepy:2015}.}
Since $W(s)\in C^2[c_0,c_1]$, the sequence $\{W(\rho_n)\}$ is uniformly bounded and converges to $W(\rho^\ast)$ almost everywhere in
$\Om$. Then by Lebesgue's dominated convergence theorem \cite[Theorem 1.19]{EvansGariepy:2015}, we have $\int_\Om W(\rho_n)\dx\to\int_\Om
W(\rho^\ast)\dx$. By Lemma \ref{lem:fm_cont_phase} and the weak lower semi-continuity of the $H^1(\Omega)$ norm, we further obtain
\[
\J_\gamma(\rho^\ast)\leq \liminf_{n\to\infty}\J_{\gamma}(\rho_n)\leq \lim_{n\to\infty}\J_{\gamma}(\rho_n)=m_\gamma.
\]
Upon recalling $\rho^\ast\in\widetilde{\mathcal{A}}$, this completes the proof of the theorem.
\end{proof}
	
The minimizer pair $(\rho^{\ast},\bold{u}^\ast)$ to problem \eqref{eqn:min-phase}--\eqref{eqn:var-form} satisfies the following necessary optimality system:
\begin{equation}\label{optsys}
	\left\{	\begin{aligned}
		\widetilde{\beta}\mathcal{F}_\gamma'(\rho^\ast)(\phi-\rho^\ast)- \int_{\Omega} p (\rho^\ast)^{p-1} \mathcal{C}_0 \strain (\bold{u}^\ast) :
		\strain (\bold{u}^\ast) (\phi -\rho^\ast)  \dx  \geq 0, \quad \forall\phi \in \Atilde,    \\
		\int_\Omega (\rho^\ast)^p\C_0 \strain(\bold{u}^\ast) : \strain(\bold{v}) \dd x = \int_\Omega \bold{f} \cdot \bold{v} \dd x +
		\int_{\Gamma_N} \bold{g}\cdot \bold{v} \dd s,\quad  \forall\bold{v} \in  \bold{V}^{D},
	\end{aligned} \right.
\end{equation}
with the directional derivative $\mathcal{F}_\gamma'(\rho^\ast)(\phi-\rho^\ast)$ given by
\begin{equation*}
	\mathcal{F}_\gamma'(\rho^\ast)(\phi-\rho^\ast)=\gamma\int_\Omega \bold{\nabla} \rho^\ast \cdot \bold{\nabla} (\phi - \rho^\ast) \dx
	+ \dfrac{1}{\gamma}\int_{\Omega} W'(\rho^\ast)(\phi - \rho^\ast) \dx.
\end{equation*}
{In \eqref{optsys}, the variational derivative of the term $\int_\Omega \rho^p\mathcal{C}_0\varepsilon(\boldsymbol{u}):\varepsilon(\boldsymbol{u})\dx$ can be computed directly following the definition; see e.g. \cite[Appendix A]{TakezawaNishiwakiKitamura:2010} for the details.}
The necessary optimality system will play an important role in deriving the adaptive algorithm and analyzing its convergence.
	
\subsection{Finite element approximation}
Let $\cT$ be a shape regular conforming triangulation of the closure $\overline{\Omega}$ into closed triangles. Over the triangulation $\cT$, we define a continuous piecewise linear finite element space
\begin{equation*}
	  V_\cT = \left\{v\in C(\overline{\Omega}): v|_T\in P_1(T),\ \ \forall T\in~\cT\right\},
\end{equation*}
where $P_1(T)$ is the set of linear polynomials on $T$. The space $\bold{V}_\cT^D: = (V_\cT)^2 \cap\bold{H}_D^1(\Omega)$ is used for approximating the displacement field $\bold{u}$, whereas the discrete admissible set $\widetilde{\mathcal{A}}_\mathcal{T}$ for approximating the density field $\rho$ is given by
\begin{equation*}
   	\widetilde{\mathcal{A}}_\cT:= V_\cT \cap {\widetilde{\mathcal{A}}}.
\end{equation*}
Then the discrete approximation of problem \eqref{eqn:min-phase}--\eqref{eqn:var-form} is given by
\begin{subequations}
\begin{align}
	\underset{\rho_\cT \in \Atilde_\cT}{\min}\left\{\J_{\gamma,\cT}(\rho_\cT)= \displaystyle \int_{\Omega} \rho_\cT^p\C_0
	\strain(\bold{u}_\cT) : \strain(\bold{u}_\cT) \dd x  + \widetilde{\beta}\left( \dfrac{\gamma}{2}\int_\Omega
	|\bold{\nabla}\rho_\cT|^2 + \dfrac{1}{\gamma}W(\rho_\cT)\dx  \right)  \right\},\label{dismin}\\
	\mbox{subject to }~\bold{u}_\cT\in \bold{V}_\cT^D:~\displaystyle \int_\Omega \rho_\cT^p\C_0 \strain(\bold{u}_\cT) :
	\strain(\bold{v}_\cT) \dd x  = \int_\Omega \bold{f} \cdot \bold{v}_\cT \dd x + \int_{\Gamma_N} \bold{g}\cdot
	\bold{v}_\cT \dd x,\quad  \forall\bold{v}_\cT \in \bold{V}_\cT^D.\label{dismin-state}
\end{align}
\end{subequations}
Like in the continuous case, the discrete problem \eqref{dismin}--\eqref{dismin-state} has at least one minimizer, since it is a finite-dimensional optimization problem and
the objective is continuous and coercive. Further, for any given $\rho_\cT \in \Atilde_\cT$, it follows directly from Korn's inequality and Lax-Milgram theorem that
\begin{equation}\label{stab:disvp}
	\|\bold{u}_\cT\|_{H^1(\Omega)} \leq c (\|\bold{f}\|_{L^2(\Omega)}+ \|\bold{g}\|_{L^2(\Gamma_N)}).
\end{equation}
	
The discrete optimization problem \eqref{dismin}--\eqref{dismin-state} can be solved by any standard optimization algorithms from	mathematical programming, e.g., method of moving asymptotes \cite{Svanberg:1987} and CONLIN \cite{Fleury:1989}. The computational complexity of these algorithms depends directly on the number of design parameters (for parameterizing the	density field $\rho$), and implicitly on the degree of freedom of the finite element system. Since the optimal design $\rho$ is nearly piecewise constant, and has large jumps across the interfaces between void and material, the state variable $\boldsymbol{u}$ will exhibit weak solution singularity in these regions, according to the classical elliptic regularity theory \cite{GuoBabuska:1993}, which requires a fairly refined mesh for adequate resolution. This naturally motivates the use of \textit{a posteriori} error estimator to drive an adaptive procedure.

\section{Adaptive phase--field method}\label{sec:adapt}
In this section, we propose a novel adaptive strategy for the finite element approximation \eqref{dismin}--\eqref{dismin-state}.

\subsection{Adaptive phase-field method}
Let $\cT_0$ be a shape regular conforming triangulation of the domain $\overline{\Omega}$ into closed triangles such that the boundary $\Gamma_D$ is exactly covered by the restriction of $\cT_0$ on the boundary $\partial \Omega$, and $\mathbb{T}$ be the set of all possible conforming triangulations of $\overline{\Omega}$ obtained from $\cT_0$ by the successive use of bisection. Then the set $\mathbb{T}$ is uniformly shape regular, i.e., the shape regularity of any mesh $\mathcal{T}\in\mathbb{T}$ is bounded by a constant depending only on $\cT_0$ \cite{NochettoSiebertVeeser:2009}.

To describe the new adaptive algorithm for problem \eqref{eqn:min-phase}--\eqref{eqn:var-form}, we introduce some notations. The collection of all edges (respectively all interior edges) in $\cT\in\mathbb{T}$ is denoted by $\mathcal{F}_{\cT}$ (respectively $\mathcal{F}_{\cT}^i$) and its restriction to $\Gamma_N$ by $\mathcal{F}_{\cT}^{N}$, respectively. An edge $F$ has a	fixed unit normal vector $\boldsymbol{n}_{F}$ in $\overline{\Omega}$ with $\boldsymbol{n}_{F}=\boldsymbol{n}$ on $\partial\Omega$. Over any triangulation $\cT\in\mathbb{T}$, we define a piecewise constant mesh size function $h_\cT:\overline{\Omega} \rightarrow \mathbb{R}^{+}$ by
\begin{equation}\label{meshsize_def}
	  h_\cT |_T := h_T = |T|^{1/2},\quad \forall T\in\cT.
\end{equation}
Since $\cT$ is shape regular, $h_T$ is equivalent to the diameter of any $T\in\cT$.
	
For the solution pair $(\rho_\cT^\ast,\bold{u}_\cT^\ast)$ of problem \eqref{dismin}--\eqref{dismin-state}, we define two element residuals for each	 element $T\in\cT$ and two face residuals for each face $F\in\mathcal{F}_{\cT}$ by
\[
\begin{aligned}
	R_{T,1}(\rho_\cT^\ast,\bold{u}_\cT^\ast) &= \dfrac{\widetilde{\beta}}{\gamma}W'(\rho^\ast_\cT)- p(\rho_\cT^\ast)^{p-1}\C_0 \strain(\bold{u}^\ast_\cT) : \strain(\bold{u}^\ast_\cT),\\
	\bold{R}_{T,2}(\rho_\cT^\ast,\bold{u}_\cT^\ast) &=
	\bold{f}+\bold{\nabla} \cdot (\rho_\cT^{\ast})^p \C_0 \strain(\bold{u}^\ast_\cT),\\
	J_{F,1}(\rho_\cT^\ast) & = \left\{\begin{array}{lll}
		\widetilde{\beta}\gamma[\nabla\rho_{\cT}^{\ast}\cdot\boldsymbol{n}_{F}]\quad&
		\mbox{for} ~~F\in\mathcal{F}_{\cT}^i,\\ [1ex]
		\widetilde{\beta}\gamma\nabla\rho_{\cT}^{\ast}\cdot\boldsymbol{n}\quad&
		\mbox{for} ~~F\in \mathcal{F}_\cT \setminus \mathcal{F}_\cT^i,
	\end{array}\right.\\
	\bold{J}_{F,2}(\rho_\cT^\ast,\bold{u}_\cT^\ast) & = \left\{\begin{array}{lll}
		[(\rho^\ast_\cT)^p \C_0 \strain(\bold{u}^\ast_\cT) \cdot \bold{n}_F] \quad&
		\mbox{for} ~~F\in\mathcal{F}_{\cT}^i,\\ [1ex]
		(\rho^*_\cT)^p \C_0 \strain(\bold{u}^\ast_\cT) \cdot \bold{n} -\bold{g} \quad&
		\mbox{for} ~~F\in \mathcal{F}^N_\cT,
	\end{array}\right.
\end{aligned}
\]
where $[\cdot]$ denotes the jumps across an interior face $F\in\mathcal{F}_T^i$. Then for any collection of elements $\mathcal{M}_{\cT}\subseteq\cT$, we define the following error estimator
\begin{align}
	\eta^2_{\cT}(\rho^{\ast}_{\cT},\bold{u}^{\ast}_{\cT},\mathcal{M}_{\cT})
	&:=\eta_{\cT,1}^{2}(\rho_{\cT}^{\ast},\bold{u}_{\cT}^{\ast},\mathcal{M}_{\cT})
	+\eta_{\cT,2}^{2}(\rho_{\cT}^{\ast},\bold{u}_{\cT}^{\ast},\mathcal{M}_{\cT})\nonumber\\
	&:=\sum_{T\in\mathcal{M}_{\cT}}{\left[\eta_{\cT,1}^{2}(\rho_{\cT}^{\ast},\bold{u}_{\cT}^{\ast},T)
	+\eta_{\cT,2}^{2}(\rho_{\cT}^{\ast},\bold{u}_{\cT}^{\ast},T)\right]},\label{def:estimator}
\end{align}
with the estimators $\eta_{\cT,1}^{2}(\rho_{\cT}^{\ast},\bold{u}_{\cT}^{\ast},T) $ and
$\eta_{\cT,2}^{2}(\rho_{\cT}^{\ast},\bold{u}_{\cT}^{\ast},T)$ given by
\begin{align*}
\eta_{\cT,1}^{2}(\rho_{\cT}^{\ast},\bold{u}_{\cT}^{\ast},T)
& :=h_{T}^{2}\|R_{T,1}(\rho^{\ast}_{\cT},\bold{u}^{\ast}_{\cT})\|_{L^{2}(T)}^{2}
+h_{T}\sum_{F\subset\partial T}\|J_{F,1}(\rho^{\ast}_{\cT})\|_{L^{2}(F)}^{2},\nonumber\\
\eta_{\cT,2}^{2}(\rho_{\cT}^{\ast},\bold{u}_{\cT}^{\ast},T)
&:=h_{T}^{2}\|\bold{R}_{T,2}(\rho^{\ast}_{\cT},\bold{u}^{\ast}_{\cT})\|_{\bold{L}^{2}(T)}^{2}
+h_{T}\sum_{F\subset\partial T}\|\bold{J}_{F,2}(\rho^{\ast}_{\cT},\bold{u}^{\ast}_{\cT})\|_{\bold{L}^{2}(F)}^{2}.\nonumber
\end{align*}
When $\mathcal{M}_{\cT}=\cT$, the argument $\mathcal{M}_{\cT}$ will be suppressed. The estimator $\eta_{\cT}$ consists of two error indicators $\eta_{\cT,1}$ and $\eta_{\cT,2}$. Roughly speaking, the former measures the error of approximating the first-order optimality condition of the objective functional $\J_{\gamma}$, whereas the latter quantifies the discretizaion error of the state equation \eqref{eqn:var}, similar to the standard residual estimators for elliptic PDEs \cite{AinsworthOden:2000,Verfurth:2013}.

We next motivate the two \textit{a posteriori} error estimators. Since $(\rho_\cT^\ast,\bold{u}_\cT^\ast)$ is a minimizing pair of problem \eqref{dismin}-\eqref{dismin-state}, it also satisfies the discrete optimality system
\begin{equation}\label{disoptsys}
\left\{\begin{array}{ll}
\displaystyle\widetilde{\beta}\mathcal{F}_\gamma'(\rho^\ast_\cT)(\phi_\cT-\rho^\ast_\cT)- \int_{\Omega} p (\rho_\cT^\ast)^{p-1} \mathcal{C}_0 \strain (\bold{u}_\cT^\ast) :
\strain (\bold{u}^\ast_\cT) (\phi_\cT -\rho^\ast_\cT)  \dx  \geq 0, \quad \forall \phi_\cT \in \Atilde_\cT,  \\
\displaystyle\int_\Omega (\rho_\cT^\ast)^p\C_0 \strain(\bold{u}_\cT^\ast) : \strain(\bold{v}_\cT) \dd x = \int_\Omega \bold{f} \cdot \bold{v}_\cT \dd x +
\int_{\Gamma_N} \bold{g}\cdot \bold{v}_\cT \dd s,\quad  \forall\bold{v}_\cT \in  \bold{V}^{D}_\cT
  \end{array} \right.
\end{equation}
with $\mathcal{F}_\gamma'(\rho^\ast_\cT)(\phi_\cT-\rho^\ast_\cT)$ given by $\gamma\int_\Omega \bold{\nabla} \rho^\ast_\cT \cdot \bold{\nabla} (\phi_\cT - \rho_\cT^\ast) \dx + \dfrac{1}{\gamma}\int_{\Omega} W'(\rho^\ast_\cT)(\phi_\cT - \rho_\cT^\ast) \dx$. Using the variational inequality in \eqref{disoptsys}, we may further deduce that for any $ \phi\in \Atilde$,
\begin{align*}
    &\quad\widetilde{\beta}\mathcal{F}_\gamma'(\rho^\ast_{\cT})(\phi-\rho^\ast_{\cT})
    - \int_{\Omega} p (\rho^\ast_{\cT})^{p-1} \mathcal{C}_0 \strain (\bold{u}^\ast_{\cT}) : \strain
    (\bold{u}^\ast_{\cT}) (\phi -\rho^\ast_{\cT})  \dx \nonumber\\
    &= \widetilde{\beta} \mathcal{F}_\gamma'(\rho^\ast_{k})(\phi-\phi_\cT^{I})
    - \int_{\Omega} p (\rho^\ast_{k})^{p-1} \mathcal{C}_0 \strain (\bold{u}^\ast_{\cT}) : \strain (\bold{u}^\ast_{\cT}) (\phi - \phi_\cT^I)  \dx \nonumber\\
    &\quad + \widetilde{\beta} \mathcal{F}_\gamma'(\rho^\ast_{\cT})(\phi_\cT^I-\rho^\ast_{\cT})
    - \int_{\Omega} p (\rho^\ast_{\cT})^{p-1} \mathcal{C}_0 \strain (\bold{u}^\ast_{\cT}) : \strain (\bold{u}^\ast_{\cT}) (\phi_\cT^I -\rho^\ast_{\cT})  \dx \nonumber\\
    &\geq \widetilde{\beta} \mathcal{F}_\gamma'(\rho^\ast_{\cT})(\phi-\phi_\cT^I)
    - \int_{\Omega} p (\rho^\ast_{\cT})^{p-1} \mathcal{C}_0
    \strain (\bold{u}^\ast_{\cT}) : \strain (\bold{u}^\ast_{\cT}) (\phi - \phi_\cT^I)  \dx.
\end{align*} Here $\phi_\cT^I\in \Atilde_\cT$ is an appropriate interpolant of $\phi$ (see \eqref{eqn:cn_int_def} below). Moreover by invoking the relevant error estimates in Lemma \ref{lem:err-NZ} and elementwise integration by parts (see the proof of Lemma \ref{lem:residual}), we can conclude
\[
  \left| \widetilde{\beta} \mathcal{F}_\gamma'(\rho^\ast_{\cT})(\phi-\phi_\cT^I)
          - \int_{\Omega} p (\rho^\ast_{\cT})^{p-1} \mathcal{C}_0
          \strain (\bold{u}^\ast_{\cT}) : \strain (\bold{u}^\ast_{\cT}) (\phi - \phi_\cT^I)  \dx \right| \leq c  \eta_{k,1}( \rho_{\cT}^{\ast}, \bold{u}_\cT^{\ast}) \|\bold{\nabla}\phi\|_{L^{2}(\Omega)},\quad \forall \phi\in \Atilde,
\]
{where $\cT$ has been suppressed from $\eta_{k,1}( \rho_{\cT}^{\ast}, \bold{u}_\cT^{\ast},\mathcal{T}) $ as mentioned previously.}
By the second equation of \eqref{disoptsys}, a similar argument yields
\[
   \left|\int_\Omega (\rho^\ast_\cT)^p\C_0 \strain(\bold{u}^\ast_{\cT}) : \strain(\bold{v}) \dd x - \int_\Omega \bold{f} \cdot \bold{v} \dx - \int_{\Gamma_N} \bold{g}\cdot \bold{v} \dd s\right| \leq c\eta_{k,2}(\rho_\cT^\ast,\bold{u}_\cT^\ast)\|\bold{\nabla}\bold{v}\|_{\bold{L}^2(\Omega)} , \quad\forall\bold{v}\in \bold{V}^D.
\]
Thus two computable quantities, i.e., $\eta_{\cT,1}(\rho_\cT^\ast,\bold{u}_\cT^\ast)$ and $\eta_{\cT,2}(\rho_\cT^\ast,\bold{u}_\cT^\ast)$, are available to bound the two residuals associated with the optimality system \eqref{optsys}. More importantly, the above derivation provides a crucial strategy to the convergence analysis in Section \ref{sec:conv} even though $\eta_{\cT,1}(\rho_\cT^\ast,\bold{u}_\cT^\ast)$ and $\eta_{\cT,2}(\rho_\cT^\ast,\bold{u}_\cT^\ast)$ are not reliable (upper) bounds of $H^1(\Omega)$-errors for the minimizer $\rho^\ast$ and the associated  displacement field $\bold{u}^\ast$ in any classical sense  \cite{AinsworthOden:2000, Verfurth:2013}.
	
Now we can present an adaptive algorithm for problem \eqref{eqn:min-phase}--\eqref{eqn:var-form}. All dependence on a triangulation $\cT_k$ is indicated by the iteration number $k$ in the subscript. In the module \texttt{MARK}, we adopt a separate marking by a comparison of the two estimators $\eta_{k,1}{(\rho_k^\ast,\bm{u}^\ast_k)}$ and $\eta_{k,2}{(\rho_k^\ast,\bm{u}^\ast_k)}$. The marking is performed for the greater one by a bulk criterion or D\"{o}rfler's strategy, and the adaptive algorithm is driven by the dominant term in the estimator. {Note that in practice, the two estimators $\eta_{k,1}$ and $\eta_{k,2}$ are of rather different magnitude; see the illustrations in Section \ref{sec:numer}. Thus, one may adopt suitable scaling when implementing the \texttt{MARK} module, which however is not pursued in this work.}
	

\begin{algorithm}
	\caption{AFEM for the minimal compliance problem.}\label{alg_afem_topopt}
        \LinesNumbered
        \KwIn{Specify an initial mesh $\cT_{0}$, choose parameters $\theta_1,\theta_2 \in (0,1]$ and set the maximum number $K$ of refinements.}
        \KwOut{$(\rho_k^*,\bold{u}_k^*)$}
		\For {$k=0:K-1$}{
                {(\texttt{SOLVE})} Solve \eqref{dismin}-\eqref{dismin-state} over $\cT_{k}$ for the minimizer $(\rho_k^{\ast},\bold{u}_{k}^{\ast})\in\Atilde_{k}\times\bold{V}^D_{k}$
                \;
                {(\texttt{ESTIMATE})} {Compute error indicators $\eta_{k,1}^{2}(\rho_{k}^{\ast},\bold{u}_{k}^{\ast})$ and
			$\eta_{k,2}^{2}(\rho_{k}^{\ast},\bold{u}_{k}^{\ast})$}
                \;
                {(\texttt{MARK})}
                \\
                \eIf{$\eta_{k,1}^2(\rho_k^\ast,\bold{u}_k^\ast)\geq\eta_{k,2}^2(\rho_k^\ast,\bold{u}_k^\ast) $}
                {  mark a subset $\mathcal{M}_{k}\subseteq\cT_k$ such that 
        		      \begin{equation}\label{eqn:marking1}
        			\eta_{k,1}^2(\rho_k^\ast,\bold{u}_k^\ast,\mathcal{M}_k) \geq \theta_1 \eta_{k,1}^2 (\rho_k^\ast,\bold{u}_k^\ast),
        		      \end{equation}
              }
    	    { mark a subset $\mathcal{M}_{k}\subseteq\cT_k$ such that
    	\begin{equation}\label{eqn:marking2}
    		\eta_{k,2}^2(\rho_k^\ast,\bold{u}_k^\ast,\mathcal{M}_k) \geq \theta_2 \eta_{k,2}^2 (\rho_k^\ast,\bold{u}_k^\ast).
    	\end{equation} }
            {(\texttt{REFINE})} Refine each element $T$ in $\mathcal{M}_{k}$ by bisection to get $\cT_{k+1}$\;
            Check the stopping criterion\;
        }
\end{algorithm}

We specify a maximum iteration number $K$ in the input of Algorithm \ref{alg_afem_topopt} for termination. An alternative choice as the stopping criterion might be a prescribed tolerance for $\eta_{k,1}^{2}+\eta_{k,2}^{2}$, or a given bound for the number of vertices of $\cT_k$. We have the following theorem on the convergence of the
numerical approximations $\{(\rho_k^*,\bold{u}_k^*)\}_{k\geq 0}$. The proof is lengthy and
technical and thus it is deferred to Section \ref{sec:conv}. {The theorem gives only the subsequence convergence of the adaptively generated sequence of approximations $\{(\rho_{k}^\ast,\bold{u}_{k}^\ast)\}_{k\geq0}$ to a solution of the optimality system \eqref{optsys}, from which the error estimators are derived. Due to the nonconvexity of the constrained optimization problem \eqref{eqn:min-phase}--\eqref{eqn:var-form}, it can have multiple minimizers, so does the optimality system \eqref{optsys}. Thus, one generally cannot expect a whole sequence convergence.}
\begin{theorem}\label{thm:conv_afem}
If the SIMP constant $p>1$ in $\C$ is an integer, then the sequence $\{(\rho_{k}^\ast,\bold{u}_{k}^\ast)\}_{k\geq0}$
generated by Algorithm \ref{alg_afem_topopt} contains a subsequence $\{(\rho_{k_j}^\ast,\bold{u}_{k_j}^\ast)\}_{j\geq0}$ convergent to a solution pair $(\rho^\ast,\bold{u}^\ast)$ of the continuous optimality system \eqref{optsys}:
\begin{equation}\label{conv_afem}
	\|\rho_{k_{j}}^\ast-\rho^\ast\|_{H^1(\Omega)},~   \|\bold{u}_{k_{j}}^\ast-\bold{u}^\ast\|_{H^1(\Omega)}\to 0 \quad \text{as}~j\to\infty.
\end{equation}
\end{theorem}

\subsection{Implementation details of the module SOLVE}

Now we present the details for solving the constrained minimization problem \eqref{eqn:min-phase}--\eqref{eqn:var-form} based on a gradient flow type algorithm. To handle the volume constraint $\int_\Omega \rho\dx = V^0$, we use the standard augmented Lagrangian method \cite{ItoKunisch:2008}. To this end, we define a Lagrangian $\mathcal{L}(\rho,\ell)$ by
\begin{equation}\label{eq:20}
\mathcal{L}(\rho,\ell)=\mathcal{J}_\gamma(\rho)+\ell G(\rho)+\frac{\alpha}{2}G(\rho)^2,
\end{equation}
where $G(\rho):=\int_\Omega \rho \dx-V^0$, $\ell$ is the associated Lagrange multiplier, and $\alpha>0$ is a penalty parameter.
The first-order optimality condition for an optimal phase-field function $\rho^*$ satisfies the following nonlinear system:
for all $\rho \in \bar{\mathcal{A}}:=\{{\rho \in H^1(\Omega)}: \rho \in [\underline{\rho},1]~\mbox{ a.e. in}~\Omega\}$
\begin{align*}
0&\leqslant\int_\Omega\frac{\partial \mathcal{L}}{\partial\rho}(\rho-\rho^*)\dx\\
&=-\int_\Omega p\left(\rho^{*}\right)^{p-1} \mathcal{C}_{0} \varepsilon\left(\boldsymbol{u}^{*}\right): \varepsilon\left(\boldsymbol{u}^{*}\right)(\rho-\rho^*)\dx+\tilde{\beta}\mathcal{F}_{\gamma}^{\prime} \left(\rho^{*}\right)(\rho-\rho^*)+\int_\Omega(\ell+\alpha G(\rho^*))(\rho-\rho^*)\dx.
\end{align*}
In the computation, we update the variables $\rho$, $\ell$, and $\alpha$ alternatingly. We denote by
$\rho_{n}$, $\ell_{n}$ and $\bm{u}_{n}$ the approximations at the $n$th iteration, respectively. We first solve for
the $n$th linear elasticity $\bm{u}_n$, and then obtain $\rho_{n+1}$ and $\ell_{n+1}$ successively.
Clearly, the crucial step is to minimize the Lagrangian $\mathcal{L}$ with respect to the density $\rho$. This step
is performed by solving the following parabolic PDE, a gradient-descent flow:
\begin{equation}\label{timePDE}
\left\{\begin{aligned}
\frac{\partial \rho}{\partial t}&=-\frac{\partial \mathcal{L}}{\partial\rho}, \quad &&{\rm in}\ \Omega\times\mathbb{R}^+,\\
\rho(0)&=\rho_0,\quad &&{\rm in}\ \Omega,\\
\frac{\partial\rho}{\partial \bm{n}}&=0, \quad &&{\rm on}\ \partial \Omega,
\end{aligned}\right.
\end{equation}
where $t>0$ is the pseudo time and $\rho_0\in L^2(\Omega)$ is the initial guess. To solve \eqref{timePDE}, we employ
an inner iteration index $m=0,1,\cdots$ such that the whole optimization algorithm has a nested inner-outer iteration
framework. {The analysis of the gradient formulation, e.g., convergence, is an important issue; see the recent work \cite{Plotnikov:2023} for in-depth discussions on related shape optimization.} The semi-implicit {in time} variational formulation of problem \eqref{timePDE} is to find $\rho^{[m+1]}_n\in H^1(\Omega)$ such that for all $ \varphi\in H^1(\Omega)$:
\begin{equation}\label{evo}
\begin{aligned}
&\tau^{-1}\big(\rho^{[m+1]}_{n}-\rho^{[m]}_{n},\varphi\big)+\tilde{\beta}\gamma(\nabla \rho^{[m+1]}_{n},\nabla\varphi)
= \Big(-p(\rho_n^{[m]})^{p-1} \mathcal{C}_{0} \varepsilon(\boldsymbol{u}_n): \varepsilon(\boldsymbol{u}_n)+\tilde{\beta}^{-1}{\gamma}W'(\rho_n^{[m]})+\ell_n+\alpha_n G(\rho_n^{[m]}),\varphi\Big),
\end{aligned}
\end{equation}
with the initial data $\rho^{[0]}_{n}$ approximated by the standard $L^2(\Omega)$ projection of $\rho_0$,
where $\tau$ is a small time step size. {Since $\tau$ is not related to the physical time, it is commonly known as pseudo time step size.} This semi-implicit scheme is typically evolved several (e.g., $M$) steps for fixed $n$. Numerically we observe that the density function $\rho^{[M]}_n$ may do not preserve the box constraint $[\underline{\rho},1]$. To remedy the issue, we perform a projection step:
\begin{equation*}
\rho^{[M]}_n\leftarrow\min\{\max\{\underline{\rho},\rho^{[M]}_n\},1\}.
\end{equation*}
Set $\rho_{n+1}=\rho_{n}^M$ after $M$ steps of \eqref{evo}. The Lagrange multiplier $\ell$ is updated using a Uzawa type scheme:
\begin{equation*}
\ell_{n+1} = \ell_{n}+{\alpha}_{n}G(\rho_{n+1}).
\end{equation*}
The value of $\ell$ increases (respectively decreases) as the current structure
volume is larger (respectively smaller) than the target volume $V^0$. The penalty
parameter is updated during iterations:
\begin{equation*}
\alpha_{n+1}=\alpha_{n}/\xi, \quad \mbox{with } \xi \leq 1.
\end{equation*}
Now we can state the whole algorithm on a fixed mesh in Algorithm \ref{Alga}. {In the algorithm, the integer $N$ denotes the number of outer iterations to solve the linear elasticity problem and update the Lagrangian multiplier while the integer $M$ is the number of pseudo time steps (for the evolution of the gradient flow) in each inner iteration.} 
Then, Algorithm \ref{alg_afem_topopt} for the adaptive phase field optimization model
can be described specifically with the module SOLVE implemented by Algorithm \ref{Alga}.

\begin{algorithm}
    \caption{Phase-field topology optimization on fixed mesh}\label{Alga}
    \LinesNumbered
    \KwIn{parameters $\ell_{0}$, $\alpha_0$, initial guess $\rho_0^{[0]}$, and integers $N$ and $M$}
    \KwOut{the phase--field function $\rho_N^{[M]}$}
    \For {$n=1:N$}{
        Solve the state problem and update the Lagrange multiplier $\ell_n$\; 
        Compute the gradient\;
        \For {$m=1:M$}{
                Semi-implicit updating step $\rho_n^{[m]}$\;
                Projection $\rho_n^{[m]}$ to $[\underline{\rho},1]$\;
        }
        Update the penalty parameter $\alpha_{n+1}$\;
    }
\end{algorithm}

\section{Numerical results and discussions} \label{sec:numer}

In this section, we present several examples in topology optimization to demonstrate the performance of Algorithm \ref{alg_afem_topopt}. All numerical simulations are performed using MATLAB {R2023a on a personal computer with a 13th Gen Intel(R) Core(TM) i7-13700 2.10 GHz CPU and 32GB memory}. The Lam\'{e} constants $\mu$ and $\lambda$ in the elasticity tensor $\mathcal{C}_0$ of the isotropic medium are taken to be
\begin{equation*}
\mu=\frac{E}{2(1+\nu)} \quad\text{and}\quad \lambda=\frac{E \nu}{(1+\nu)(1-2 \nu)},
\end{equation*}
with Young modulus $E=1$ and Poisson's ratio $\nu =0.3$. Throughout, we fix the penalization parameter $p=3$. Problem \eqref{dismin}--\eqref{dismin-state} is solved using the augmented Lagrangian method, in order to handle the volume constraint; see Algorithm \ref{Alga}. We fix $K=6$ steps to refine the mesh, and set $\theta_1 = \theta_2 = 0.95$ in the module \texttt{MARK}. The red and blue colors in the design plots below denote the material and void, respectively.

\begin{figure}[hbt!]
    \centering
    \setlength{\tabcolsep}{0pt}
        \begin{tabular}{cccc}
        \includegraphics[width=.2\textwidth]{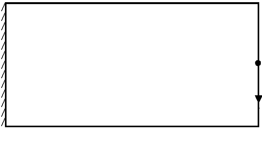}
        &\includegraphics[width=.195\textwidth]{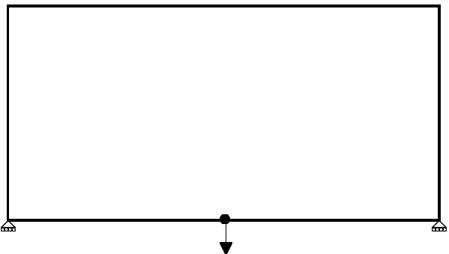}
        & \multirow{3}{*}[+11ex]{\includegraphics[width=.185\textwidth]{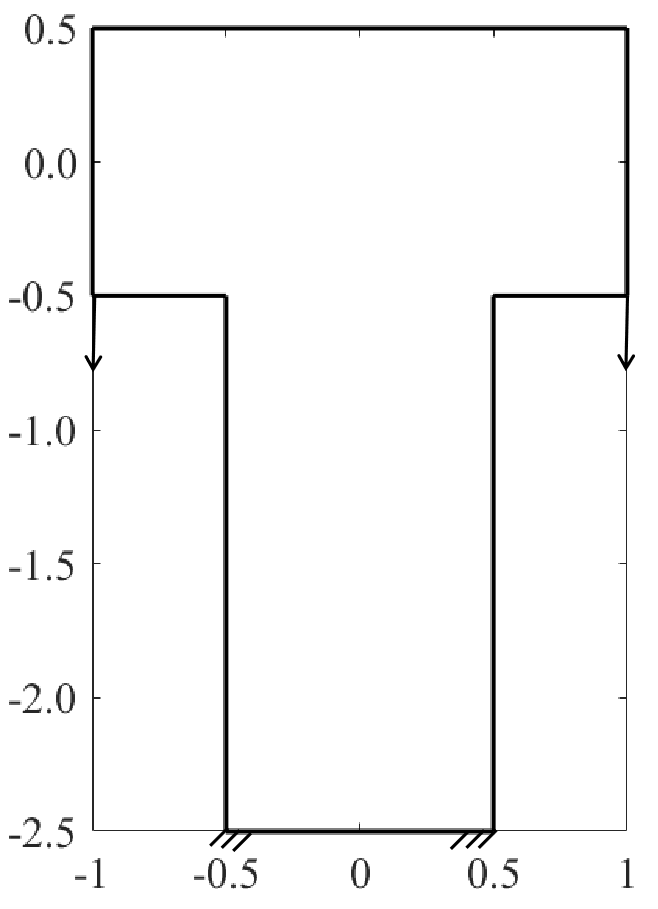}}
        & \multirow{3}{*}[+11ex]{\includegraphics[width=.185\textwidth]{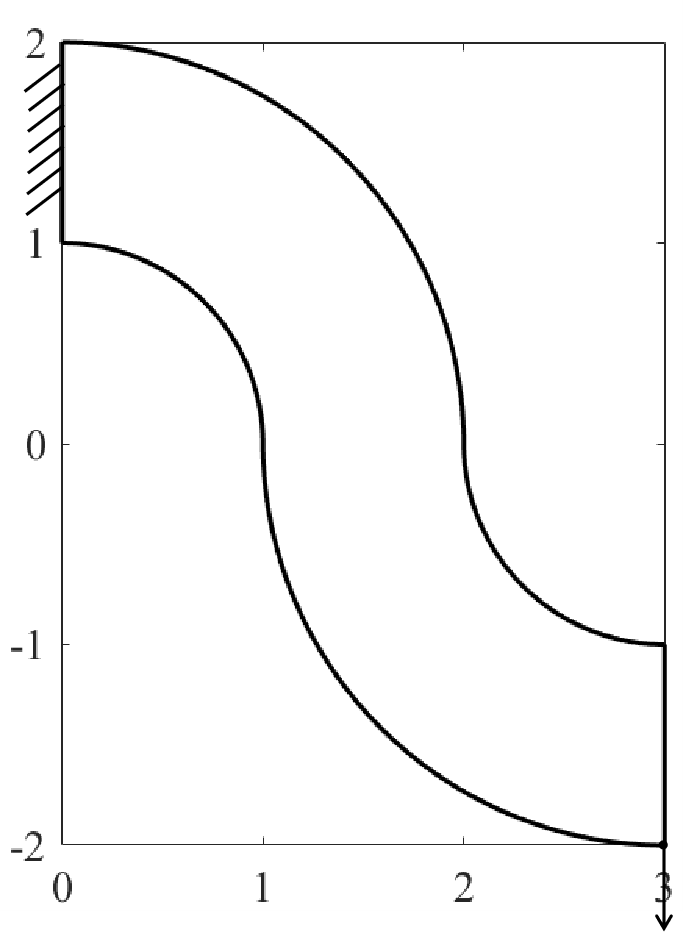}}\\
         (a) & (b) &  &  \\
        \includegraphics[width=.2\textwidth]{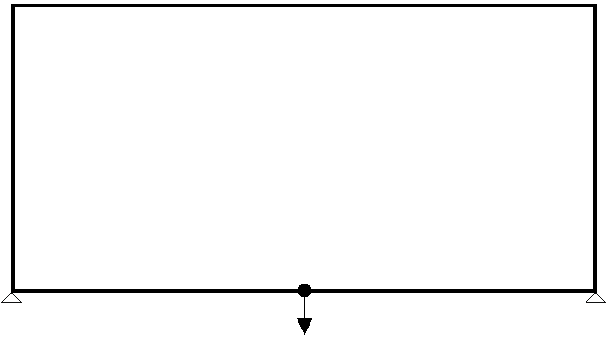}
        & \includegraphics[width=.175\textwidth]{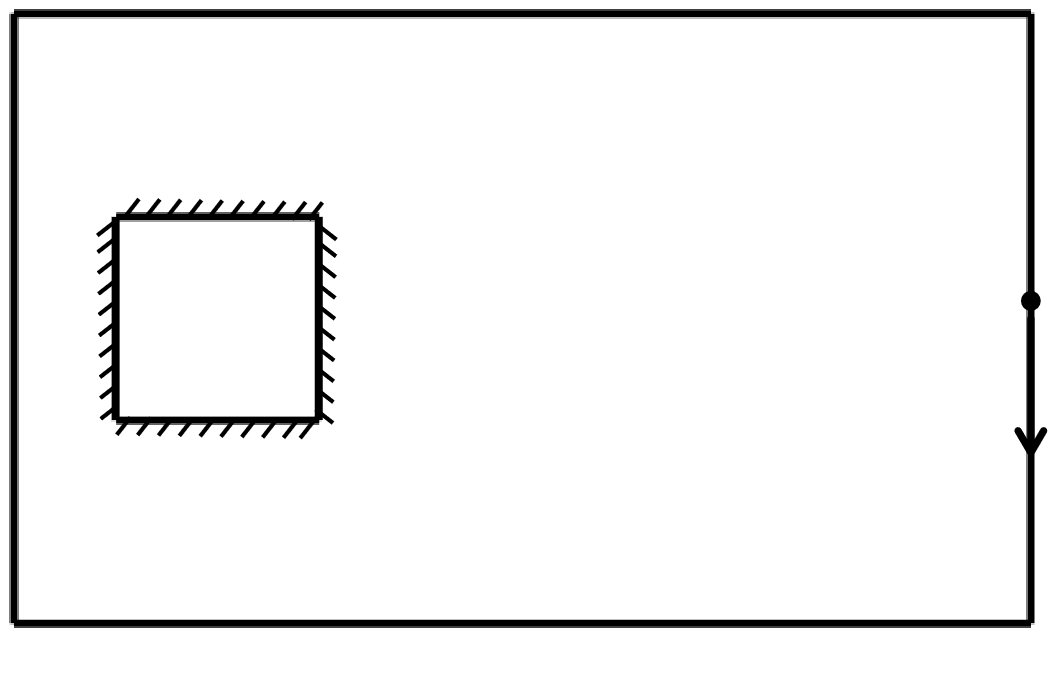}
        &  &  \\
        \multicolumn{1}{c}{(c)} & \multicolumn{1}{c}{(d)} & \multicolumn{1}{c}{(e)} & \multicolumn{1}{c}{(f)}
        \end{tabular}
    \caption{The schematic illustration of the geometry, loading and boundary conditions of the examples.\label{fig:geometry}}
\end{figure}

In the numerical experiments, we consider the following six examples. All these examples are very classical in topology optimization.
See Fig. \ref{fig:geometry} for a schematic illustration of the geometry, loading and boundary conditions of the examples.
\begin{enumerate}
  \item[(a)] The cantilever problem with the design domain $\Omega=(-1,1)\times(-0.5,0.5)$, with $V^0=1$. The displacement field on the left side is fixed at zero, and the free traction on the rest except a vertical point load $(0,-1)^\top$ is applied to the middle of the right side.
  \item[(b)] The Messerschmidt-Bolkow-Blohm (MBB) beam problem, with $V^0=0.8$. The settings are identical with (a) except now the low corners are free to move in the $x$-direction and fixed in the $y$-direction, and a unit vertical load is applied in the middle of the lower part and a traction free condition is imposed on the rest.
  \item[(c)] The cantilever problem as (a), with $V^0=0.4$, except that the low corners are fixed.
  \item[(d)] The square hole, the domain $[-2.5,2.5]\times (-1.5,1.5)\backslash[-2,-1]\times(-0.5,0.5)$, $V^0=2.2$. The inner square boundary is equipped with a zero Dirichlet boundary condition, and the free traction on the rest except a vertical point load $(0,-1)^\top$ is applied in the middle of the right side.
  \item[(e)] The T-shaped domain with height 3, width 1 at the bottom and 2 at the top, with $V^0=0.8$. The two lower corners are fixed while unit vertical load is applied at the lower corners of the horizontal branch of T.
  \item[(f)] The curved domain, comprising of four arcs and two line segments, and $V^0$ accounts for 40\% of the total volume. The left vertical edge is equipped with a zero Dirichlet boundary condition, and the free traction on the rest except a vertical point load $(0,-1)^\top$ is applied in the bottom of the right vertical edge.
\end{enumerate}

\begin{figure}[hbt!]
\centering
\begin{tabular}{cc}
\includegraphics[width=\figwid]{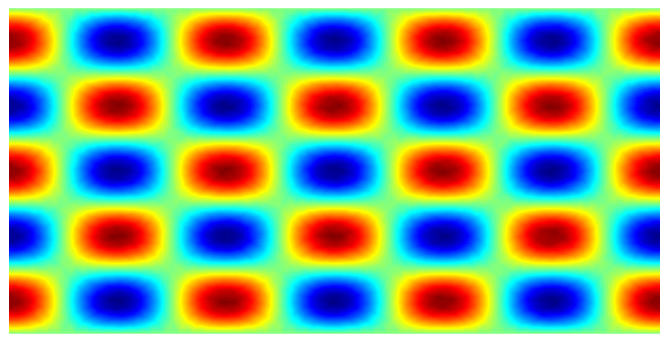} & \includegraphics[width=\figwid]{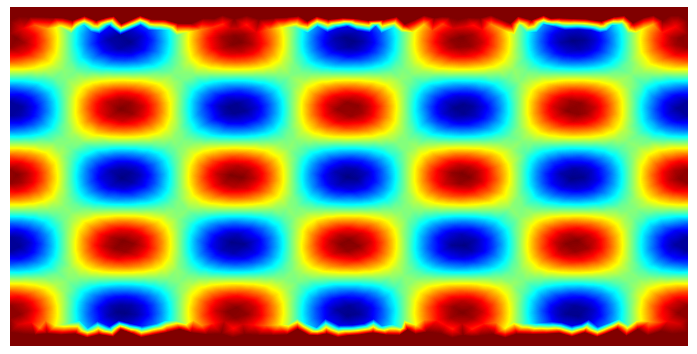}\\
\subfigure[Examples (a) and (c)]{\includegraphics[width=\figwid]{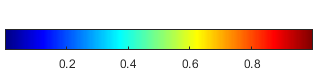}} & \subfigure[Examples (a) and (c)]{\includegraphics[width=\figwid]{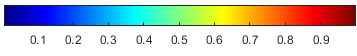}}
\end{tabular}
\caption{The initial phase--field functions $\rho_0$ for Examples (a)-(c).}\label{fig:InitialPhaseFieldFunction}
\end{figure}

The hyperparameters in Algorithm \ref{Alga} for the numerical experiments are summarized in Table \ref{tab:hyper}. In the table, $N$ denotes the number of {outer iterations} to update the elastic displacement, and $M$ denotes the number of inner iterations {or pseudo time steps for the gradient flow} to update the density field $\rho$. The scalar $\tau$ is the phase-field pseudo time step size. The parameters $\ell_0$ and $\alpha_0$ are the initial Lagrange multiplier, and initial penalty weight for the augmented Lagrangian terms. The two parameters $\beta$ and $\gamma$ appear in the phase-field model \eqref{eqn:min-phase}. The lower bound $\underline\rho$ of the density $\rho$ is set to 1e-4 except for Example (d), for which it is set to  1e-3; The factor $\xi$ is fixed at $0.99$ for all examples except for Examples (d) and (e), for which it is set to $0.95$. {Due to the nonconvexity of the constrained optimization problem \eqref{eqn:min-phase}--\eqref{eqn:var-form}, it can have multiple global minimizers, and the convergence of the sequence generated by the gradient flow formulation will depend on the initial phase-field function $\rho$.} In the experiment, the initial phase--field functions $\rho$ for Examples (a)-(c) are shown in Fig. \ref{fig:InitialPhaseFieldFunction}, and that for Examples (d)-(f) are taken to be uniform, i.e., 0.3, 0.4 and 0.6, respectively. {To the best of our knowledge, there is still no general guideline for choosing these hyper-parameters. They are determined in a trial-and-error manner in the present work.}

\begin{table}[]
\centering
\begin{threeparttable}
\caption{The hyper-parameters for Algorithm \ref{Alga} and the phase--field model \eqref{eqn:min-phase} used in the experiments. \label{tab:hyper}}
\begin{tabular}{c|cccccc}
\toprule
 Example & $(N,M)$ & $\tau$ & $\ell_0$ & $\alpha_0$ & $\beta$ & $\gamma$ \\
\midrule
 (a) & $(20,3)$ & 3.5e-2 & $0.8$ & $0.8$  & 1e-5 & 1e-2\\
 (b) & $(30,10)$& 4.5e-2 & $0.4$ & $0.2$  & 1e-5 & 1e-2\\
 (c) & $(50,5)$ & 4.0e-2 & $0.8$ & $0.2$  & 5e-5 & 1e-2\\
 (d) & $(50,5)$ & 6.5e-3 & $1.5$ & $0.65$ & 5e-4 & 1e-3\\
 (e) & $(50,5)$ & 3.0e-2 & $0.8$ & $0.03$ & 5e-5 & 1e-2\\
 (f) & $(50,5)$ & 2.5e-2 & $0.8$ & $0.1$  & 5e-5 & 2e-3\\
\bottomrule
\end{tabular}
\end{threeparttable}
\end{table}

\begin{figure}[hbt!]
\centering\setlength{\tabcolsep}{0pt}
\begin{tabular}{ccccccc}
	\includegraphics[width=\figwid]{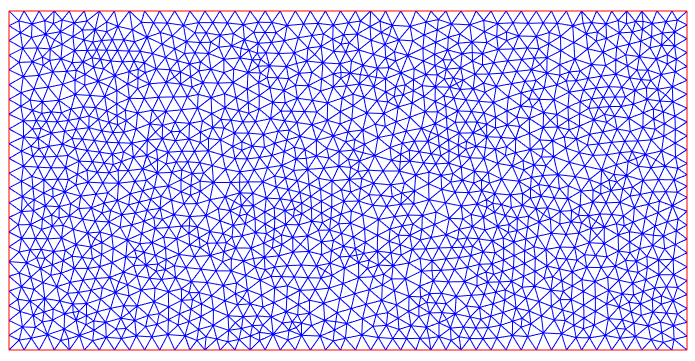}& \includegraphics[width=\figwid]{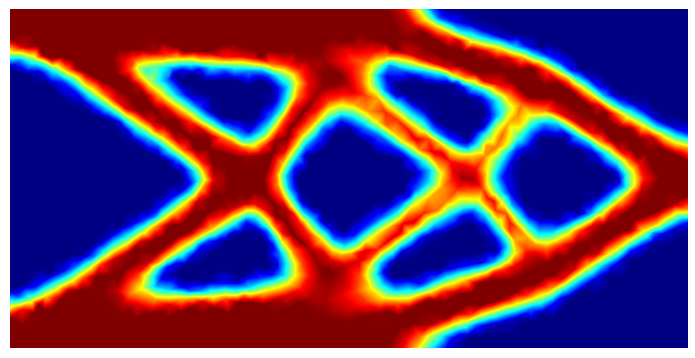}
    & \includegraphics[width=\figwid]{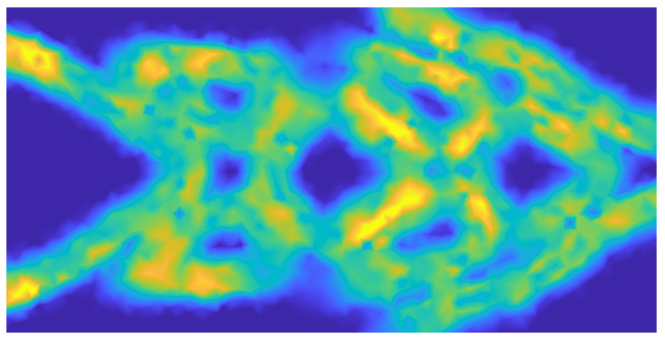} & \includegraphics[width=\figwid]{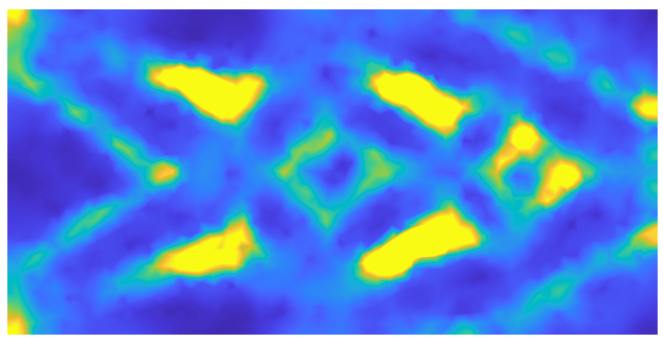}\\
   \includegraphics[width=\figwid]{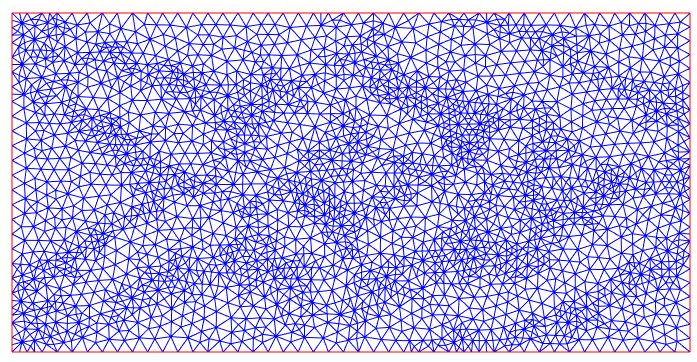}& \includegraphics[width=\figwid]{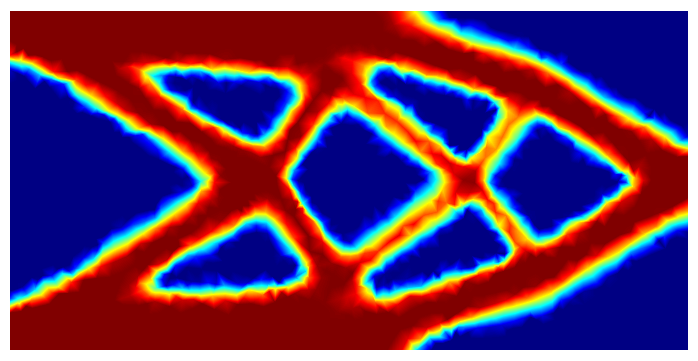}
    &\includegraphics[width=\figwid]{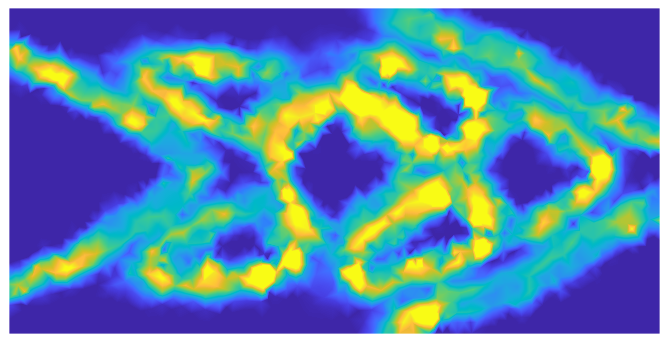}&\includegraphics[width=\figwid]{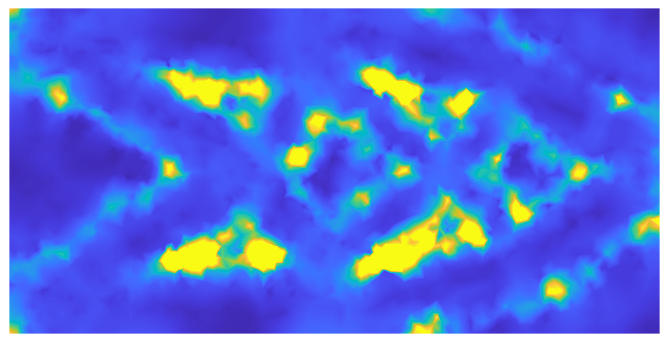}\\
	\includegraphics[width=\figwid]{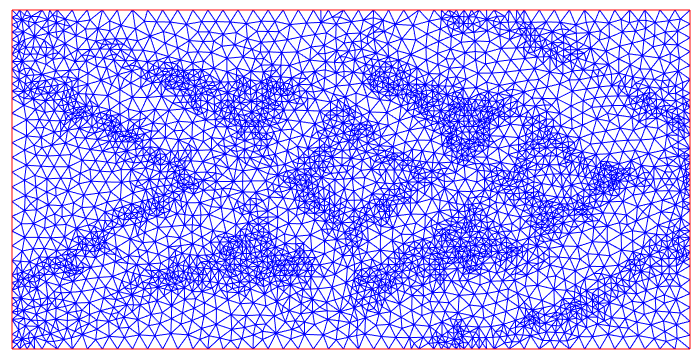}& \includegraphics[width=\figwid]{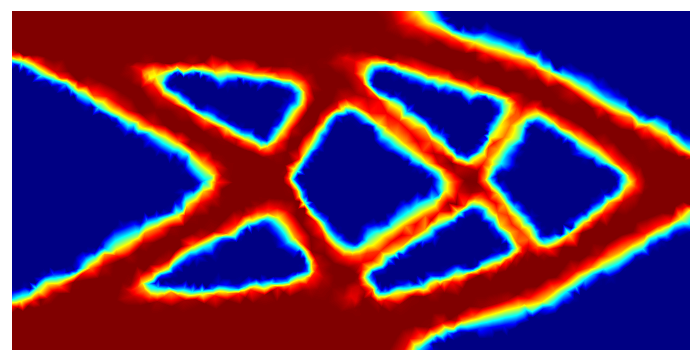}
	&\includegraphics[width=\figwid]{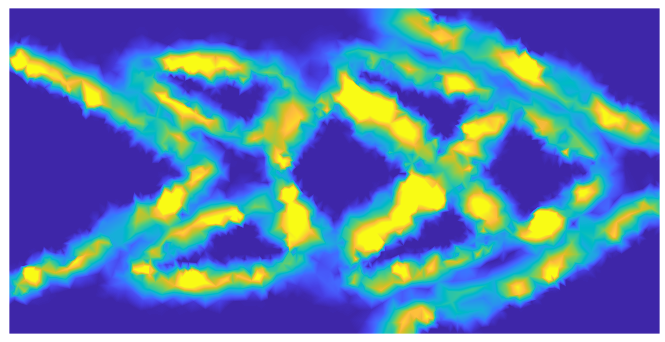}&\includegraphics[width=\figwid]{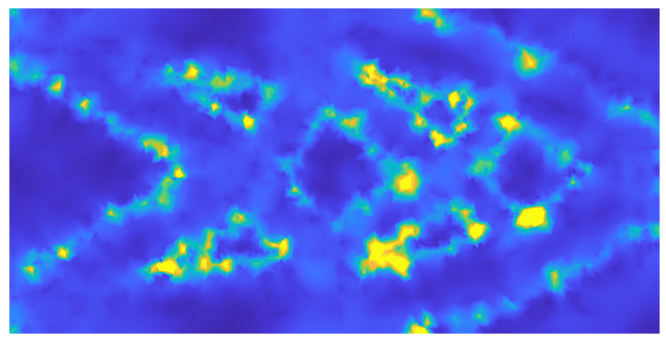}\\
    \includegraphics[width=\figwid]{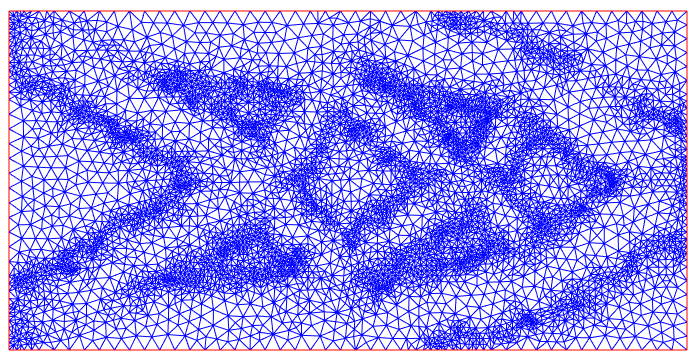}& \includegraphics[width=\figwid]{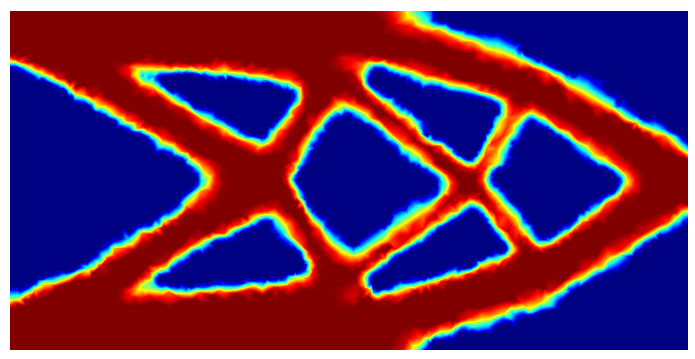}
	&\includegraphics[width=\figwid]{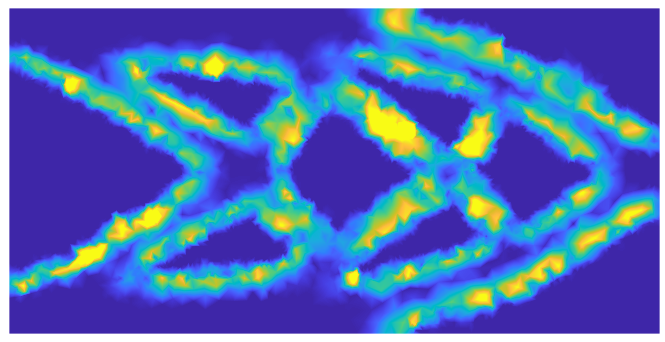}&\includegraphics[width=\figwid]{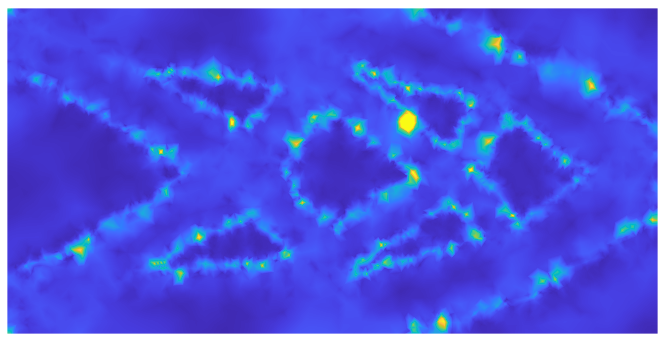}\\
    \includegraphics[width=\figwid]{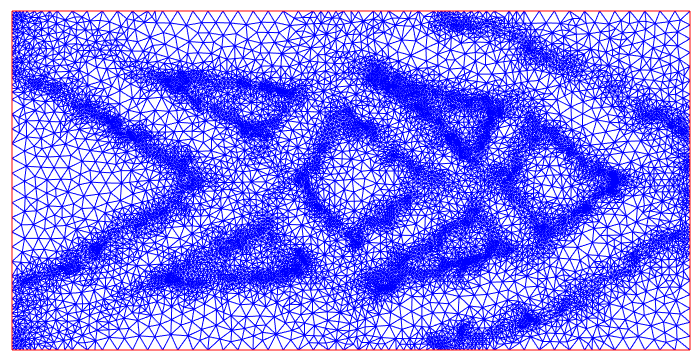}& \includegraphics[width=\figwid]{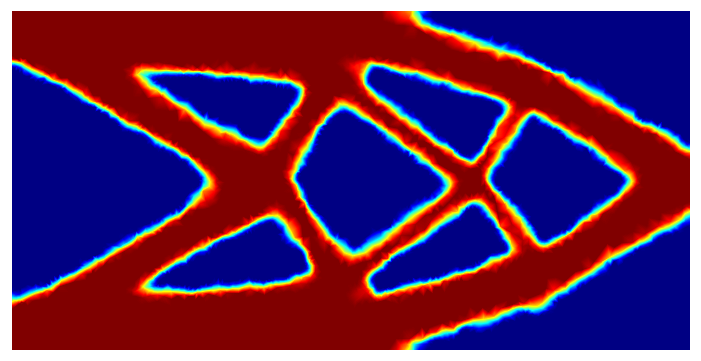}
    &\includegraphics[width=\figwid]{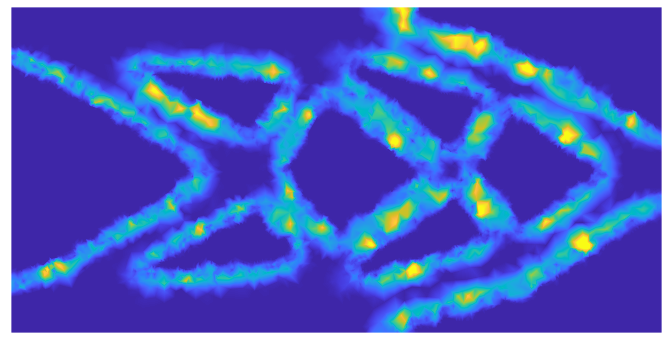}&\includegraphics[width=\figwid]{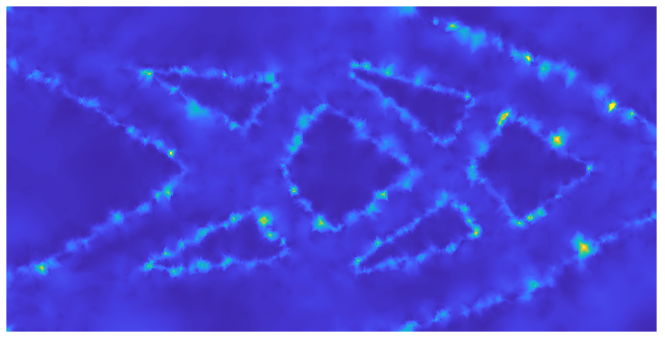}\\
    \includegraphics[width=\figwid]{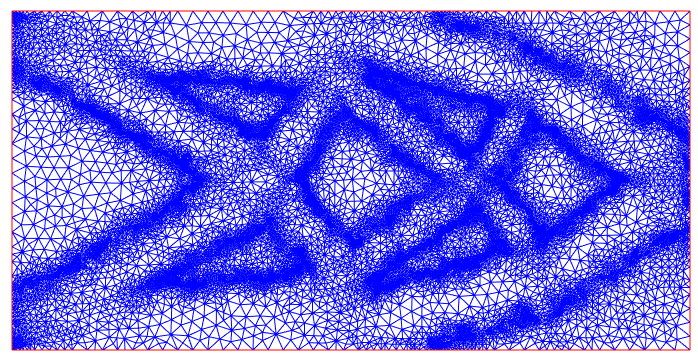}& \includegraphics[width=\figwid]{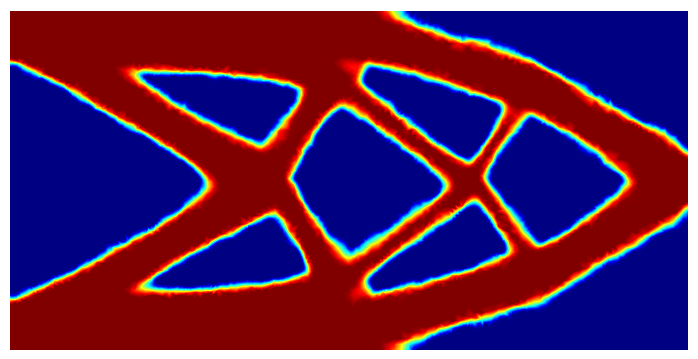}
    & \includegraphics[width=\figwid]{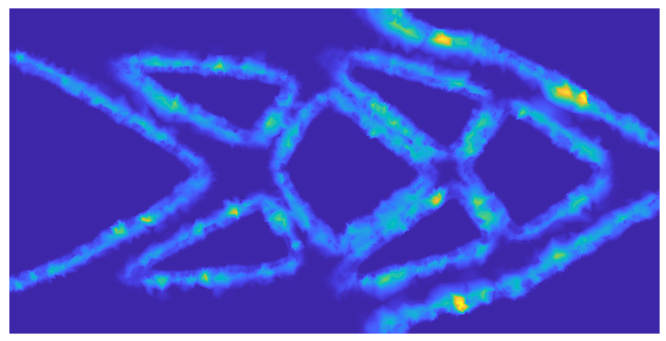}&\includegraphics[width=\figwid]{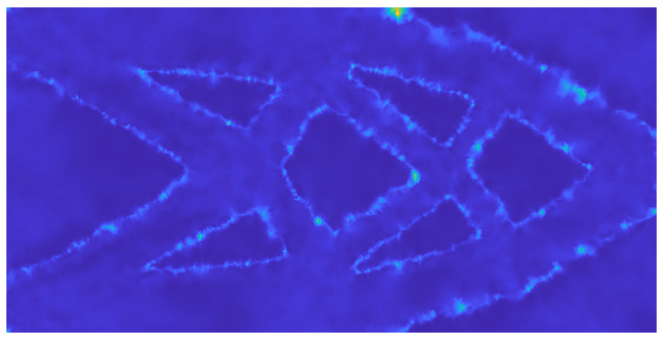}\\
    & & \includegraphics[width=\figwid]{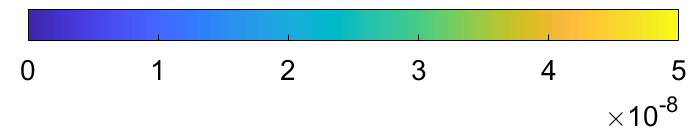}
    &\includegraphics[width=\figwid]{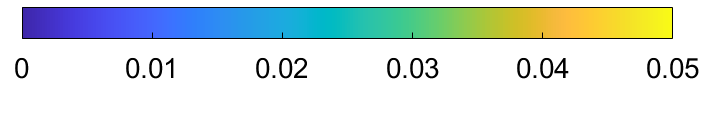}
\end{tabular}
\caption{The evolution of the mesh during the adaptive process, from $k=0$ (initial) to 5 for Example (a), with the number of vertices of each mesh being 1718, 2782, 4039, 6464, 10684 and 16687. The second, third and last columns show the optimized design $\rho_k^*$, the error indicators $\eta_{k,1}$ and $\eta_{k,2}$ respectively.} \label{fig:Cantilever_MehsAndDeigsn}
\end{figure}

\begin{figure}[htb!]
\centering \setlength{\tabcolsep}{0pt}
\begin{tabular}{ccccccc}
	\includegraphics[width=\figwid]{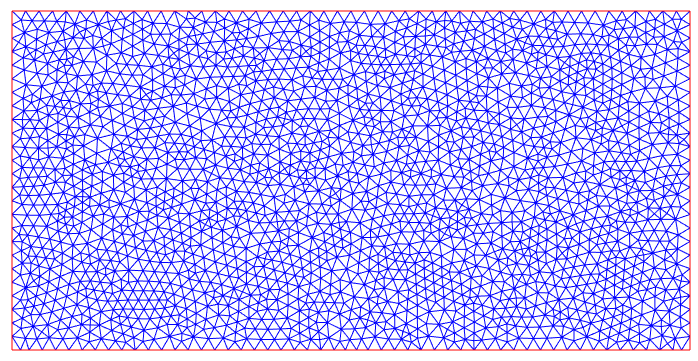} &\includegraphics[width=\figwid]{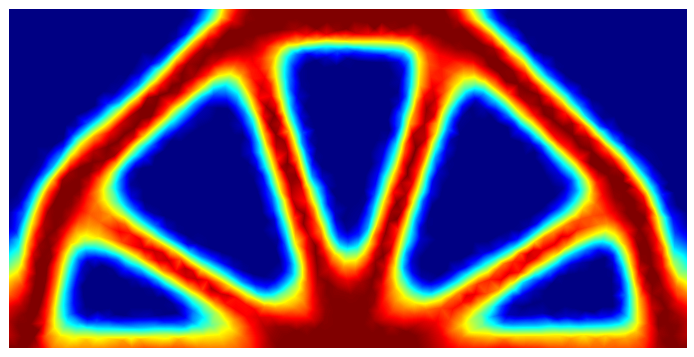}
    &\includegraphics[width=\figwid]{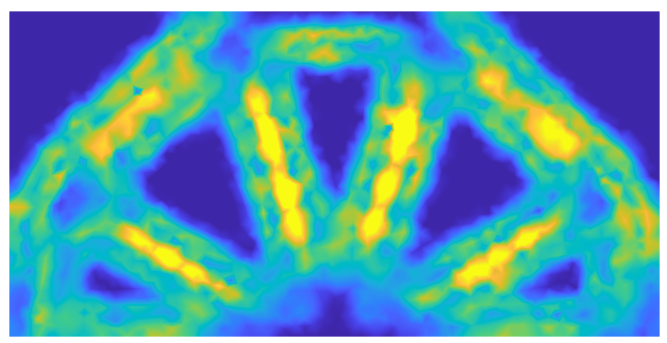}&\includegraphics[width=\figwid]{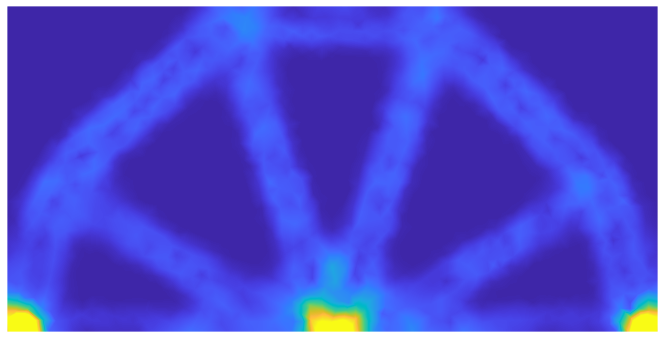}  \\
    \includegraphics[width=\figwid]{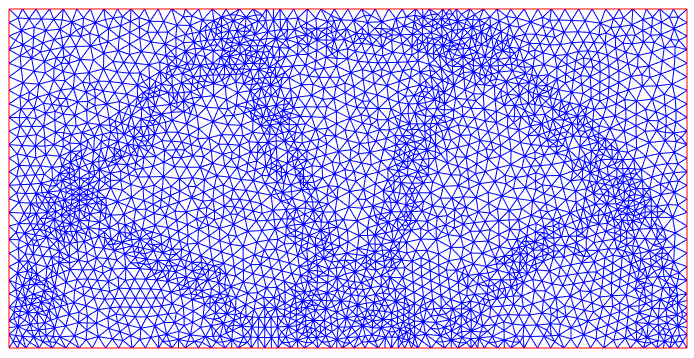} &\includegraphics[width=\figwid]{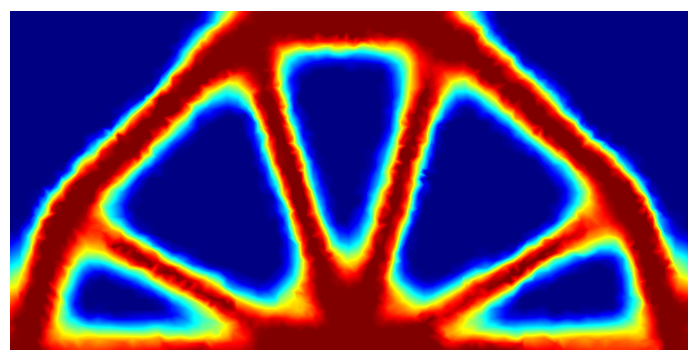}
    &\includegraphics[width=\figwid]{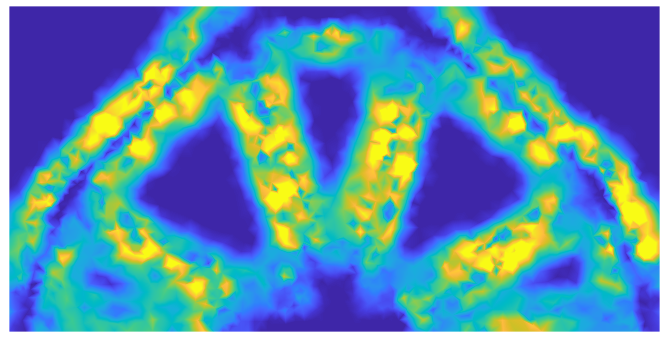}&\includegraphics[width=\figwid]{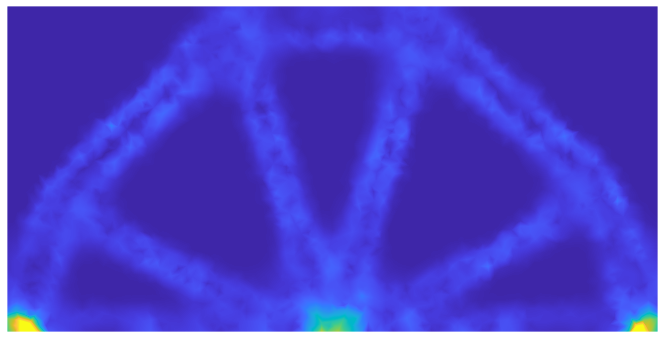}  \\
	\includegraphics[width=\figwid]{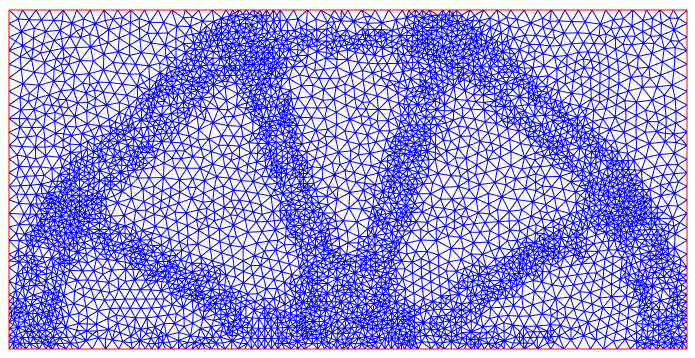}  &\includegraphics[width=\figwid]{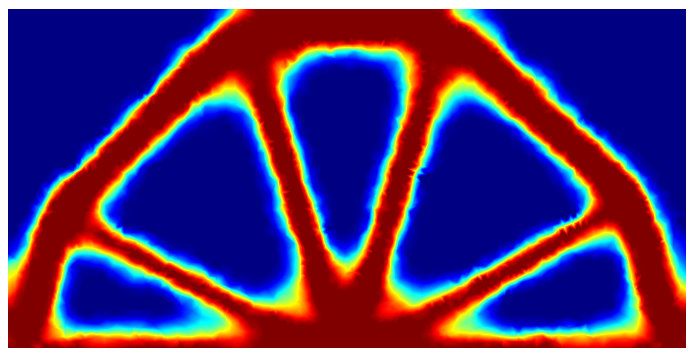}
    &\includegraphics[width=\figwid]{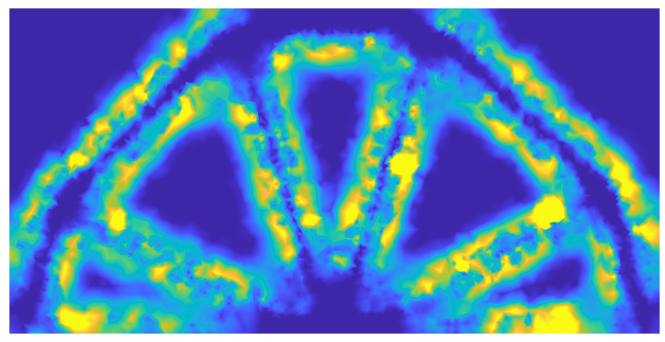}&\includegraphics[width=\figwid]{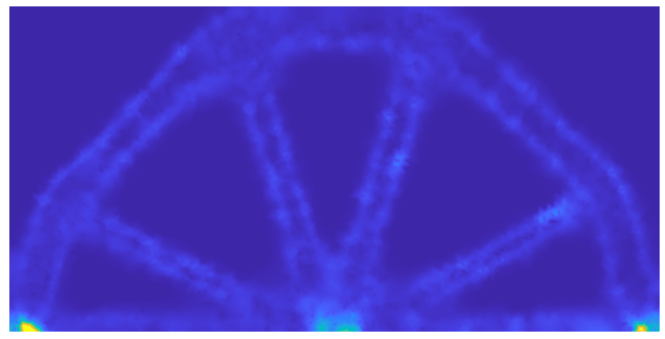}  \\
    \includegraphics[width=\figwid]{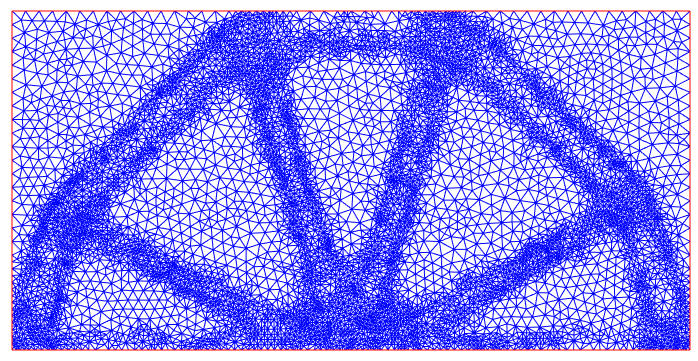} &\includegraphics[width=\figwid]{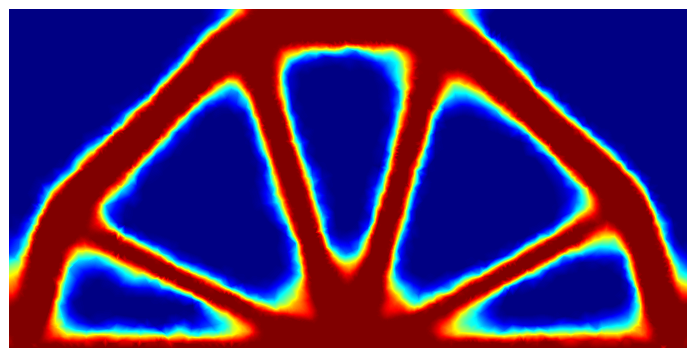}
    &\includegraphics[width=\figwid]{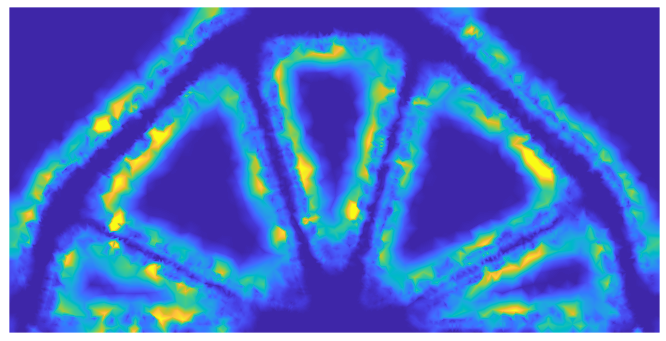}&\includegraphics[width=\figwid]{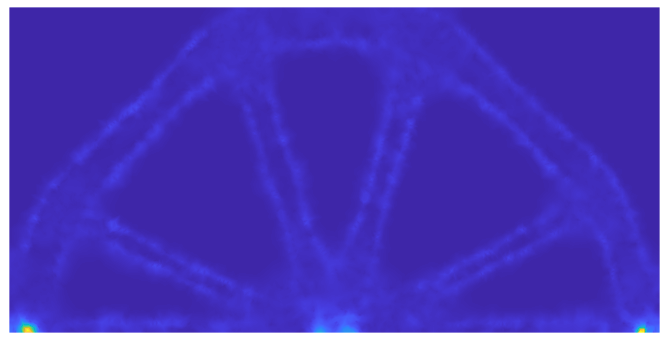}  \\
	\includegraphics[width=\figwid]{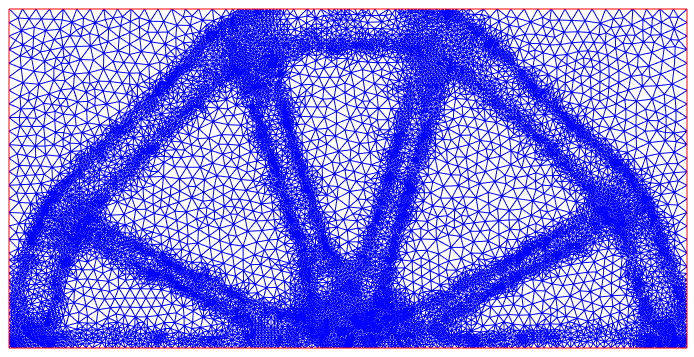}  &\includegraphics[width=\figwid]{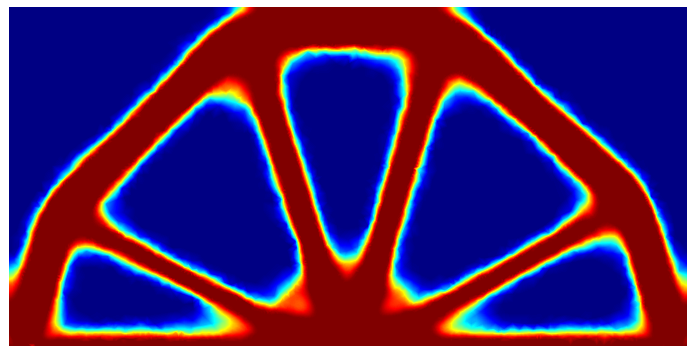}
    &\includegraphics[width=\figwid]{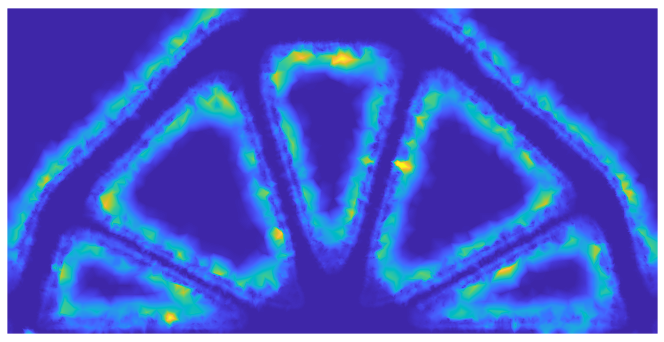}&\includegraphics[width=\figwid]{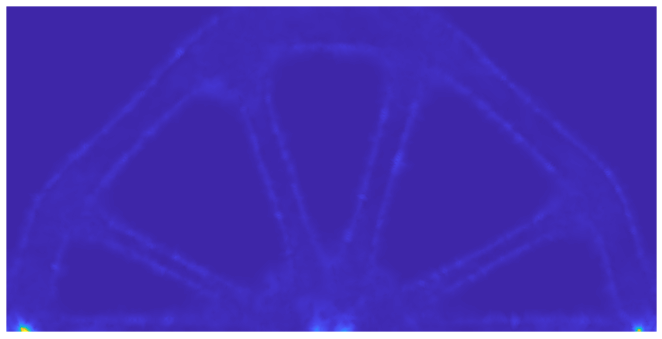}  \\
    \includegraphics[width=\figwid]{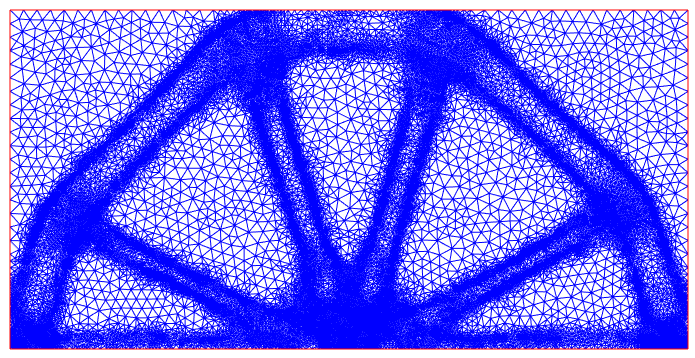} &\includegraphics[width=\figwid]{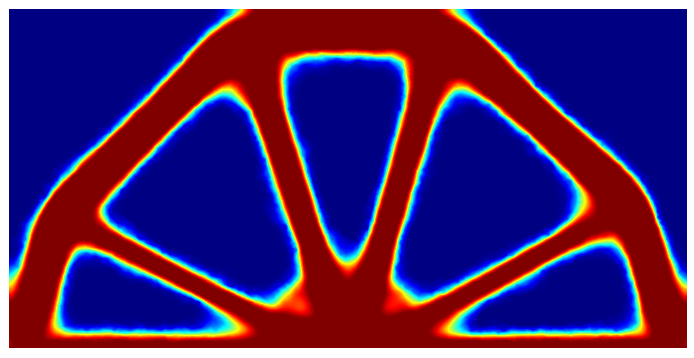}
    &\includegraphics[width=\figwid]{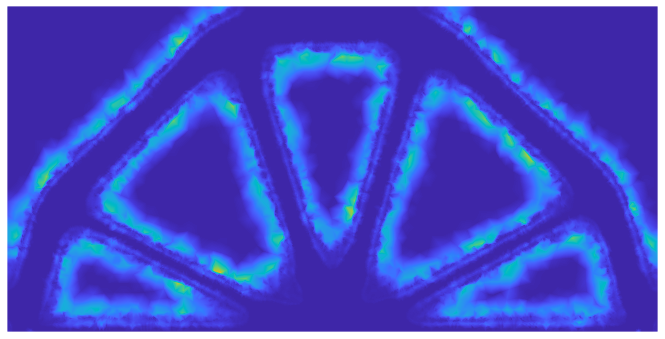}&\includegraphics[width=\figwid]{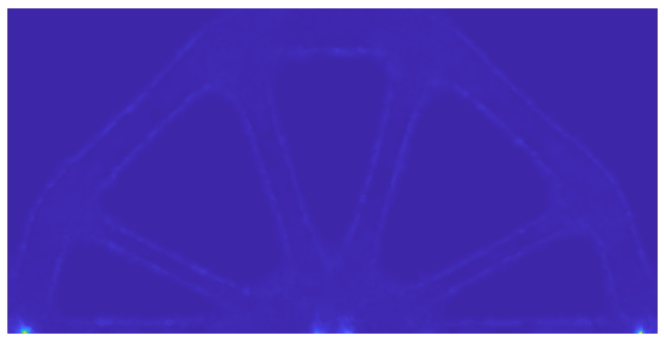}\\
    & &\includegraphics[width=\figwid]{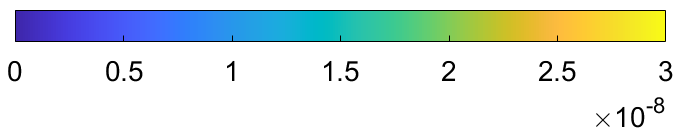} &\includegraphics[width=\figwid]{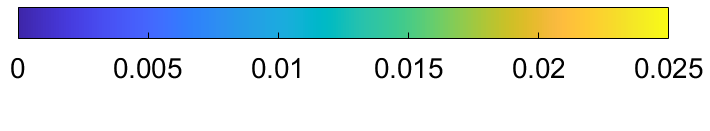}
\end{tabular}
\caption{The evolution of the mesh during the adaptive procedure, from $k=0$ (initial) to 5 for Example (b) (MBB), with the number of vertices of each mesh being 2230, 3338, 5237, 8525, 14084 and 23314. The second, third, and fourth columns show the optimized design $\rho_k^*$,
error indicators $\eta_{k,1}$ and $\eta_{k,2}$ respectively.}\label{fig:MBB2_MehsAndDeigsn}
\end{figure}

\begin{figure}[htb!]
\centering \setlength{\tabcolsep}{0pt}
\begin{tabular}{ccccccc}
\includegraphics[width=\figwid]{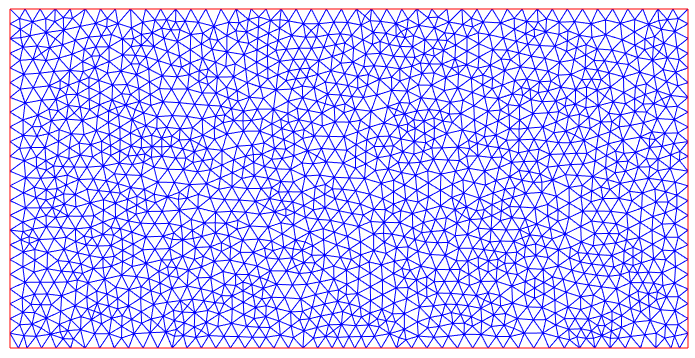} &\includegraphics[width=\figwid]{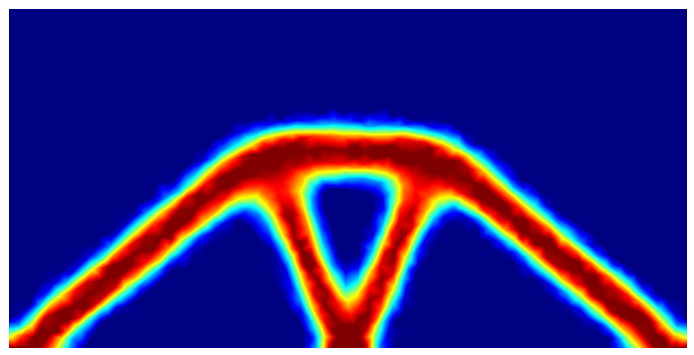}&\includegraphics[width=\figwid]{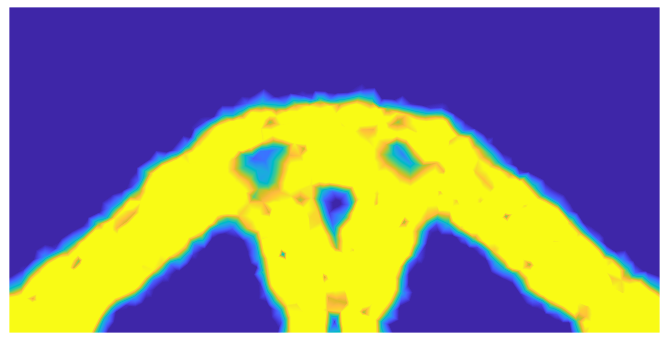} &\includegraphics[width=\figwid]{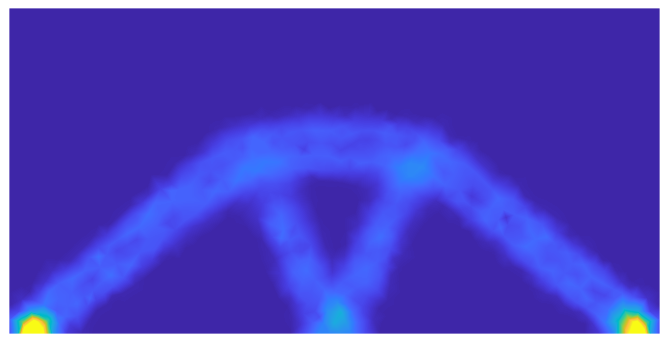}\\
\includegraphics[width=\figwid]{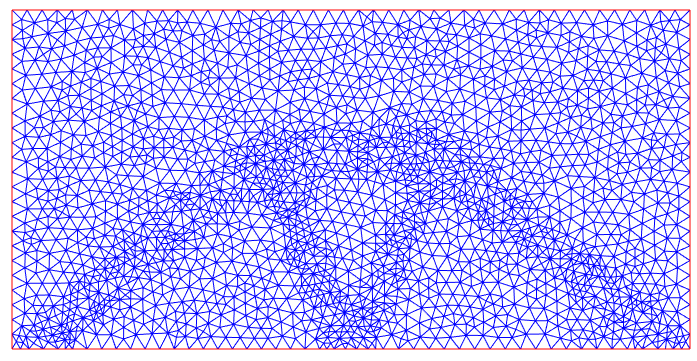} &\includegraphics[width=\figwid]{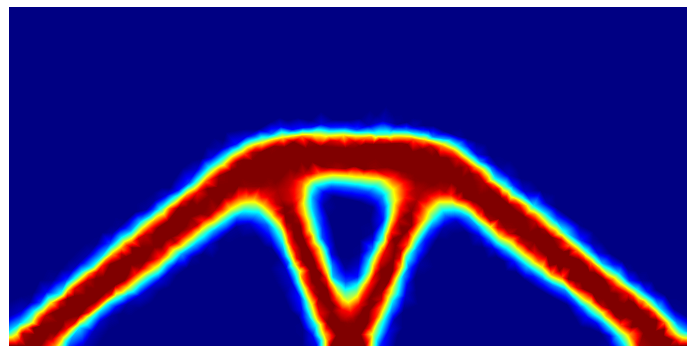}&\includegraphics[width=\figwid]{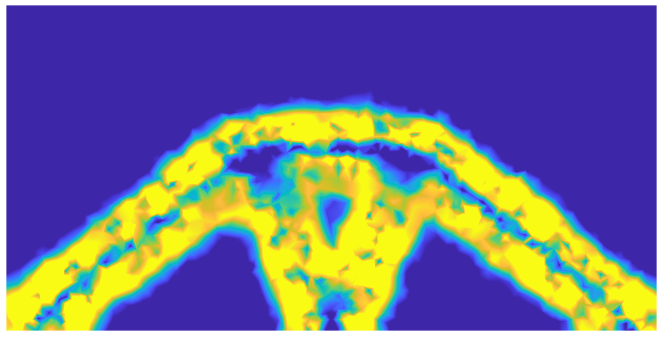} &\includegraphics[width=\figwid]{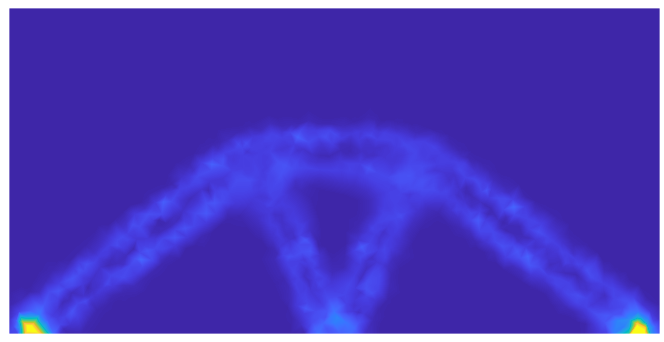}\\
\includegraphics[width=\figwid]{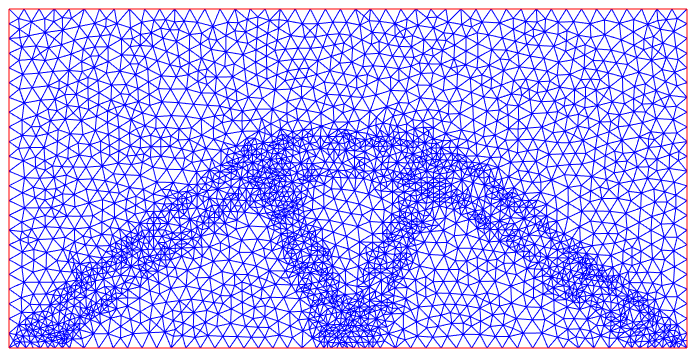} &\includegraphics[width=\figwid]{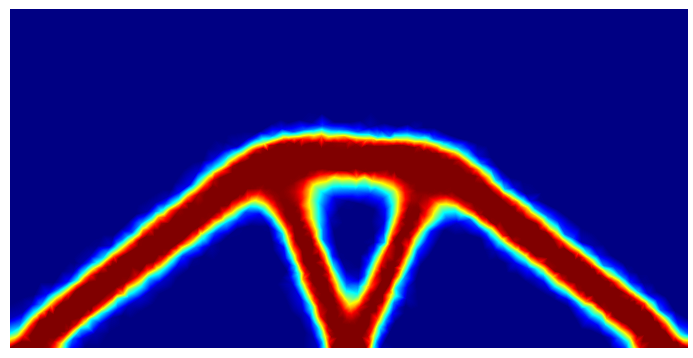}&\includegraphics[width=\figwid]{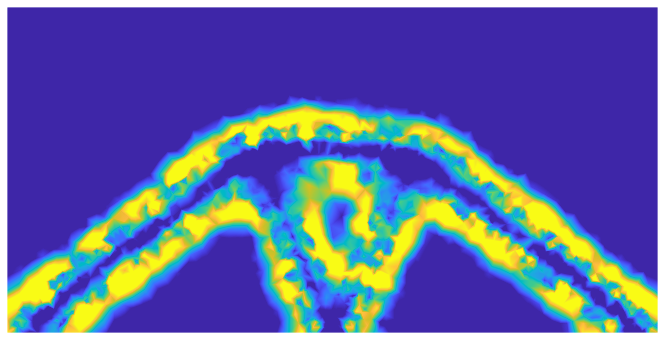} &\includegraphics[width=\figwid]{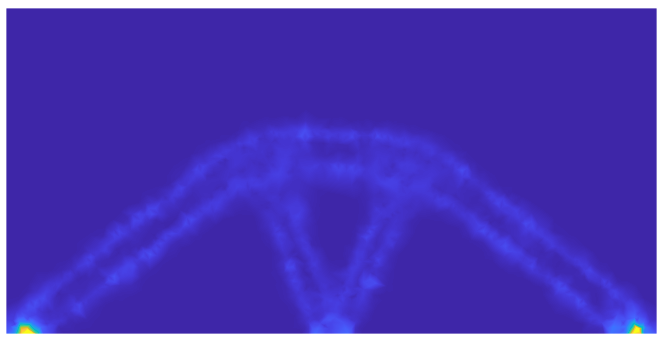}\\
\includegraphics[width=\figwid]{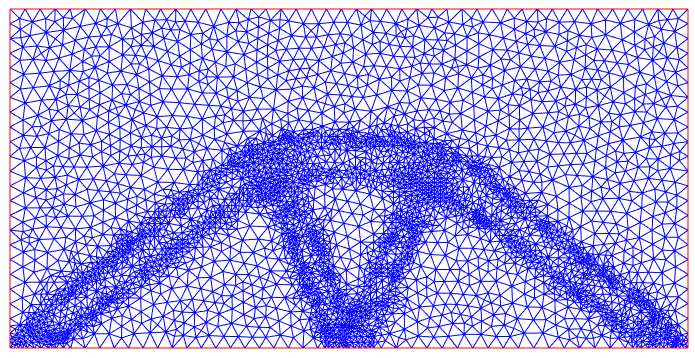} &\includegraphics[width=\figwid]{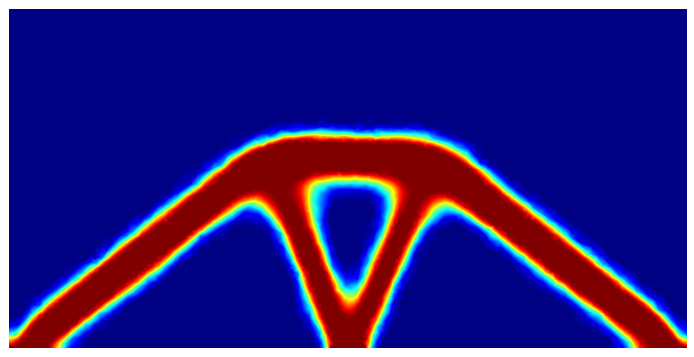}&\includegraphics[width=\figwid]{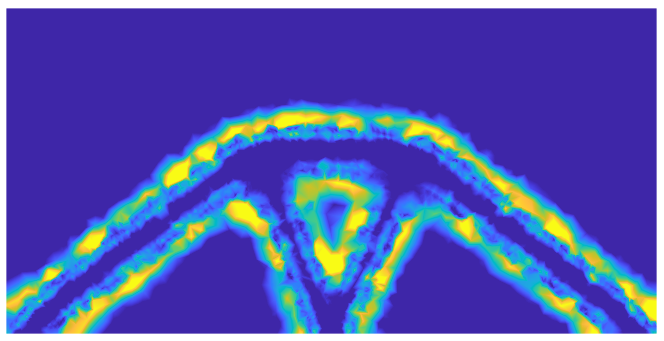} &\includegraphics[width=\figwid]{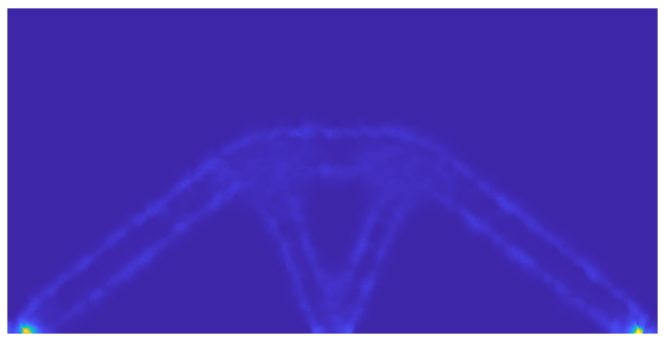}\\
\includegraphics[width=\figwid]{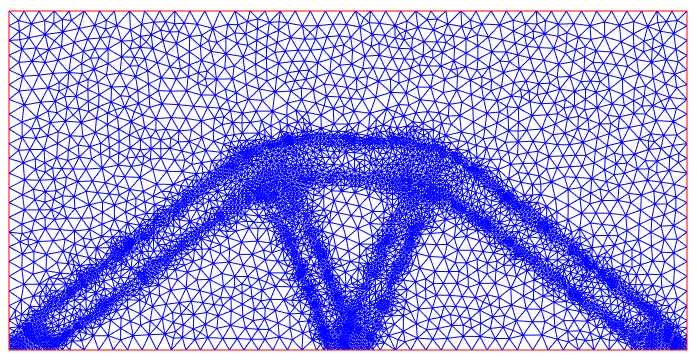} &\includegraphics[width=\figwid]{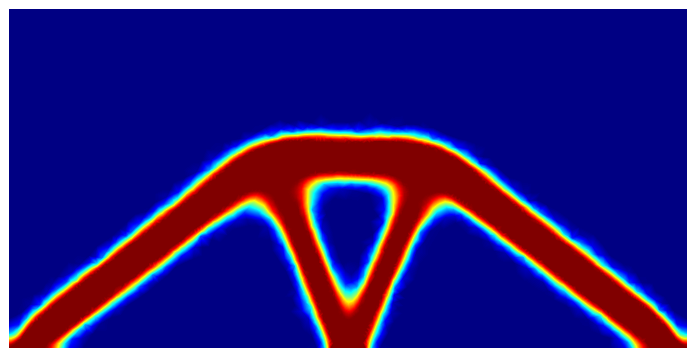}&\includegraphics[width=\figwid]{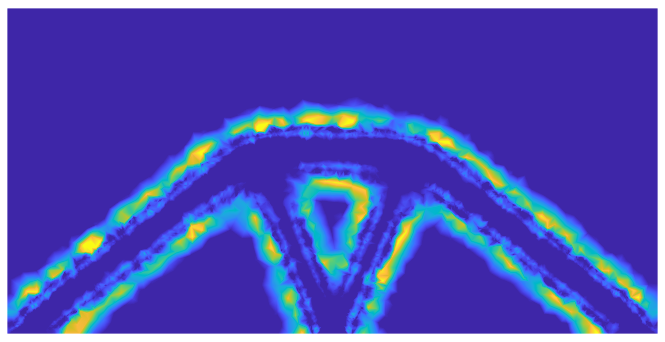} &\includegraphics[width=\figwid]{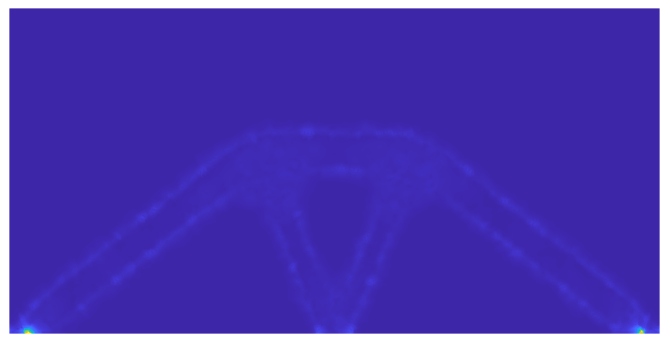}\\
\includegraphics[width=\figwid]{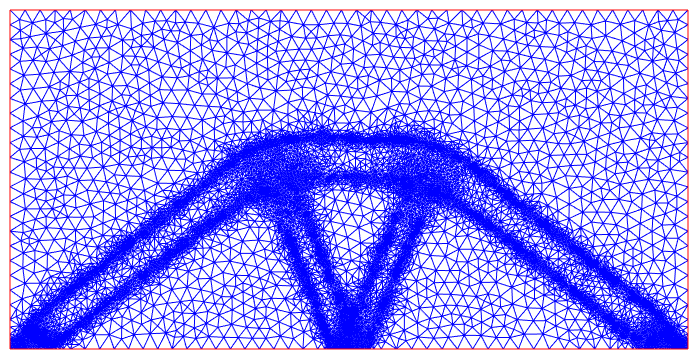} &\includegraphics[width=\figwid]{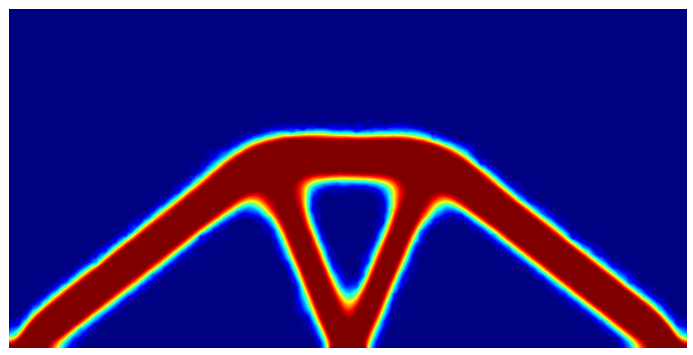}
&\includegraphics[width=\figwid]{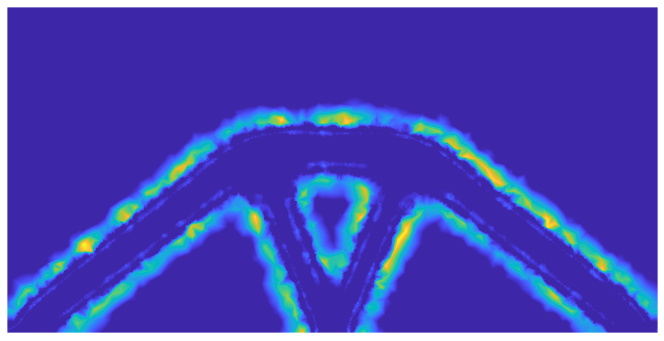} &\includegraphics[width=\figwid]{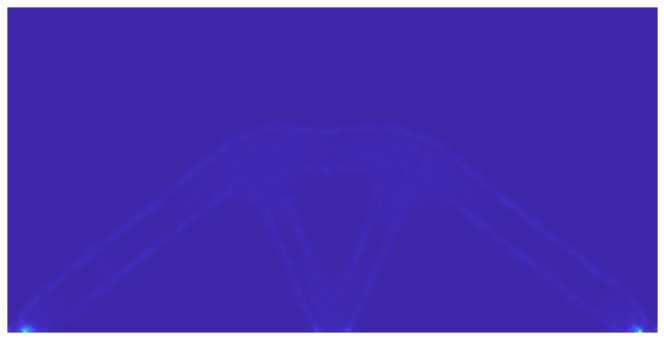}\\
& &\includegraphics[width=\figwid]{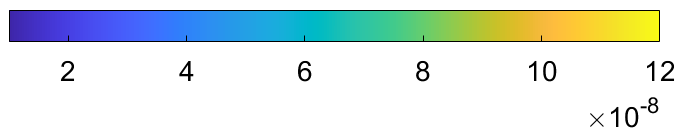}
&\includegraphics[width=\figwid]{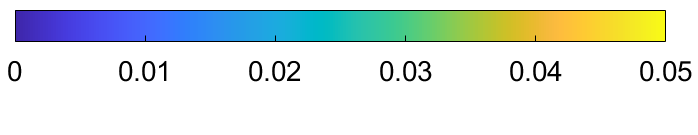}
\end{tabular}
\caption{The evolution of the mesh during the adaptive refinement procedure from $k=0$ (initial) to 5 for Example (c), with the number of vertices of each mesh being 1767, 2269, 3094, 4473, 6572 and 9959. The second, third and fourth columns show the optimized design $\rho_k^*$, error indicators $\eta_{k,1}$ and $\eta_{k,2}$ respectively.} \label{fig:Bridge_MeshAndDeisgn}
\end{figure}

\begin{figure}[htb!]
\centering \setlength{\tabcolsep}{0pt}
\begin{tabular}{cccccc}
\includegraphics[width=\figwid]{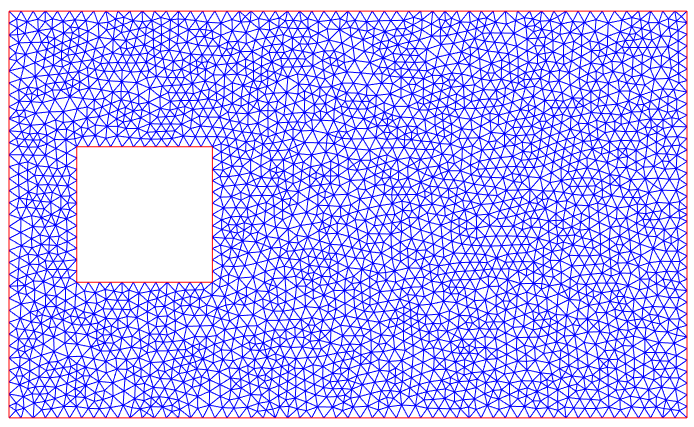} & \includegraphics[width=\figwid]{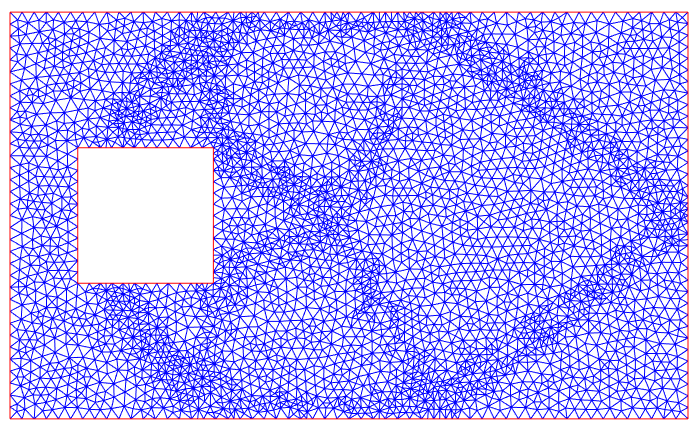}
&\includegraphics[width=\figwid]{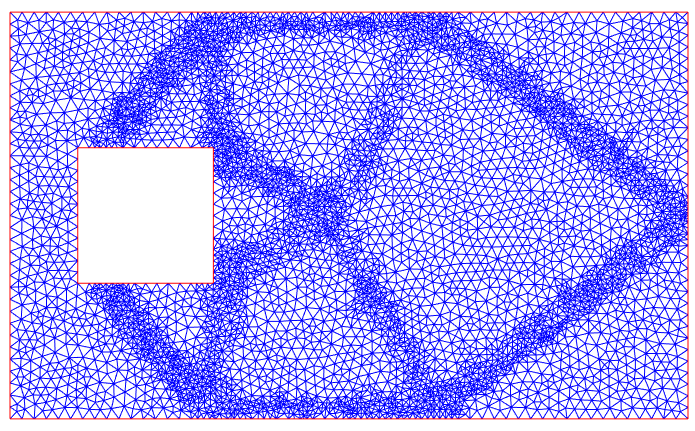} &\includegraphics[width=\figwid]{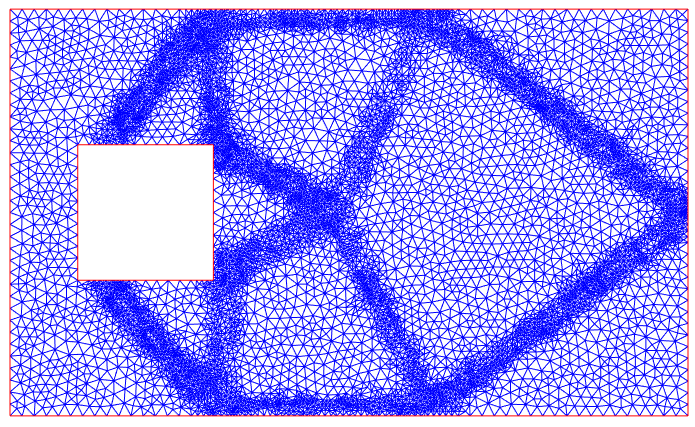}\\
\includegraphics[width=\figwid]{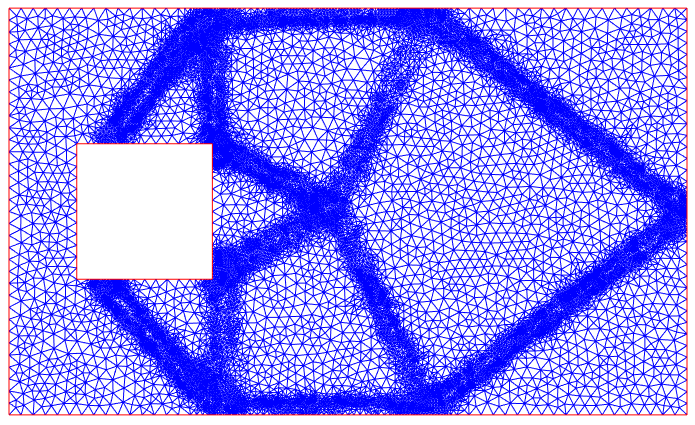}  &\includegraphics[width=\figwid]{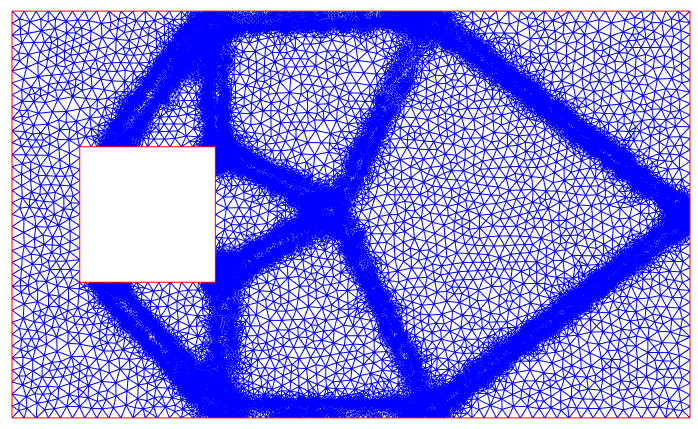}& \includegraphics[width=\figwid]{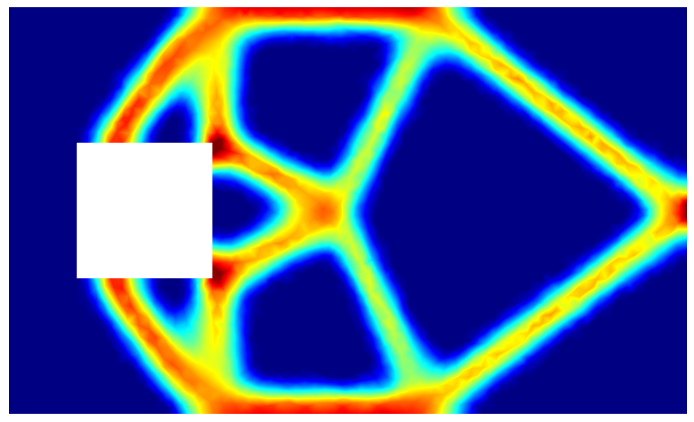} &\includegraphics[width=\figwid]{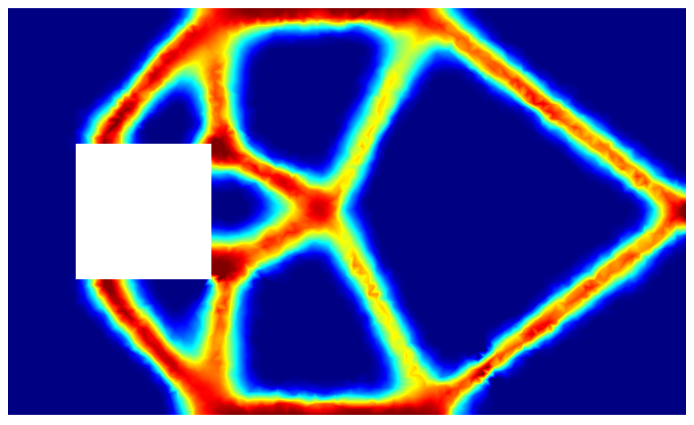}\\
\includegraphics[width=\figwid]{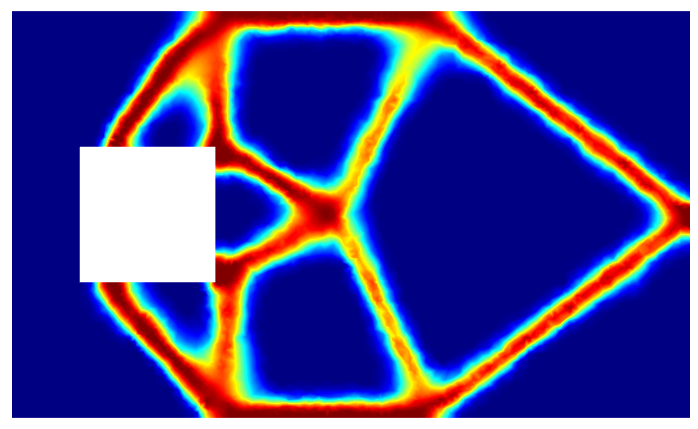}&\includegraphics[width=\figwid]{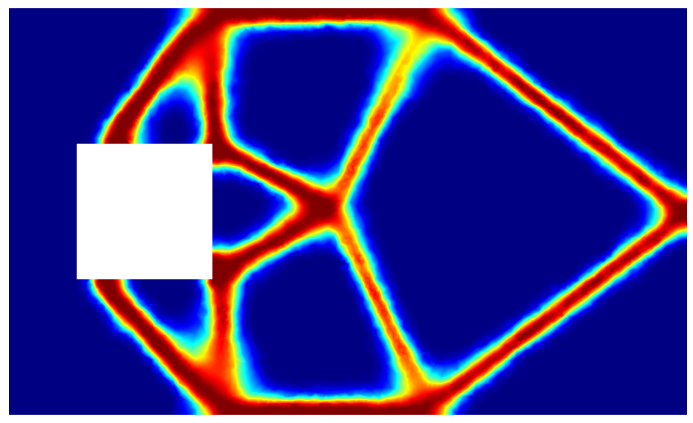}
&\includegraphics[width=\figwid]{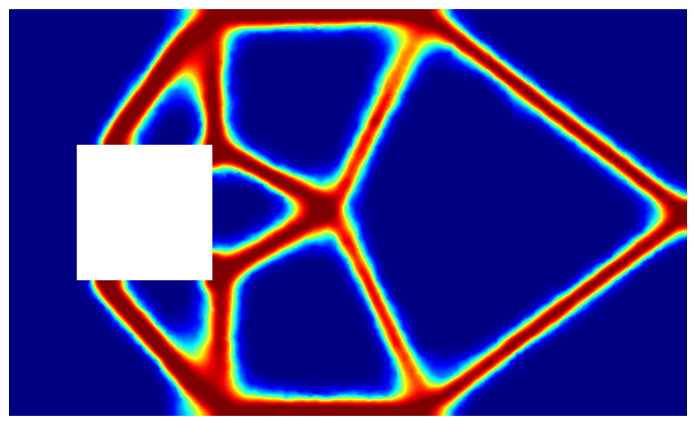}&\includegraphics[width=\figwid]{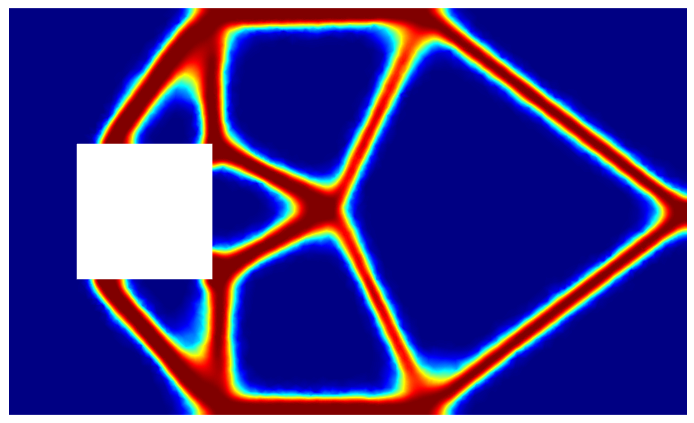}
\end{tabular}
\caption{The evolution of the mesh during the adaptive process and the optimized design $\rho_k^*$, from $k=0$ (initial) to 5 for Example (d), with the number of vertices of each mesh being 3051, 4089, 5790, 8828, 14327 and 24123.} \label{fig:SquareHole_MeshAndDesign}
\end{figure}

\begin{figure}[htb!]
\centering \setlength{\tabcolsep}{0pt}
\begin{tabular}{cccccc}
	 \includegraphics[width=\figwid]{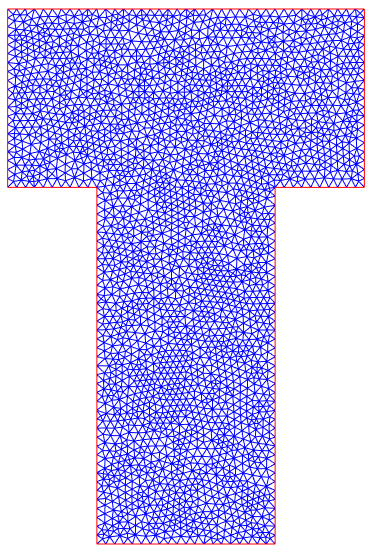} &\includegraphics[width=\figwid]{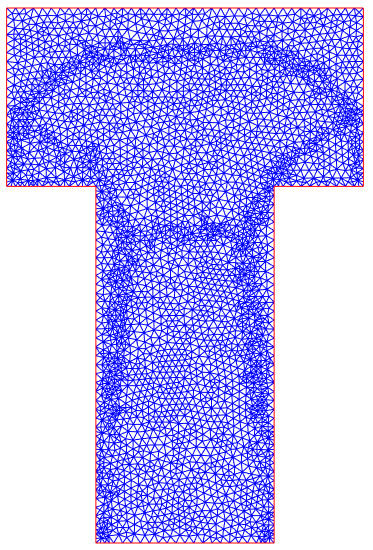}
	&\includegraphics[width=\figwid]{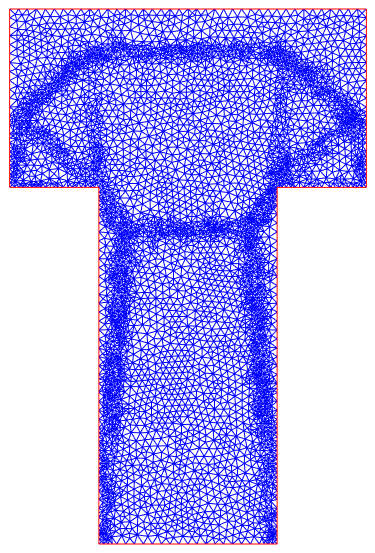} &\includegraphics[width=\figwid]{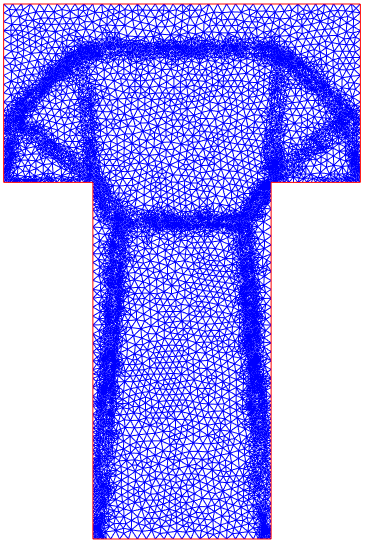}\\
	\includegraphics[width=\figwid]{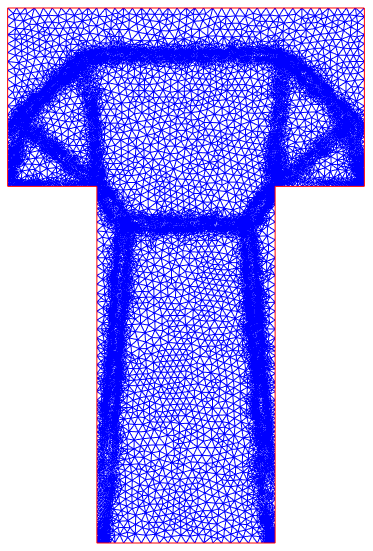}  &\includegraphics[width=\figwid]{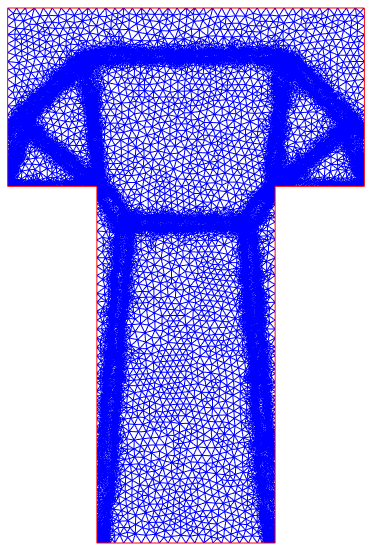}
	&\includegraphics[width=\figwid]{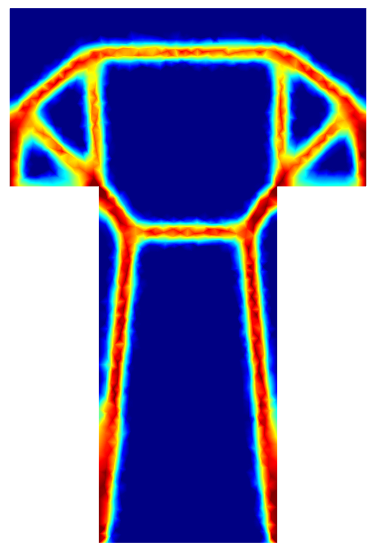}&\includegraphics[width=\figwid]{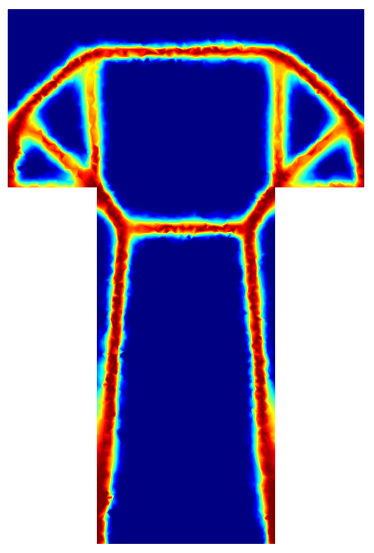}\\
	\includegraphics[width=\figwid]{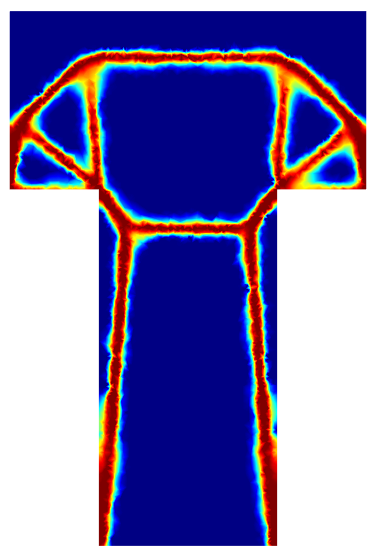} &\includegraphics[width=\figwid]{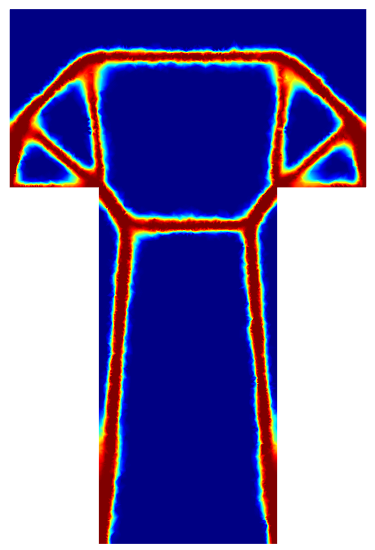}
	&\includegraphics[width=\figwid]{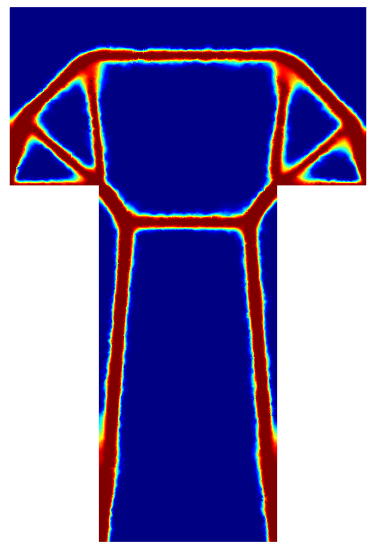}&\includegraphics[width=\figwid]{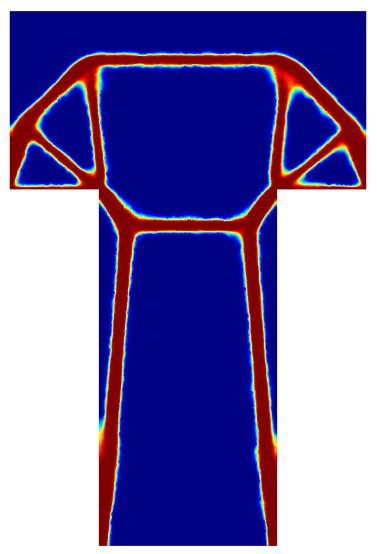}
\end{tabular}
\caption{The evolution of the mesh and optimal design $\rho_k^*$ during the adaptive process, from $k=0$ (initial) to 5 for Example (e), with the number of vertices of each mesh being 2878, 3749, 5231, 7692, 12125 and 19079.}
   \label{fig:T_MeshAndDesign}
\end{figure}

\begin{figure}[htb!]
\centering \setlength{\tabcolsep}{0pt}
\begin{tabular}{cccccc}
\hspace{-15mm}
\includegraphics[width=\figwid]{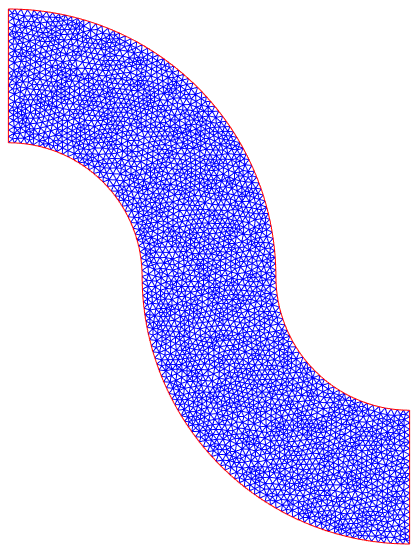} &\includegraphics[width=\figwid]{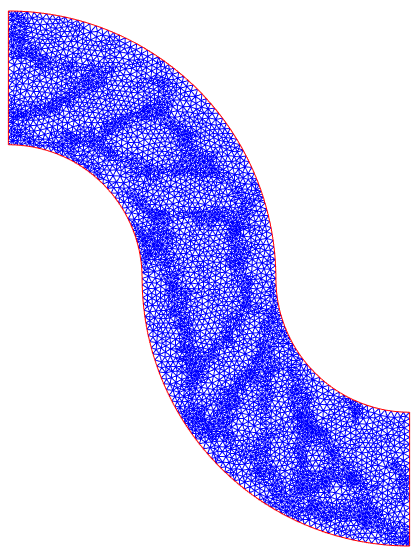}
&\includegraphics[width=\figwid]{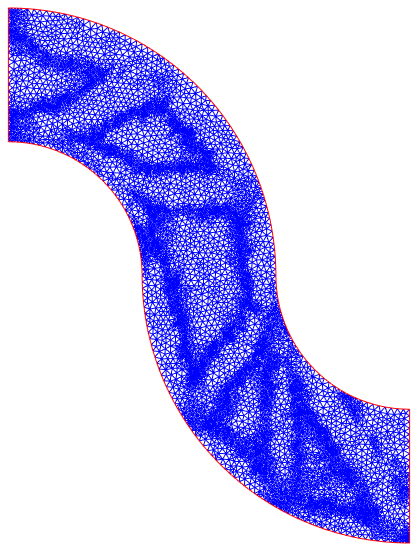}&\includegraphics[width=\figwid]{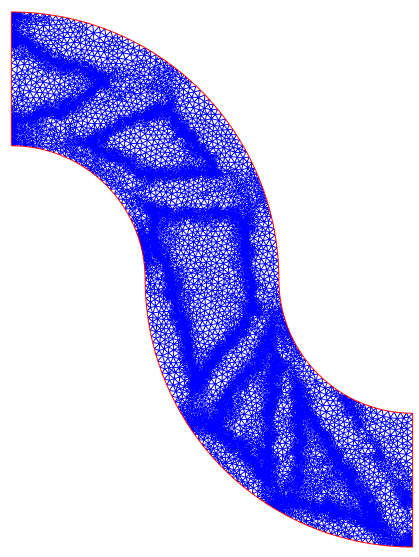} \\
\includegraphics[width=\figwid]{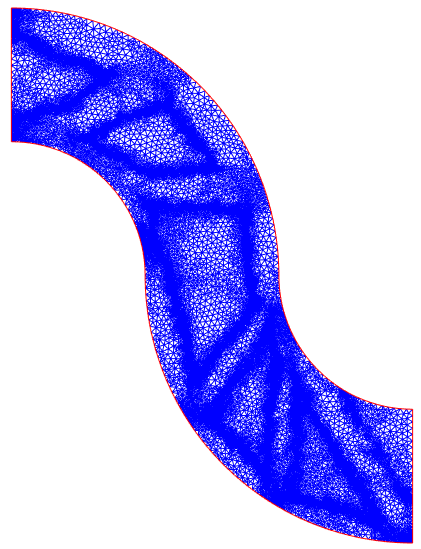} &\includegraphics[width=\figwid]{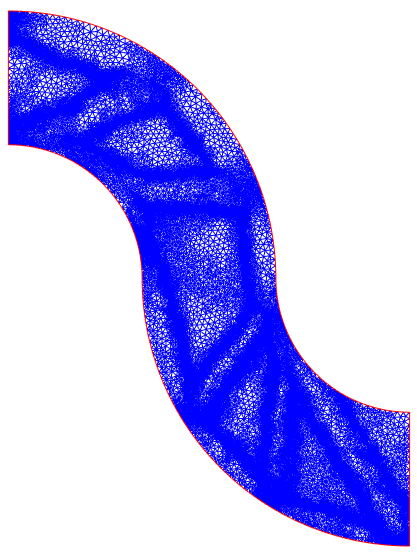} &\includegraphics[width=\figwid]{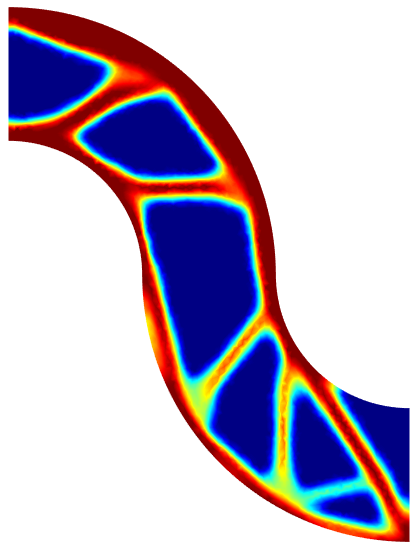}&\includegraphics[width=\figwid]{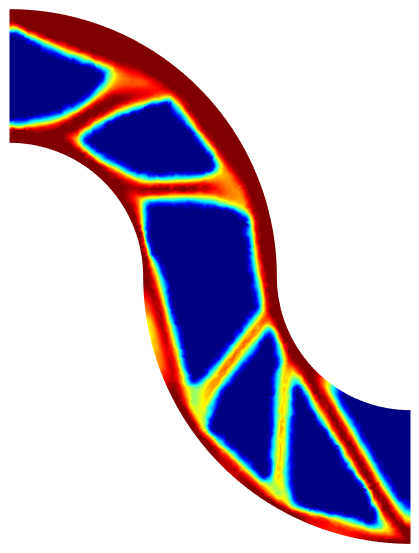}\\
\includegraphics[width=\figwid]{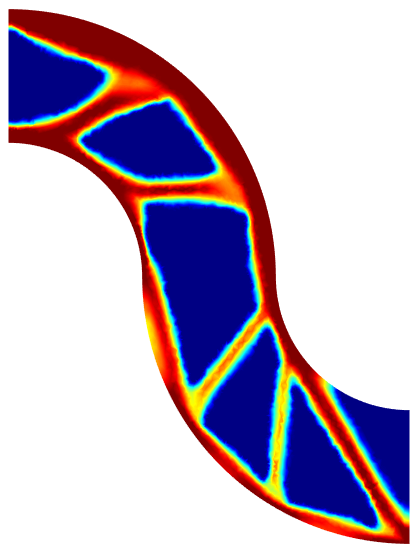}&\includegraphics[width=\figwid]{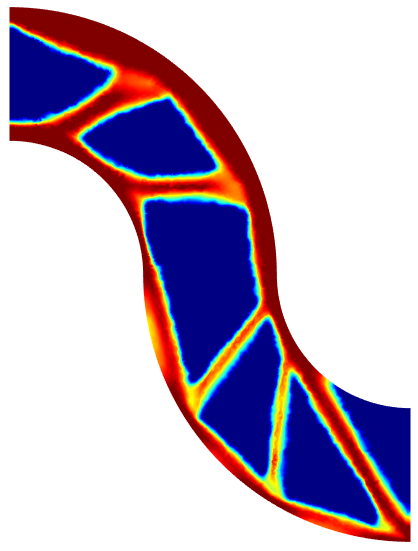}
&\includegraphics[width=\figwid]{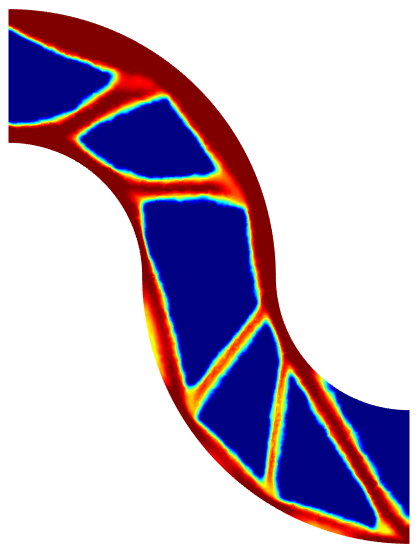}&\includegraphics[width=\figwid]{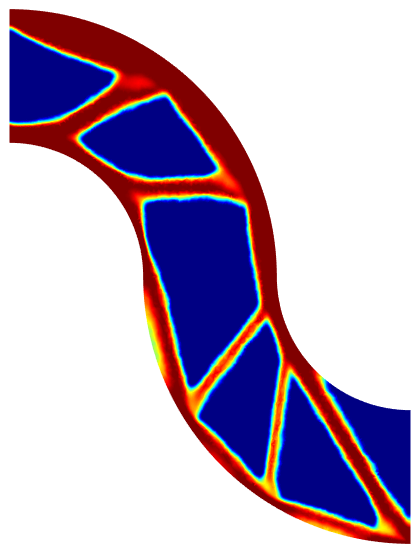}
\end{tabular}
\caption{The evolution of the mesh and the optimized design $\rho_k^*$ during the adaptive process, from $k=0$ (initial) to 5 for Example (f), with the number of vertices of each mesh being 3485, 5295, 7929, 12010, 18431 and 27439.} \label{fig:Curved_MeshAndDesign}
\end{figure}

In Figs. \ref{fig:Cantilever_MehsAndDeigsn}, \ref{fig:MBB2_MehsAndDeigsn} and \ref{fig:Bridge_MeshAndDeisgn}, we summarize the results for
Examples (a)--(c), including the evolution of meshes, optimal design and the error estimators $\eta_{k,1}(\rho_k^*,\bm{u}_k^*,T)$ and $\eta_{k,2}(\rho_k^*,\bm{u}_k^*,T)$ during the adaptive refinement process.
It is observed that as the adaptive loop proceeds, Algorithm \ref{alg_afem_topopt} can accurately capture interfaces between the material and the void regions, around which local refinements
are performed. More precisely, the first adaptive iterations $(k=0)$ yield a rough shape
of the desired design, and with the subsequent adaptive refinement, the algorithm produces improved boundary interfaces of the solid
material, due to the better resolution of the state variable $\bm{u}$. Notably, the additional refinements mainly take place around the
interfaces during the adaptive refinement procedure, which shows clearly the effectiveness of the proposed adaptive algorithm. These observations {also hold for} Examples (d)--(f); see Figs. \ref{fig:SquareHole_MeshAndDesign}, \ref{fig:T_MeshAndDesign} and \ref{fig:Curved_MeshAndDesign} for the evolution of the mesh and optimal design during the adaptive refinement process.

In Figs. \ref{fig:Cantilever_MehsAndDeigsn}, \ref{fig:MBB2_MehsAndDeigsn} and \ref{fig:Bridge_MeshAndDeisgn}, we show also the evolution of the error estimators $\eta_{k,1}$ and $\eta_{k,2}$ for Examples (a)--(c). By the very construction, the two estimators $\eta_{k,1}$ and $\eta_{k,2}$ play different roles: the estimator $\eta_{k,1}$ indicates the presence of the design edges, whereas {the  estimator $\eta_{k,2}$ concentrates in the region where the displacement field $\boldsymbol{u}$ varies most (or equivalently the gradient field $\boldsymbol\nabla \boldsymbol{u}$ exhibits the most dramatic change).} The magnitude of $\eta_{k,1}$ is observed to be much smaller than $\eta_{k,2}$ for all three examples (and also for the remaining examples, although not shown). The precise mechanism for this interesting phenomenon is still unknown. Nonetheless, this observation clearly indicates that one may only use $\eta_{k,2}$ as the error indicator to drive the adaptive procedure in practical computation. Also it is observed that the magnitudes of the estimators decrease steadily to zero throughout the adaptive process, concurring with the vanishing property numerically, cf. Lemma \ref{lem:err->0}.

\begin{figure}[htbp]
\centering\setlength{\tabcolsep}{0pt}
\begin{tabular}{c|c|c}
\toprule
\includegraphics[width=0.33\textwidth]{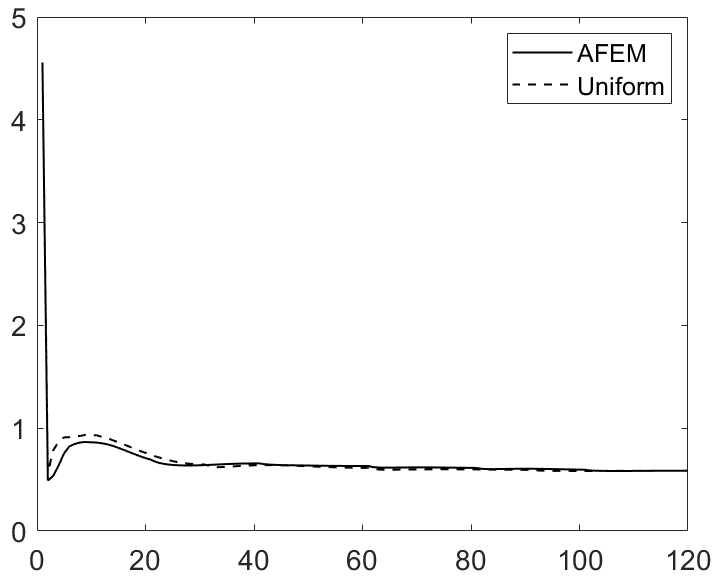}  & \includegraphics[width=0.33\textwidth]{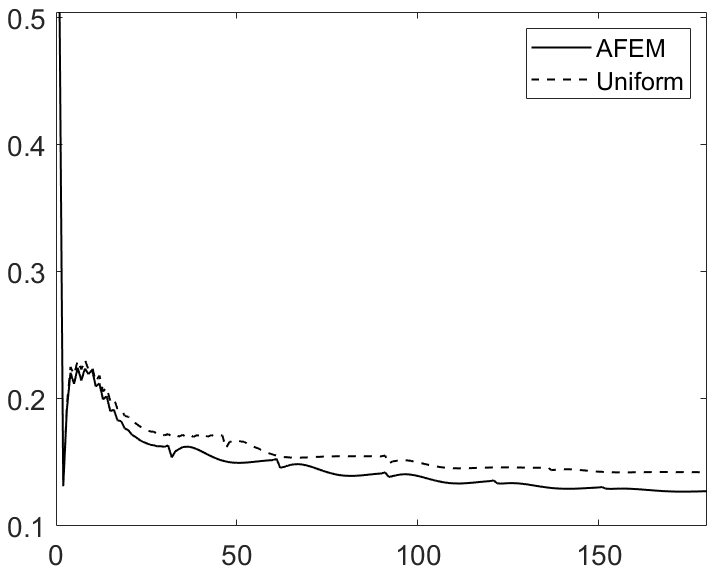} & \includegraphics[width=0.33\textwidth]{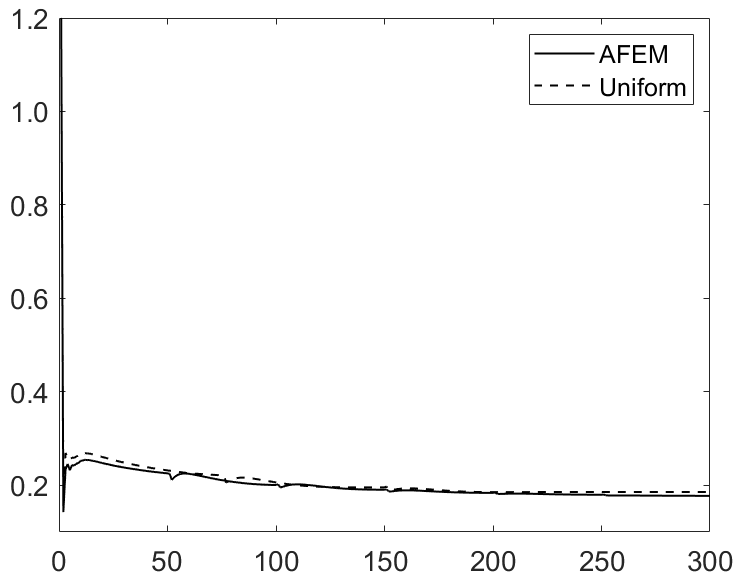} \\
\includegraphics[width=0.33\textwidth]{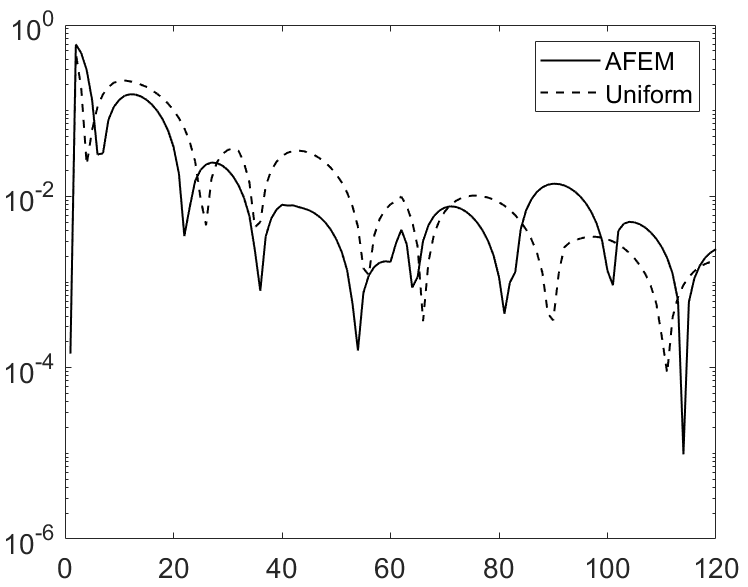} &\includegraphics[width=0.33\textwidth]{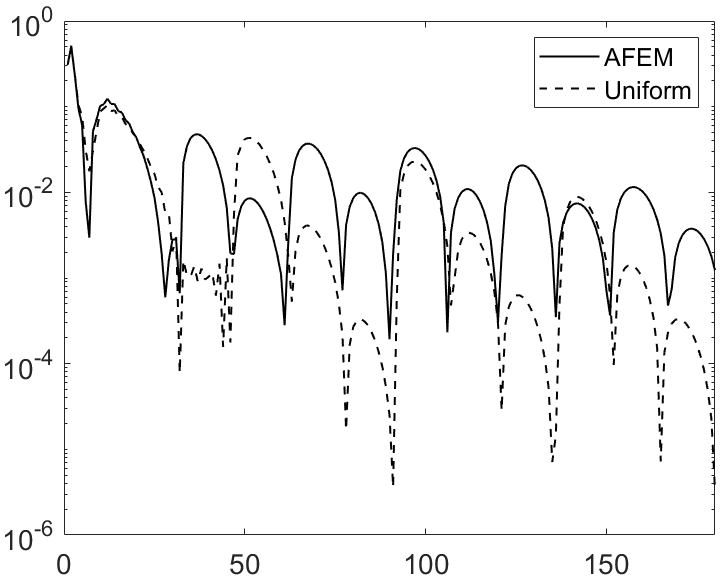} & \includegraphics[width=0.33\textwidth]{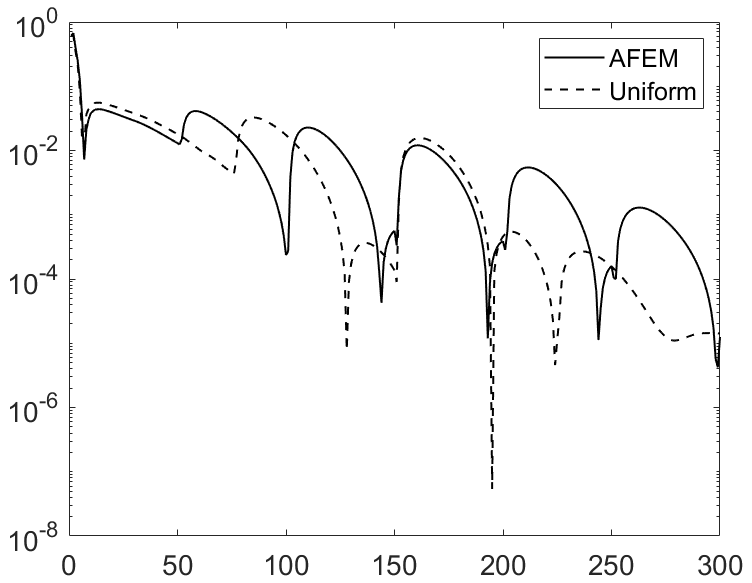} \\
(a) & (b) & (c)\\
\midrule
\includegraphics[width=0.33\textwidth]{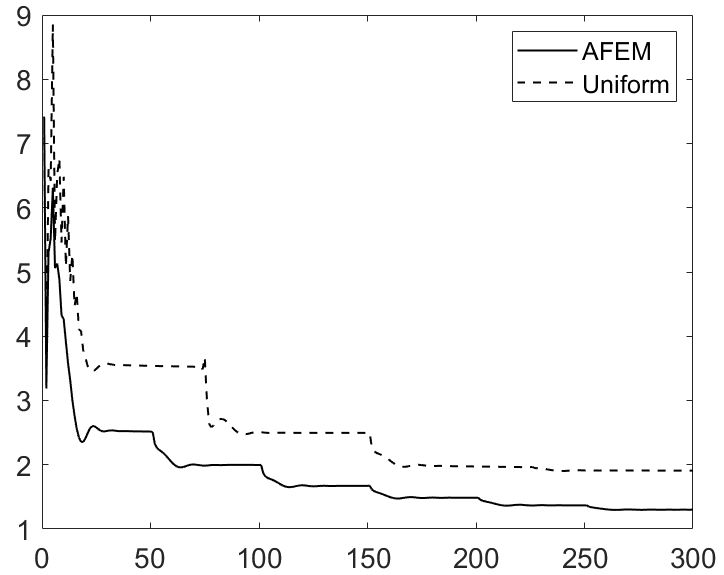} & \includegraphics[width=0.33\textwidth]{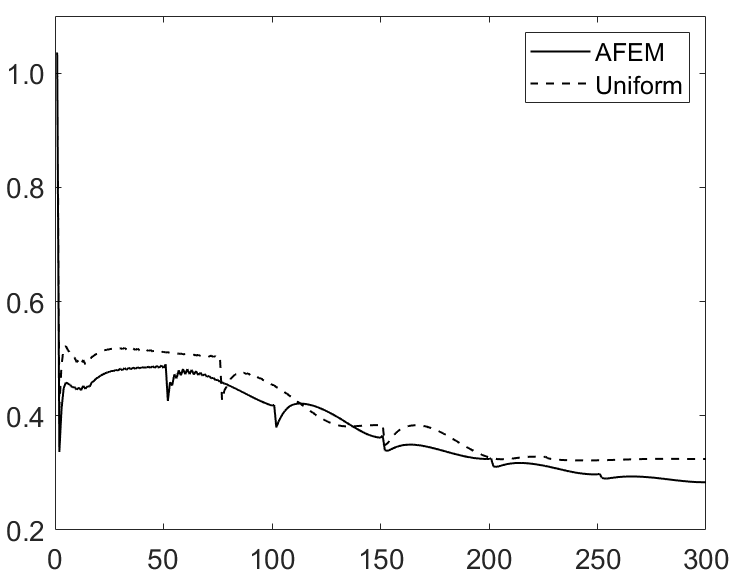} & \includegraphics[width=0.33\textwidth]{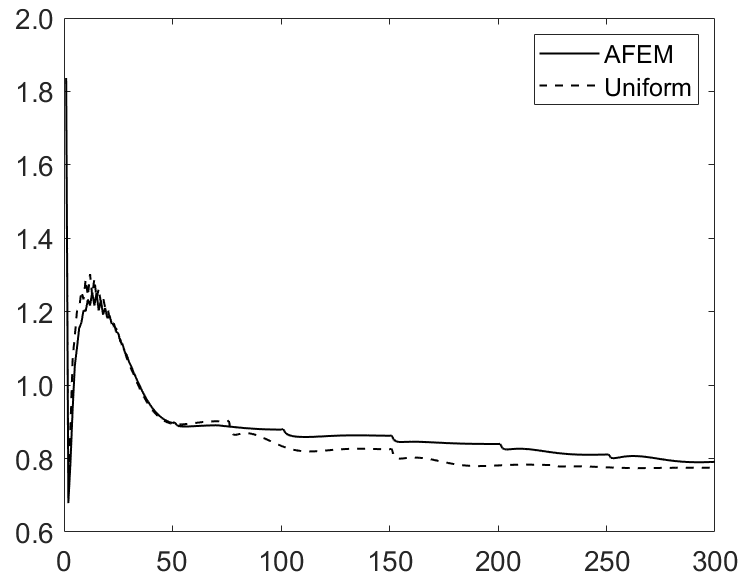} \\
\includegraphics[width=0.33\textwidth]{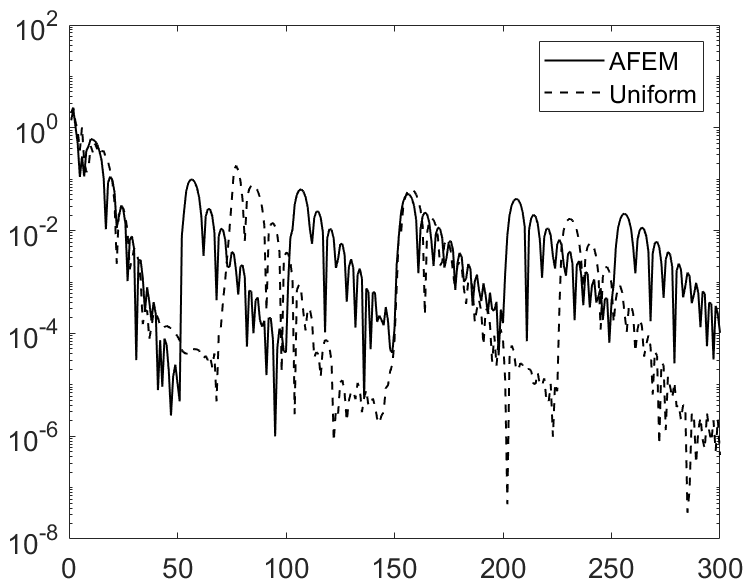} & \includegraphics[width=0.33\textwidth]{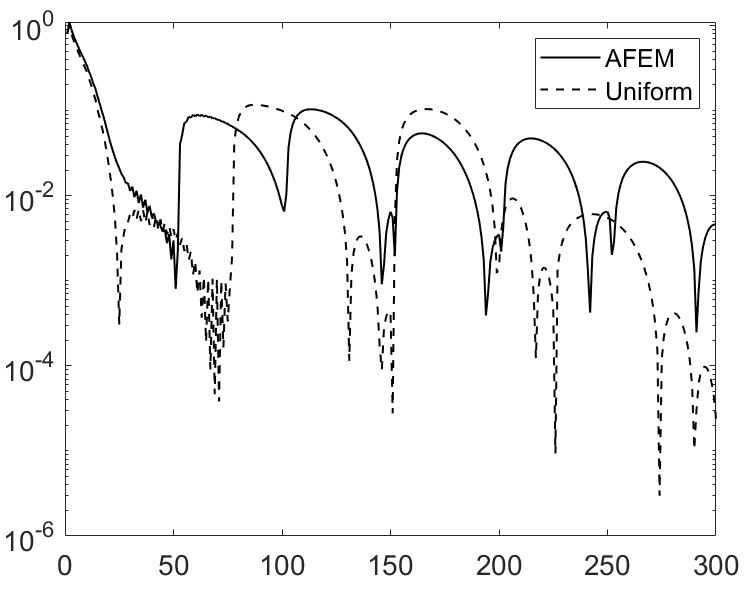} & \includegraphics[width=0.33\textwidth]{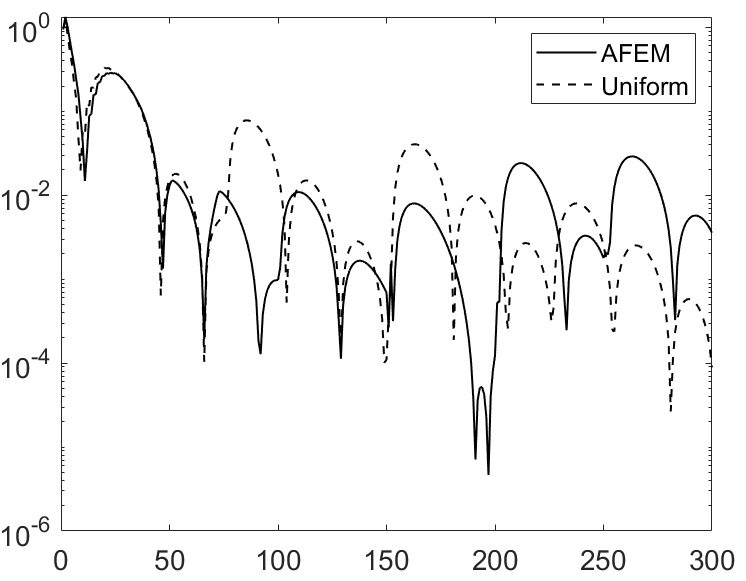}\\
(d) & (e) & (f)\\
\bottomrule
\end{tabular}
\caption{The convergence history of the objective value and the volume error as a function of iteration numbers in the augmented Lagrangian method.
All the examples take five and three refinements for the adaptive and uniform strategies respectively. For Example (a), the algorithm takes 20 (adaptive) and 30 (uniform) iterations per mesh, for Example (b), 30 (adaptive) and 45 (uniform), and for Examples (c)--(f), 50 (adaptive) and 75 (uniform) iterations.
\label{fig:conv}}
\end{figure}

\begin{figure}[hbt!]
\centering\setlength{\tabcolsep}{0pt}
\begin{threeparttable}
\begin{tabular}{|cc|cc|}
    \toprule
    \multicolumn{2}{c}{(a)}&
    \multicolumn{2}{c}{(b)}\\
    \cmidrule(lr){1-2} \cmidrule(lr){3-4}
    \includegraphics[width=0.24\textwidth]{CantileverAdap6thOptimizedDeisgn2.png}
	&\includegraphics[width=0.24\textwidth]{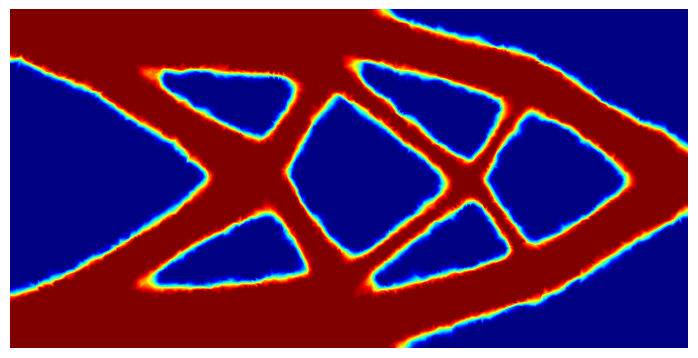}
   & \includegraphics[width=0.24\textwidth]{MBB2Adap6thOptimizedDeisgn2.png}
	&\includegraphics[width=0.24\textwidth]{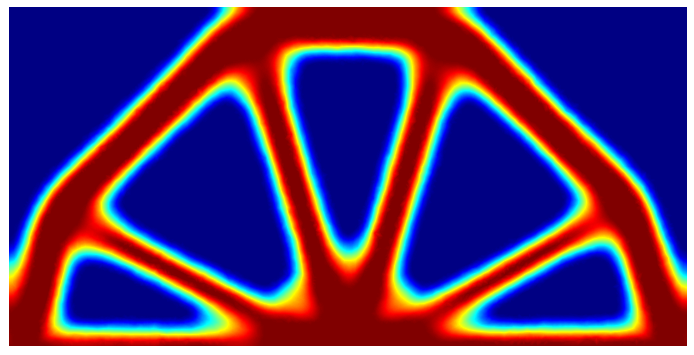}\\
 \midrule
 \multicolumn{2}{c}{(c)}&
    \multicolumn{2}{c}{(d)}\\
    \cmidrule(lr){1-2} \cmidrule(lr){3-4}
    \includegraphics[width=0.24\textwidth]{BridgeAdap6thOptimizedDeisgn2.png}
    &\includegraphics[width=0.24\textwidth]{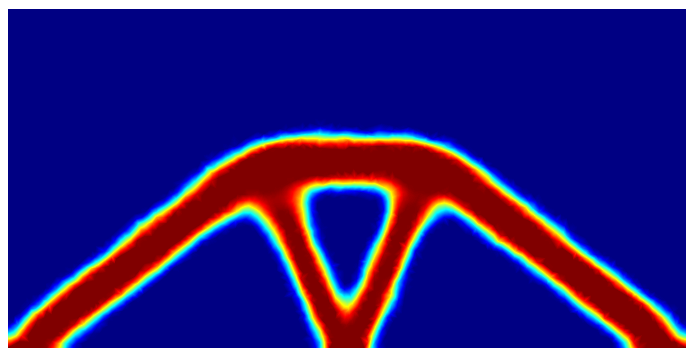}
    &\includegraphics[width=0.24\textwidth]{SquareHoleAdap6thOptimizedDeisgn2.png}
	&\includegraphics[width=0.24\textwidth]{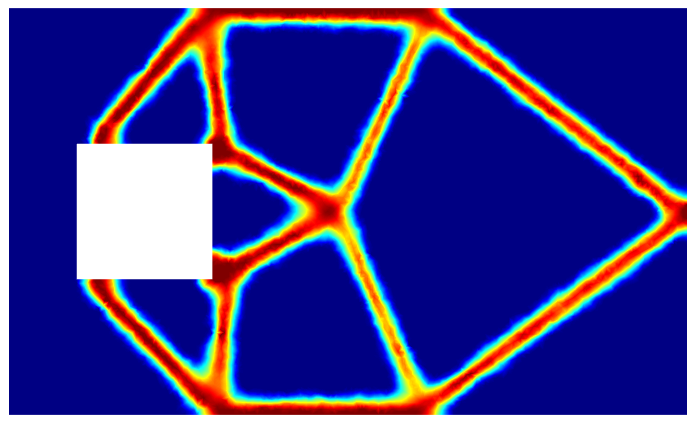}\\
 \midrule
 \multicolumn{2}{c}{(e)}&
    \multicolumn{2}{c}{(f)}\\
    \cmidrule(lr){1-2} \cmidrule(lr){3-4}
    \includegraphics[width=0.24\textwidth]{TAdap6thOptimizedDeisgn2.png}
    &\includegraphics[width=0.24\textwidth]{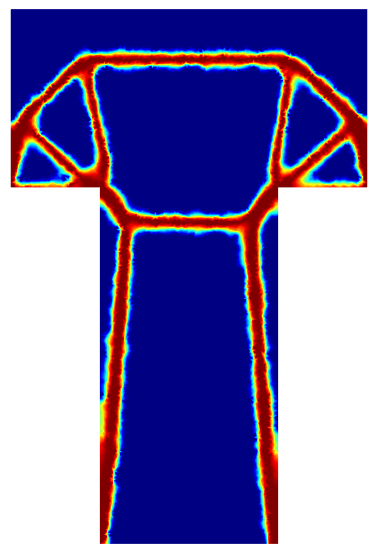}
    &\includegraphics[width=0.24\textwidth]{CurvedAdap6thOptimizedDeisgn2.png}
	&\includegraphics[width=0.24\textwidth]{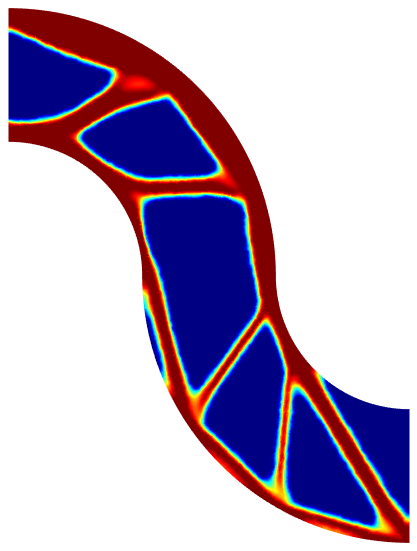}\\
    \bottomrule
    \end{tabular}
\end{threeparttable}
\caption{Qualitative results of the optimal designs on the final meshes by the adaptive and uniform refinements for Examples (a)--(f). Within each pair, the left and right designs are for the adaptive and uniform refinements respectively. \label{tab:results-fig}}
\end{figure}

The use of the adaptive algorithm does influence also the performance of Algorithm \ref{Alga}. In Fig. \ref{fig:conv}, we compare the decay of the objective value $\mathcal{J}_\gamma(\rho_n^{[k]})$ and the volume error by the adaptive and uniform refinements, where the refinement steps and iteration numbers of the augmented Lagrangian method for adaptive and uniform refinements are adjusted such that the total numbers of iterations and the numbers of vertices on the final meshes are comparable. It is observed that for both refinement strategies, the objective values perform equally well, and the volume constraint is well preserved during the refinement process by both strategies. Also it is consistently observed that at the beginning of the iteration, the objective value decreases rapidly, and then it stabilizes. In particular, the final objective values are always very close to each other for both refinement strategies, which agree with the fact that converged optimal designs are visually similar, cf. Fig. \ref{tab:results-fig}. {The objective value experiences slight spikes during the iteration when changing from one mesh to another one. This is attributed to the interpolation and projection errors when initiating the displacement and density fields on the fine mesh using that on the coarse mesh.}

Table \ref{tab:results} presents more quantitative results for both refinement strategies: the values of the objective functional 
(achieved at the final meshes generated by the uniform and adaptive refinements), the number of vertices on the final mesh, and the computing time 
(in seconds). Interestingly, the objective functional for the adaptive refinement is very close to the uniform one, or slightly smaller. Moreover, 
the computing time is comparable: the adaptive algorithm saves between 20\% and 200\% of the time compared to the uniform refinement. {To ensure a fair comparison, both uniform and adaptive refinements start from an identical initial coarse mesh and stop with nearly the same number of vertices after several mesh refinement steps, and both use the same number of total outer iterations (i.e., number of times to call an elasticity solver), cf. Fig. \ref{fig:conv}, and the same number of phase field inner iterations. The saving in computing time mainly stems from the fact that the increment of the degree of freedom of the mesh is milder for the adaptive refinement when compared with the uniform refinement, which implies that the associated  optimization problem is less expensive to minimize.} 
Fig. \ref{tab:results-fig} shows the optimized designs of the two methods and it can be seen that the boundary of adaptive algorithm is smoother, with fewer small oscillations along the interfaces, than that by the uniform refinement. This is attributed to the better resolution of the interface and the displacement field $\boldsymbol{u}$ by the adaptive strategy.

\begin{table}[hbt!]
\centering
\begin{threeparttable}
\caption{Quantitative results for the examples: the number of vertices of the final mesh,
the objective value and the total computing time (in second). Saving refers to the
saving in the overall computing time.\label{tab:results}}
\begin{tabular}{|c|ccc|ccc|c|}
    \toprule
    \multicolumn{1}{c}{}&
    \multicolumn{3}{c}{adaptive}&\multicolumn{3}{c}{uniform}\\
    \cmidrule(lr){2-4} \cmidrule(lr){5-7}
      Example & vertices & objective & time (sec) & vertices & objective & time (sec) & saving\\
    \midrule
    (a) & 16687 & 0.5868  & 394.03   & 16599 & 0.5863  & 509.05   &29.19\%\\
    (b) & 23314 & 0.1274  & 1220.36  & 22935 & 0.1423  & 1637.09  &34.15\%\\
    (c) & 9959  & 0.1774  & 316.48   & 9946  & 0.1855  & 701.86   &121.77\%\\
    (d) & 24123 & 1.3008  & 1549.34  & 23965 & 1.9089  & 4212.33  &171.88\%\\
    (e) & 19079 & 0.2836  & 1374.52  & 18564 & 0.3244  & 2632.16  &91.50\%\\
    (f) & 27439 & 0.7910  & 3814.19  & 27218 & 0.7748  & 5843.00  &53.19\%\\
    \bottomrule
    \end{tabular}
\end{threeparttable}
\end{table}

\section{Convergence analysis}\label{sec:conv}

Now we analyze the convergence of Algorithm \ref{alg_afem_topopt}.
First, we introduce an auxiliary minimum compliance problem over a limiting admissible set $\widetilde{\mathcal{A}}_\infty$ given by
Algorithm \ref{alg_afem_topopt} (cf. \eqref{medmin} below). Then using techniques in nonlinear
optimization, the sequence of adaptively generated solutions is shown to contain a
subsequence convergent to the solution of the auxiliary problem in $H^1(\Omega)$.
Then the desired result follows by proving that the solution also satisfies the necessary
optimality system \eqref{optsys}. Note that two error indicators $\eta_{k,1}$ and
$\eta_{k,2}$ are involved in Algorithm \ref{alg_afem_topopt}, however, in the process
of our analysis, only one subsequence convergence is obtained.

The convergence analysis employs the following auxiliary problem:
\begin{subequations}
\begin{align}
\underset{\rho_\infty \in \widetilde{\mathcal{A}}_\infty}\min \left\{
\J_{\gamma,\infty} (\rho_\infty) = \displaystyle \int_{\Omega} \rho_\infty^p\C_0 \strain(\bold{u}_\infty) :
\strain(\bold{u}_\infty) \dd x + \widetilde{\beta}\left( \dfrac{\gamma}{2}\int_\Omega |\bold{\nabla}\rho_\infty|^2 + \dfrac{1}{\gamma}W(\rho_\infty)\dx  \right)  \right\},\label{medmin}\\
\mbox{subject to} \displaystyle \int_\Omega \rho_\infty^p\C_0 \strain(\bold{u}_\infty) : \strain(\bold{v}_\infty) \dd x  = \int_\Omega \bold{f} \cdot \bold{v}_\infty \dd x
+ \int_{\Gamma_N} \bold{g}\cdot \bold{v}_\infty \dd s\quad  \forall \bold{v}_\infty \in \bold{V}_\infty^D,\label{medmin-state}
\end{align}
\end{subequations}
where $\Atilde_\infty$ and $\bold{V}_\infty^D$ are closures of $\{\Atilde_k,\bold{V}_k^D\}_{k\geq 0}$ given by Algorithm \ref{alg_afem_topopt}, i.e.,
\[
\widetilde{\mathcal{A}}_{\infty} :=\overline{\bigcup_{k \geq 0} \widetilde{\mathcal{A}}_{k}}\quad\text { in } H^{1}(\Omega) \text {-norm}\quad\mbox{and}\quad \bold{V}^D_{\infty} :=\overline{\bigcup_{k \geq 0} \bold{V}_{k}^D}\quad \text {in } \bold{H}^1(\Omega) \text{-norm}.
\]
{By the very construction of the space $\boldsymbol{V}^D_\infty$, the sequence of spaces $\{\boldsymbol{V}_k^D\}_{k\geq0}$ is dense in $\bm{V}_\infty^D$ in the $\boldsymbol{H}^1(\Omega)$-norm. This density result will be used frequently below.}

\begin{lemma}\label{lem:convexclosed}
$\bold{V}_\infty^D$ is a closed subspace of $\bold{V}^D$ and $\Atilde_\infty$ is a closed convex subset of $\Atilde$.
\end{lemma}
\begin{proof}
By definition, $\Atilde_\infty$ (respectively $\bold{V}_\infty^D$) is closed in $H^1(\Omega)$ (respectively
$\bold{V}^D$). For any $\rho$ and $\nu$ in $\Atilde_{\infty}$, there exist two sequences $\{\rho_{k}\}_{k\geq0}$ and $\{\nu_{k}\}_{k\geq0}
\subset\bigcup_{k\geq 0}\widetilde{\mathcal{A}}_{k}$ such that $\rho_{k}\rightarrow\rho$ and $\nu_{k}\rightarrow\nu$ in $H^{1}(\Omega)$.
By the convexity of the set $\Atilde_{k}$, $\{t\rho_{k}+(1-t)\nu_{k}\}_{k\geq0}\subset\bigcup_{k\geq 0}
\Atilde_{k}$ for any $t\in(0,1)$. Then $t\rho_{k}+(1-t)\nu_{k}\rightarrow t\rho+(1-t)\nu$ in
$H^{1}(\Omega)$, i.e. $t\rho+(1-t)\nu\in\Atilde_{\infty}$ for any $t\in (0,1)$  since $\Atilde_\infty$ is closed. Hence $\Atilde_{
\infty}$ is convex. Moreover, we have $\int_\Omega \rho_k \dd x \to \int_\Omega \rho \dd x$ and $\rho_{k}\rightarrow\rho$ almost
everywhere in $\Omega$ after (possibly) passing to a subsequence, which, along with the constraint $\int_\Omega \rho_k \dd x =
V^0$ and $\underline{\rho}\leq\rho_k\leq 1$ almost everywhere in $\Omega$, indicates that $\int_\Omega \rho \dd x = V^0$ and
$\underline{\rho}\leq\rho\leq 1$ almost everywhere in $\Omega$. Lastly, the fact that $\Atilde_{\infty} \subset H^1(\Omega)$ implies
$\Atilde_{\infty}\subseteq\Atilde$.
\end{proof}

\begin{theorem}\label{thm:medmin}
There exists at least one minimizer to problem \eqref{medmin}-\eqref{medmin-state}.
Moreover, the sequence $\{(\rho_k^\ast,\bold{u}_k^\ast)\}_{k\geq 0}$
of discrete minimizers to problem \eqref{dismin}--\eqref{dismin-state} generated by Algorithm \ref{alg_afem_topopt} contains a subsequence $\{(\rho_{k_j}^\ast,
\bold{u}_{k_j}^\ast)\}_{j\geq 0}$ convergent to a minimizer $(\rho^\ast_\infty,\bold{u}_\infty^\ast)$ to problem \eqref{medmin}-\eqref{medmin-state}:
\begin{equation*}
\rho_{k_{j}}^\ast\rightarrow\rho_{\infty}^\ast\quad\mbox{ in}~H^1(\Omega),
\quad \rho_{k_{j}}^\ast\rightarrow\rho_{\infty}^\ast\quad\mbox{a.e. in}~\Omega,
\quad  \bold{u}_{k_{j}}^\ast\rightarrow \bold{u}_{\infty}^\ast\quad\text{in}~\bold{H}^1(\Omega).
\end{equation*}
\end{theorem}
\begin{proof}
Since the constant function $\rho\equiv \frac{V^0}{|\Omega|} \in \Atilde_k$ for all $k$, the minimizing property of $\J_{\gamma,k}
    (\rho_k^\ast)$ over $\Atilde_k$ and the stability \eqref{stab:disvp} for the associated discrete displacement field $\bold{u}_k \left(\frac{V^0}{|\Omega|}\right)$ imply that
\[
    \J_{\gamma,k}
    (\rho_k^\ast)\leq \J_{\gamma,k}\left(\frac{V^0}{|\Omega|}\right)\leq c.
\]
The box constraint in $\Atilde_k$ implies that $\| \rho_k^\ast
\|_{L^2(\Om)}\leq c$, and thus $\{\rho_{k}^\ast\}_{k\geq 0}$ is uniformly bounded in $H^1(\Omega)$. By Lemma \ref{lem:convexclosed}
and Sobolev compact embedding theorem, there exist a subsequence, denoted by $\{\rho_{k_j}^\ast\}_{j\geq 0}$, and some $\rho^\ast\in
\Atilde_\infty$ such that
\begin{equation}\label{pf:medmin_01}
\rho_{k_j}^\ast \rightharpoonup \rho^\ast\quad \mbox{in}~H^1(\Om),\quad
\rho_{k_j}^\ast \to \rho^\ast\quad \mbox{in}~L^2(\Om),\quad
\rho_{k_j}^\ast \to \rho^\ast\quad \mbox{a.e. in}~\Om.
\end{equation}
Next we define a discrete analogue of the variational problem in \eqref{medmin} with $\rho_\infty=\rho^\ast$: find
$\bold{u}_{k_j}\in \bold{V}^D_{k_j}$ such that
\begin{equation}\label{pf:medmin_02}
\int_\Omega (\rho^\ast)^p \C_0 \strain(\bold{u}_{k_j}) : \strain (\bold{v}) \dx = \int_\Omega \bold{f} \cdot \bold{v} \dx
+ \int_{\Gamma_N} \bold{g} \cdot \bold{v} \dd s, \quad \forall \bold{v}\in\bold{V}_{k_j}^D.
\end{equation}
That is, $\bold{u}_{k_j}$ is the finite element approximation of the solution $\bold{u}_\infty^\ast$ to the auxiliary problem
\eqref{medmin-state}) (with $\rho_\infty=\rho^\ast$).
By Korn's inequality, C\'{e}a's lemma and {the density of $\bigcup_{k\geq0}\bold{V}_k^D$ in $\bold{V}_\infty^D$}, we have
\begin{equation}\label{pf:medmin_03}
\|\bold{u}^\ast_\infty-\bold{u}_{k_j}\|_{\boldsymbol{H}^1(\Omega)}\leq c
\inf_{\bold{v}\in\bold{V}_{k_j}}\|\bold{u}_{\infty}^\ast-\bold{v}\|_{\boldsymbol{H}^1(\Omega)}.
\end{equation}
Meanwhile, taking $\bold{v}=\bold{u}_{k_j} - \bold{u}^\ast_{k_j}$ in \eqref{pf:medmin_02} and \eqref{dismin-state} yields
\begin{align*}
&\quad
\displaystyle\int_\Omega (\rho^\ast_{k_j})^p \C_0 \strain(\bold{u}_{k_j} - \bold{u}^\ast_{k_j}): \strain(\bold{u}_{k_j} -
\bold{u}^\ast_{k_j})  \dx
= \int_\Omega ((\rho^\ast_{k_j})^p-(\rho^{\ast})^p) \C_0 \strain (\bold{u}_{k_j}) :
\strain(\bold{u}_{k_j} - \bold{u}^\ast_{k_j}) \dx \\
&= \int_\Omega ((\rho^\ast_{k_j})^p-(\rho^{\ast})^p) \C_0 \strain (\bold{u}_{k_j} - \bold{u}^\ast_{\infty}) :
\strain(\bold{u}_{k_j} - \bold{u}^\ast_{k_j}) \dx  + \int_\Omega ((\rho^\ast_{k_j})^p-(\rho^{\ast})^p) \C_0
\strain (\bold{u}_{\infty}^\ast) : \strain(\bold{u}_{k_j} - \bold{u}^\ast_{k_j}) \dx.\\
&\leq
(\|((\rho^\ast_{k_j})^p-(\rho^{\ast})^p) \C_0 \strain (\bold{u}_{k_j} - \bold{u}^\ast_{\infty})\|_{L^{2\times2}(\Omega)}+ \|((\rho^\ast_{k_j})^p-(\rho^{\ast})^p) \C_0 \strain (\bold{u}_{\infty}^\ast)\|_{L^{2\times2}(\Omega)} )\|\strain(\bold{u}_{k_j} - \bold{u}^\ast_{k_j})\|_{L^{2\times2}(\Omega)}.
\end{align*}
By the pointwise convergence in \eqref{pf:medmin_01} and Lebesgue's dominated convergence theorem \cite[Theorem 1.19]{EvansGariepy:2015}, we get
\begin{align*}
\|((\rho^\ast_{k_j})^p-(\rho^{\ast})^p) \C_0 \strain (\bold{u}_{k_j} - \bold{u}^\ast_{\infty})\|_{L^{2\times2}(\Omega)} &\leq
\|\C_0 \strain (\bold{u}_{k_j} - \bold{u}^\ast_{\infty})\|_{L^{2\times2}(\Omega)} \to 0,\\
\|((\rho^\ast_{k_j})^p-(\rho^{\ast})^p) \C_0 \strain (\bold{u}^\ast_{\infty})\|_{L^{2\times2}(\Omega)}  & \to 0,
\end{align*}
which, together with the inequality $\rho\geq \underline{\rho}$ and Korn's inequality, implies
$    \|\bold{u}^\ast_{k_j}-\bold{u}_{k_j}\|_{\bold{H}^1(\Omega)}\to 0$. This, \eqref{pf:medmin_03}, and
the triangle inequality show
\begin{equation}\label{pf:medmin_04}
\|\bold{u}^\ast_\infty-\bold{u}^\ast_{k_j}\|_{\bold{H}^1(\Omega)}
\to 0.
\end{equation}
Since $\bold{u}^\ast_{k_j}$ solves \eqref{dismin-state} with $\rho_\cT = \rho_{k_j}^\ast$ and $\bold{u}^\ast_\infty$
solves \eqref{medmin-state} with $\rho_\infty = \rho^\ast$, it follows from \eqref{pf:medmin_04} that
\begin{equation}\label{pf:medmin_05}
\begin{aligned}
	&\quad\int_\Omega (\rho_{k_j}^\ast)^p \C_0 \strain(\bold{u}_{k_j}^\ast):\strain(\bold{u}_{k_j}^\ast) \dx = \int_{\Omega} \bold{f} \cdot \bold{u}_{k_j}^\ast \dx
	+ \int_{\Gamma_N} \bold{g} \cdot \bold{u}_{k_j}^\ast \dd s \\
	& \to   \int_{\Omega} \bold{f} \cdot \bold{u}_{\infty}^\ast \dx
	+ \int_{\Gamma_N} \bold{g} \cdot \bold{u}_{\infty}^\ast \dd s = \int_\Omega (\rho^\ast)^p \C_0 \strain(\bold{u}_{\infty}^\ast):\strain(\bold{u}_{\infty}^\ast) \dx.
\end{aligned}
\end{equation}
Using the pointwise convergence in \eqref{pf:medmin_01} and Lebesgue's dominated convergence theorem again gives
\begin{equation}\label{pf:medmin_06}
\int_\Omega W(\rho_{k_j}^\ast) \dx \to \int_\Omega W(\rho^\ast) \dx.
\end{equation}
Now, for any $\rho_\infty \in \Atilde_\infty$ there exists a sequence $\{\rho_k\}_{k\geq 0}$ such that $\rho_k\to \rho_\infty$
in $H^1(\Omega)$. Using the above argument and a standard subsequence argument gives
$$\int_\Omega (\rho_{k})^p \C_0 \strain(\bold{u}_{k})
:\strain(\bold{u}_{k})\dx \to \int_\Omega (\rho_{\infty})^p \C_0 \strain(\bold{u}_{\infty}):\strain(\bold{u}_{\infty})\dx\quad\mbox{and}
\int_\Omega W(\rho_{k}) \dx \to \int_\Omega W(\rho_{\infty}) \dx.$$
Now combining \eqref{pf:medmin_05} with \eqref{pf:medmin_06},
the weak lower semicontinuity of $H^1(\Omega)$-seminorm and the minimizing property of $\rho^\ast_{k}$ to $\J_{\gamma,k}$ over $\Atilde_k$, we get
\[
\begin{aligned}
\J_{\gamma,\infty}(\rho^\ast)&\leq \liminf_{j\to\infty}\J_{\gamma,k_j}(\rho^\ast_{k_j})\leq
\limsup_{j\to\infty}\J_{\gamma,k_j}(\rho^\ast_{k_j})\nonumber\\
&\leq\limsup_{k\to\infty}\J_{\gamma,k}(\rho^\ast_{k})\leq
\limsup_{k\to\infty}\J_{\gamma,k}(\rho_{k})=\J_{\gamma,\infty}(\rho_\infty),\quad\forall\rho_\infty\in \Atilde_{\infty}.
\end{aligned}
\]
Since $\rho^\ast\in \Atilde_\infty$, $\rho_\infty^\ast:=\rho^\ast$ is a minimizer of $\J_{\gamma,\infty}$ over $\Atilde_{\infty}$.
Moreover, taking $\rho_{\infty}=\rho_\infty^\ast$ in the last inequality, we obtain
$$\J_{\gamma,k_j}(\rho_{k_j}^\ast)\to\J_{\gamma,\infty}(\rho^\ast_\infty).$$
It follows from \eqref{pf:medmin_05} and \eqref{pf:medmin_06} that
$$\|\nabla\rho_{k_j}^\ast
\|^2_{\bold{L}^2(\Om)}\to \|\nabla\rho_{\infty}^\ast\|^2_{\bold{L}^2(\Om)}.$$
This directly implies  $\rho_{k_{j}}^\ast\rightarrow\rho_{\infty}^\ast$ in $H^1(\Omega)$,
and completes the proof of the theorem.
\end{proof}

In view of Theorem \ref{thm:medmin}, it suffices to prove that the tuple $(\rho_\infty^\ast,\bold{u}_\infty^\ast)$
satisfies \eqref{optsys}. We make use of the marking strategy \eqref{eqn:marking1}-\eqref{eqn:marking2}.
In the proof, we denote by $D_{T}$ the union of all elements in $\cT$ with a non-empty
intersection with an element $T\in\cT$ and by $\omega_F$ the union of elements in $\cT$ sharing
a common edge with $F\in\mathcal{F}_\cT$. 

\begin{lemma}\label{lem:err->0}
Let $\{\cT_k,\Atilde_{k}\times\boldsymbol{V}_{k}^D,(\rho^{\ast}_k,\bold{u}^{\ast}_k)\}_{k\geq 0}$ be the sequence of
meshes, discrete admissible sets, finite element spaces and discrete solutions generated by Algorithm
\ref{alg_afem_topopt}. If $p>1$ is an integer, then for the convergent subsequence $\{(\rho^{\ast}_{k_j},\bold{u}^{\ast}_{k_j})\}_{j\geq0}$
in Theorem \ref{thm:medmin}, the corresponding sequences of the estimators converge to zero, i.e.,
\begin{equation}\label{err->0}
\lim_{j\to \infty} \eta_{k_j,1}(\rho_{k_j}^\ast,\bold{u}^\ast_{k_j})=\lim_{j\to \infty} \eta_{k_j,2}(\rho_{k_j}^\ast,\bold{u}^\ast_{k_j})= 0.
\end{equation}
\end{lemma}

\begin{proof}
We relabel $\{(\rho^{\ast}_{k_j},\bold{u}^{\ast}_{k_j})\}_{j\geq0}$ as $\{(\rho^{\ast}_{k},\bold{u}^{\ast}_{k})
\}_{k\geq0}$. Since the module \textsf{MARK} in Algorithm \ref{alg_afem_topopt} involves comparing two estimators $\eta_{k,1}$ and $\eta_{k,2}$, the convergent subsequence $\{(\rho^\ast_{k},\bold{u}^\ast_{k})\}_{k\geq 0}$ consists of two subsequences (relabeled as) $\{(\rho_k^{\ast1},\bold{u}_k^{\ast1})\}_{k\geq0}$ and $\{(\rho_k^{\ast2},\bold{u}_k^{\ast2})\}_{k\geq0}$ with
\begin{equation}\label{pf:err->0_01}
\eta^2_{k,1}(\rho_k^{\ast1},\bold{u}_k^{\ast1})\geq \eta^2_{k,2}(\rho_k^{\ast1},\bold{u}_k^{\ast1}),
\end{equation}
\begin{equation}\label{pf:err->0_02}
\eta^2_{k,1}(\rho_k^{\ast2},\bold{u}_k^{\ast2})< \eta^2_{k,2}(\rho_k^{\ast2},\bold{u}_k^{\ast2}).
\end{equation}
Next we divide the lengthy proof into three steps.

\noindent {\bf Step 1.}
For any $T\in\cT_{k+1}$ and any $\delta>0$, Young's inequality implies
\begin{equation}\label{pf:err->0_03}
\begin{aligned}
\eta_{k+1,1}^2(\rho_{k+1}^{\ast1},\bold{u}_{k+1}^{\ast1},T)
&\leq (1+\delta)\eta_{k,1}^2(\rho_{k}^{\ast1},\bold{u}_{k}^{\ast1},T)
+2(1+\delta^{-1})(  h_T^2(\tfrac{\widetilde{\beta}}{\gamma})^2\|W'(\rho^{\ast1}_{k+1})-W'(\rho^{\ast1}_{k})\|^2_{L^2(T)}\\
& \quad +h_T^2p^2\|(\rho_{k+1}^{\ast1})^{p-1}\C_0\strain(\bold{u}_{k+1}^{\ast1}) :\strain(\bold{u}_{k+1}^{\ast1})-(\rho_{k}^{\ast1})^{p-1}\C_0\strain(\bold{u}_{k}^{\ast1}) :\strain(\bold{u}_{k}^{\ast1})\|_{L^2(T)}^2) \\
&\quad  + (1+\delta^{-1})h_T\sum_{F\subset\partial T}\|J_{F,1}(\rho_{k+1}^{\ast1})-J_{F,1}(\rho_{k}^{\ast1})\|^2_{L^2(F)}.
\end{aligned}
\end{equation}
Since $p>1$ is an integer, $(\rho_{k+1}^{\ast1})^{p-1}\C_0\strain(\bold{u}_{k+1}^{\ast1}) :\strain(\bold{u}_{k+1}^{\ast1})-(\rho_{k}^{\ast1})^{p-1}\C_0\strain(\bold{u}_{k}^{\ast1}) :\strain(\bold{u}_{k}^{\ast1})$ is a polynomial of degree $(p-1)$ over $T$. By the element inverse estimate for finite element functions,
\begin{equation}\label{eqn:poly-p-inverse}
\begin{aligned}
&\quad h_T^2\|(\rho_{k+1}^{\ast1})^{p-1}\C_0\strain(\bold{u}_{k+1}^{\ast1}) :\strain(\bold{u}_{k+1}^{\ast1})-(\rho_{k}^{\ast1})^{p-1}\C_0\strain(\bold{u}_{k}^{\ast1}) :\strain(\bold{u}_{k}^{\ast1})\|_{L^2(T)}^2 \\
&\leq
c \|(\rho_{k+1}^{\ast1})^{p-1}\C_0\strain(\bold{u}_{k+1}^{\ast1}) :\strain(\bold{u}_{k+1}^{\ast1})-(\rho_{k}^{\ast1})^{p-1}\C_0\strain(\bold{u}_{k}^{\ast1}) :\strain(\bold{u}_{k}^{\ast1})\|_{L^1(T)}^2\\
& \leq
c \|(\rho_{k+1}^{\ast1})^{p-1}\C_0\strain(\bold{u}_{k+1}^{\ast1}) :\strain(\bold{u}_{k+1}^{\ast1})-(\rho_{k}^{\ast1})^{p-1}\C_0\strain(\bold{u}_{k}^{\ast1}) :\strain(\bold{u}_{k}^{\ast1})\|_{L^1(T)} \\
&\quad \times
\|(\rho_{k+1}^{\ast1})^{p-1}\C_0\strain(\bold{u}_{k+1}^{\ast1}) :\strain(\bold{u}_{k+1}^{\ast1})-(\rho_{k}^{\ast1})^{p-1}\C_0\strain(\bold{u}_{k}^{\ast1}) :\strain(\bold{u}_{k}^{\ast1})\|_{L^1(\Omega)},
\end{aligned}
\end{equation}
and by the scaled trace theorem \cite{Verfurth:2013} and the element inverse estimate, we deduce
\[
h_T\sum_{F\subset\partial T}\|J_{F,1}(\rho_{k+1}^{\ast1})-J_{F,1}(\rho_{k}^{\ast1})\|^2_{L^2(F)}
\leq c \sum_{F\subset\partial T} \|\bold{\nabla}(\rho_{k+1}^{\ast1}-\rho_{k}^{\ast1})\|_{\bold{L}^2(\omega_F)}^2.
\]
Summing the estimate \eqref{pf:err->0_03} over all elements $T\in\cT_{k+1}$ and using the last two estimates and the finite
overlapping property of the neighborhood $\omega_F$ give
\begin{align}
&\eta_{k+1,1}^2(\rho_{k+1}^{\ast1},\bold{u}_{k+1}^{\ast1})
\leq (1+\delta)\sum_{T\in\cT_{k+1}}\eta_{k}^2(\rho_{k}^{\ast1},\bold{u}_{k}^{\ast1},T)\nonumber\\
&\quad + c(1+\delta^{-1})(\|W'(\rho_{k+1}^{\ast1})-W'(\rho_{k}^{\ast1})\|_{L^2(\Omega)}^2+ \|\bold{\nabla}(\rho_{k+1}^{\ast1}-\rho_{k}^{\ast1})\|_{\bold{L}^2(\Omega)}^2 \nonumber\\
&\quad + \|(\rho_{k+1}^{\ast1})^{p-1}\C_0\strain(\bold{u}_{k+1}^{\ast1}) :\strain(\bold{u}_{k+1}^{\ast1})-(\rho_{k}^{\ast1})^{p-1}\C_0\strain(\bold{u}_{k}^{\ast1}) :\strain(\bold{u}_{k}^{\ast1})\|_{L^1(\Omega)}^2 ).\label{pf:err->0_04}
\end{align}
Next we bound the terms separately. By Theorem \ref{thm:medmin}, $W(s)\in C^2[\underline{\rho},1]$ and  Lebesgue's dominated
convergence theorem, we deduce
\begin{equation}\label{pf:err->0_05}
\|W'(\rho^{\ast1}_{k+1})-W'(\rho^{\ast1}_{k})\|_{L^2(\Omega)}^2\rightarrow 0\quad \mbox{and}\quad \|\bold{\nabla}(\rho^{\ast1}_{k+1}-\rho^{\ast1}_{k})\|^2_{\bold{L}^2(\Omega)}\rightarrow 0.
\end{equation}
For the last term, by the triangle inequality, we have
\[
\begin{aligned}
&\quad\|(\rho_{k+1}^{\ast1})^{p-1}\C_0\strain(\bold{u}_{k+1}^{\ast1}) :\strain(\bold{u}_{k+1}^{\ast1})-(\rho_{k}^{\ast1})^{p-1}\C_0\strain(\bold{u}_{k}^{\ast1}) :\strain(\bold{u}_{k}^{\ast1})\|_{L^1(\Omega)}\\ &\leq \|(\rho_{k+1}^{\ast1})^{p-1}(\C_0\strain(\bold{u}_{k+1}^{\ast1}) :\strain(\bold{u}_{k+1}^{\ast1}) - \C_0\strain(\bold{u}_{k}^{\ast1}) :\strain(\bold{u}_{k}^{\ast1}))\|_{L^1(\Omega)}\\
&\quad + \|\left((\rho_{k+1}^{\ast1})^{p-1}-(\rho_{k}^{\ast1})^{p-1}\right)\C_0\strain(\bold{u}_{k}^{\ast1}) :\strain(\bold{u}_{k}^{\ast1})\|_{L^1(\Omega)} \\
&\leq
\|(\rho_{k+1}^{\ast1})^{p-1}(\C_0\strain(\bold{u}_{k+1}^{\ast1}) :\strain(\bold{u}_{k+1}^{\ast1}) - \C_0\strain(\bold{u}_{k}^{\ast1}) :\strain(\bold{u}_{k}^{\ast1}))\|_{L^1(\Omega)}\\
&\quad + \|\left((\rho_{k+1}^{\ast1})^{p-1}-(\rho_{k}^{\ast1})^{p-1}\right)\C_0\strain(\bold{u}_{\infty}^{\ast1}) :\strain(\bold{u}_{\infty}^{\ast1})\|_{L^1(\Omega)} \\
&\quad + \|\left((\rho_{k+1}^{\ast1})^{p-1}-(\rho_{k}^{\ast1})^{p-1}\right)\left(\C_0\strain(\bold{u}_{k}^{\ast1}) :\strain(\bold{u}_{k}^{\ast1})-\C_0\strain(\bold{u}_{\infty}^{\ast1}) :\strain(\bold{u}_{\infty}^{\ast1})\right)\|_{L^1(\Omega)}.
\end{aligned}
\]
It follows from Theorem \ref{thm:medmin} that the sequence $\{\|\bold{u}_{k}^{\ast1}\|_{\bold{H}^1(\Omega)}\}_{k\geq0}$ is
uniformly bounded. Then by the symmetry of the tensor $\C_0$ and the Cauchy-Schwarz inequality,
\[
\begin{aligned}
& \|(\rho_{k+1}^{\ast1})^{p-1}(\C_0\strain(\bold{u}_{k+1}^{\ast1}) :\strain(\bold{u}_{k+1}^{\ast1}) - \C_0\strain(\bold{u}_{k}^{\ast1}) :\strain(\bold{u}_{k}^{\ast1}))\|_{L^1(\Omega)} \\
\leq&
\|(\C_0\strain(\bold{u}_{k+1}^{\ast1})-\C_0\strain(\bold{u}_{k}^{\ast1})):(\strain(\bold{u}_{k+1}^{\ast1})+\strain(\bold{u}_{k}^{\ast1}))\|_{L^1(\Omega)}\\
\leq&
\|\C_0\strain(\bold{u}_{k+1}^{\ast1}) - \C_0\strain(\bold{u}_{k}^{\ast1})\|_{\bold{L}^2(\Omega)} \|\strain(\bold{u}_{k+1}^{\ast1})+\strain(\bold{u}_{k}^{\ast1})\|_{\bold{L}^2(\Omega)}\\
\leq& c \|\bold{u}_{k+1}^{\ast1}-\bold{u}_{k}^{\ast1}\|_{\bold{H}^1(\Omega)}\to 0.
\end{aligned}
\]
Similarly,
\[
\|\left((\rho_{k+1}^{\ast1})^{p-1}-(\rho_{k}^{\ast1})^{p-1}\right)\left(\C_0\strain(\bold{u}_{k}^{\ast1}) :\strain(\bold{u}_{k}^{\ast1})-\C_0\strain(\bold{u}_{\infty}^{\ast1}) :\strain(\bold{u}_{\infty}^{\ast1})\right)\|_{L^1(\Omega)}\to 0.
\]
Upon noting the pointwise convergence of $\{\rho_{k}^\ast\}_{k\geq0}$ in Theorem \ref{thm:medmin} and $\C_0\strain(\bold{u}_{\infty}^{\ast1}) :\strain(\bold{u}_{\infty}^{\ast1}) \in L^1(\Omega)$, and applying Lebesgue's dominated convergence theorem, we find
\begin{equation*}
\|\left((\rho_{k+1}^{\ast1})^{p-1}-(\rho_{k}^{\ast1})^{p-1}\right)\C_0\strain(\bold{u}_{\infty}^{\ast1}) :\strain(\bold{u}_{\infty}^{\ast1})\|_{L^1(\Omega)}\to 0.
\end{equation*}
The last three vanishing limits imply
\begin{equation}\label{pf:err->0_06}
\|(\rho_{k+1}^{\ast1})^{p-1}\C_0\strain(\bold{u}_{k+1}^{\ast1}) :\strain(\bold{u}_{k+1}^{\ast1})-(\rho_{k}^{\ast1})^{p-1}\C_0\strain(\bold{u}_{k}^{\ast1}) :\strain(\bold{u}_{k}^{\ast1})\|_{L^1(\Omega)}^2 \to 0.
\end{equation}
Meanwhile, any $T\in\cT_{k+1}\setminus\cT_{k}$ is generated by applying bisection to some $T'\in\cT_k\setminus\cT_{k+1}$
at least once. Thus $h_T=\sqrt{|T|} \leq \sqrt{|T'|/2} = h_{T'}/\sqrt{2}$ and $J_{F,1}(\rho_{k}^{\ast1})$ is zero across a side
$F$ in the interior of $T'$. Therefore, we can split the term $\sum_{T\in\cT_{k+1}}\eta_{k,1}^2(\bold{u}_{k}^{\ast1},\rho_{k}^{\ast1},T)$ in \eqref{pf:err->0_04} as
\[
\begin{aligned}
\sum_{T\in\cT_{k+1}}\eta_{k,1}^2(\bold{u}_{k}^{\ast1},\rho_{k}^{\ast1},T)
&= \eta_{k,1}^2(\bold{u}_{k}^{\ast1},\rho_{k}^{\ast1},\cT_{k+1}\cap\cT_{k})
+\eta_{k,1}^2(\bold{u}_{k}^{\ast1},\rho_{k}^{\ast1},\cT_{k+1}\setminus\cT_{k})\\
&\leq \eta_{k,1}^2(\bold{u}_{k}^{\ast1},\rho_{k}^{\ast1},\cT_{k+1}\cap\cT_{k})
+ \tfrac{1}{\sqrt{2}}\eta_{k,1}^2(\bold{u}_{k}^{\ast1},\rho_{k}^{\ast1},\cT_{k}\setminus\cT_{k+1})\\
&= \eta_{k,1}^2(\bold{u}_{k}^{\ast1},\rho_{k}^{\ast1},\cT_{k})
-(1- \tfrac{1}{\sqrt{2}})\eta_{k,1}^2(\bold{u}_{k}^{\ast1},\rho_{k}^{\ast1},\cT_{k}\setminus\cT_{k+1}),
\end{aligned}
\]
which, along with $\mathcal{M}_k \subseteq \cT_{k}\setminus\cT_{k+1}$ and D\"{o}rfler marking criterion \eqref{eqn:marking1}, yields
\begin{equation}\label{pf:err->0_07}
\sum_{T\in\cT_{k+1}}\eta_{k,1}^2(\bold{u}_{k}^{\ast1},\rho_{k}^{\ast1},T)
\leq \big(1-\theta_1(1-\tfrac{1}{\sqrt{2}})\big) \eta_{k,1}^2(\bold{u}_{k}^{\ast1},\rho_{k}^{\ast1}).
\end{equation}
Now collecting the estimates \eqref{pf:err->0_04}--\eqref{pf:err->0_07} leads to
\begin{equation*}
\eta_{k+1,1}^2(\bold{u}_{k+1}^{\ast1},\rho_{k+1}^{\ast1})
\leq \alpha_1\eta_{k,1}^2(\bold{u}_{k}^{\ast1},\rho_{k}^{\ast1}) + \xi_{k,1},
\end{equation*}
where for a sufficiently small $\delta>0$,
\[
\alpha_1=(1+\delta)\left(1-\theta_1(1- \tfrac{1}{\sqrt{2}})\right)\in(0,1)
\]
and
\[
\begin{aligned}
\xi_{k,1}:&=     c(1+\delta^{-1})(\|W'(\rho_{k+1}^{\ast1})-W'(\rho_{k}^{\ast1})\|_{L^2(\Omega)}^2+ \|\bold{\nabla}(\rho_{k+1}^{\ast1}-\rho_{k}^{\ast1})\|_{L^2(\Omega)}^2 \\
&\quad + \|(\rho_{k+1}^{\ast1})^{p-1}\C_0\strain(\bold{u}_{k+1}^{\ast1}) :\strain(\bold{u}_{k+1}^{\ast1})-(\rho_{k}^{\ast1})^{p-1}\C_0\strain(\bold{u}_{k}^{\ast1}) :\strain(\bold{u}_{k}^{\ast1})\|_{L^1(\Omega)}^2 ) \to 0.
\end{aligned}
\]
Therefore, one can prove that (cf. \cite[Lemma 2.3]{Aurada:2012})
$$\displaystyle\lim_{k\to\infty}\eta_{k,1}(\rho_{k}^{\ast1},\bold{u}_{k}^{\ast1})=0.$$
Then due to the condition \eqref{pf:err->0_01}, we have
$$\displaystyle\lim_{k\to\infty}\eta_{k,2}(\rho_{k}^{\ast1},\bold{u}_{k}^{\ast1})=0,$$
i.e., the desired vanishing limits hold for $\{(\rho_k^{\ast1},\bold{u}_k^{\ast1})\}_{k\geq0}$.

\noindent {\bf Step 2.} For the sequence $\{(\rho_k^{\ast2},\bold{u}_k^{\ast2})\}_{k\geq0}$, by Young's inequality, we derive that for any $T\in\cT_{k+1}$ and any $\delta>0$
\begin{align*}
\eta_{k+1,2}^2(\rho_{k+1}^{\ast2},\bold{u}_{k+1}^{\ast2},T)
\leq& (1+\delta)\eta_{k,2}^2(\rho_{k}^{\ast2},\bold{u}_{k}^{\ast2},T)
 +(1+\delta^{-1})(  h_T^2\|\bold{\nabla}\cdot(\rho_{k+1}^{\ast2})^{p}\C_0\strain(\bold{u}_{k+1}^{\ast2}) -\bold{\nabla}\cdot(\rho_{k}^{\ast2})^{p}\C_0\strain(\bold{u}_{k}^{\ast2})\|_{\bold{L}^2(T)}^2 \\
& + h_T\sum_{F\subset\partial T}\|[(\rho_{k+1}^{\ast2})^{p}\C_0\strain(\bold{u}_{k+1}^{\ast2}) \cdot\bold{n}_F-(\rho_{k}^{\ast2})^{p}\C_0\strain(\bold{u}_{k}^{\ast2})\cdot\bold{n}_F ]\|^2_{\bold{L}^2(F)}).
\end{align*}
Since each component of $\C_0 \strain(\bold{u}^\ast_k)$ and $\C_0 \strain(\bold{u}^\ast_{k+1})$ is a constant and $(\rho_{k+1}^{\ast2})^{p-1} \C_0 \strain(\bold{u}^{\ast2}_{k+1}) \cdot \bold{\nabla} \rho_{k+1}^{\ast2}-(\rho_k^{\ast2})^{p-1} \C_0 \strain(\bold{u}^{\ast2}_k) \cdot \bold{\nabla} \rho_k^{\ast2}$ is also a polynomial of degree $p-1$ on any $T\in\cT_{k+1}$, as in Step 1, we obtain
\begin{align*}
&\quad h_T^2\|\bold{\nabla}\cdot(\rho_{k+1}^{\ast2})^{p}\C_0\strain(\bold{u}_{k+1}^{\ast2}) -\bold{\nabla}\cdot(\rho_{k}^{\ast2})^{p}\C_0\strain(\bold{u}_{k}^{\ast2})\|_{\bold{L}^2(T)}^2 \\
&=h_T^{2}p^2\| (\rho_{k+1}^{\ast2})^{p-1} \C_0 \strain(\bold{u}^{\ast2}_{k+1}) \cdot \bold{\nabla} \rho_{k+1}^{\ast2}-(\rho_k^{\ast2})^{p-1} \C_0 \strain(\bold{u}^{\ast2}_k) \cdot \bold{\nabla} \rho_k^{\ast2}\|_{\bold{L}^2(T)}^2 \\
&\leq c\| (\rho_{k+1}^{\ast2})^{p-1} \C_0 \strain(\bold{u}^{\ast2}_{k+1}) \cdot \bold{\nabla} \rho_{k+1}^{\ast2}-(\rho_k^{\ast2})^{p-1} \C_0 \strain(\bold{u}^{\ast2}_k) \cdot \bold{\nabla} \rho_k^{\ast2}\|_{\bold{L}^1(T)} \\
&\quad  \times \| (\rho_{k+1}^{\ast2})^{p-1} \C_0 \strain(\bold{u}^{\ast2}_{k+1}) \cdot \bold{\nabla} \rho_{k+1}^{\ast2}-(\rho_k^{\ast2})^{p-1} \C_0 \strain(\bold{u}^{\ast2}_k) \cdot \bold{\nabla} \rho_k^{\ast2}\|_{\bold{L}^1(\Omega)},
\end{align*}
and by the scaled trace theorem \cite{Verfurth:2013} and the element inverse estimate,
\begin{align*}
&h_T\sum_{F\subset\partial T}\|[(\rho_{k+1}^{\ast2})^{p}\C_0\strain(\bold{u}_{k+1}^{\ast2}) \cdot\bold{n}_F-(\rho_{k}^{\ast2})^{p}\C_0\strain(\bold{u}_{k}^{\ast2})\cdot\bold{n}_F ]\|^2_{\bold{L}^2(F)}\\
\leq& c \sum_{F\subset\partial T} \|(\rho_{k+1}^{\ast2})^{p}\C_0\strain(\bold{u}_{k+1}^{\ast2})-(\rho_{k}^{\ast2})^{p}\C_0\strain(\bold{u}_{k}^{\ast2})\|_{\bold{L}^2(\omega_F)}^2.
\end{align*}
Consequently,
\begin{align}
\eta_{k+1,2}^2(\rho_{k+1}^{\ast2},\bold{u}_{k+1}^{\ast2})
&\leq
(1+\delta)\sum_{T\in\cT_{k+1}}\eta_{k,2}^2(\rho_{k}^{\ast2},\bold{u}_{k}^{\ast2},T)
+ c(1+\delta^{-1})(\|(\rho_{k+1}^{\ast2})^{p}\C_0\strain(\bold{u}_{k+1}^{\ast2})-(\rho_{k}^{\ast2})^{p}\C_0\strain(\bold{u}_{k}^{\ast2})\|_{\bold{L}^2(\Omega)}^2\nonumber \\
&\quad +\| (\rho_{k+1}^{\ast2})^{p-1} \C_0 \strain(\bold{u}^{\ast2}_{k+1}) \cdot \bold{\nabla} \rho_{k+1}^{\ast2}-(\rho_k^{\ast2})^{p-1} \C_0 \strain(\bold{u}^{\ast2}_k) \cdot \bold{\nabla} \rho_k^{\ast2}\|_{\bold{L}^1(\Omega)}^2).\label{pf:err->0_08}
\end{align}
By D\"{o}rfler marking criterion \eqref{eqn:marking2} and repeating the argument for \eqref{pf:err->0_07}, we obtain
\begin{equation}\label{pf:err->0_09}
\sum_{T\in\cT_{k+1}}\eta_{k,2}^2(\rho_{k}^{\ast2},\bold{u}_{k}^{\ast2},T)
\leq \left(1-\theta_2(1-\tfrac{1}{\sqrt{2}})\right) \eta_{k,2}^2(\bold{u}_{k}^{\ast2},\rho_{k}^{\ast2}).
\end{equation}
By the triangle inequality,
\begin{align*}
&\quad \| (\rho_{k+1}^{\ast2})^{p-1} \C_0 \strain(\bold{u}^{\ast2}_{k+1}) \cdot \bold{\nabla} \rho_{k+1}^{\ast2}-(\rho_k^{\ast2})^{p-1} \C_0 \strain(\bold{u}^{\ast2}_k) \cdot \bold{\nabla} \rho_k^{\ast2}\|_{\bold{L}^1(\Omega)} \\
&\leq \|(\rho_{k+1}^{\ast2})^{p-1}\left(\C_0 \strain(\bold{u}^{\ast2}_{k+1}) \cdot \bold{\nabla} \rho_{k+1}^{\ast2}-\C_0 \strain(\bold{u}^{\ast2}_{k}) \cdot \bold{\nabla} \rho_{k}^{\ast2}\right)\|_{\bold{L}^1(\Omega)}\\
&\quad+ \|\left((\rho_{k+1}^{\ast2})^{p-1}-(\rho_{k}^{\ast2})^{p-1}\right)\C_0 \strain(\bold{u}^{\ast2}_{k}) \cdot \bold{\nabla} \rho_{k}^{\ast2}\|_{\bold{L}^1(\Omega)}\\
&\leq \|(\rho_{k+1}^{\ast2})^{p-1}\left(\C_0 \strain(\bold{u}^{\ast2}_{k+1}) \cdot \bold{\nabla} \rho_{k+1}^{\ast2}-\C_0 \strain(\bold{u}^{\ast2}_{k}) \cdot \bold{\nabla} \rho_{k}^{\ast2}\right)\|_{\bold{L}^1(\Omega)}\\
&\quad + \|\left((\rho_{k+1}^{\ast2})^{p-1}-(\rho_{k}^{\ast2})^{p-1}\right)\left(\C_0 \strain(\bold{u}^{\ast2}_{k}) \cdot \bold{\nabla} \rho_{k}^{\ast2}-\C_0 \strain(\bold{u}^{\ast2}_{\infty}) \cdot \bold{\nabla} \rho_{\infty}^{\ast2}\right)\|_{\bold{L}^1(\Omega)} \\
&\quad + \|\left((\rho_{k+1}^{\ast2})^{p-1}-(\rho_{k}^{\ast2})^{p-1}\right)\C_0 \strain(\bold{u}^{\ast2}_{\infty}) \cdot \bold{\nabla} \rho_{\infty}^{\ast2}\|_{\bold{L}^1(\Omega)}.
\end{align*}
Since the sequences $\{\|\bold{u}_{k}^\ast\|_{\bold{H}^1(\Omega)}\}_{k\geq 0}$ and $\{\|\rho_k^\ast\|_{H^1(\Omega)}\}_{k\geq0}$ are both uniformly
bounded due to Theorem \ref{thm:medmin}, the Cauchy-Schwarz inequality implies
\[
\begin{aligned}
&\quad \|(\rho_{k+1}^{\ast2})^{p-1}\left(\C_0 \strain(\bold{u}^{\ast2}_{k+1}) \cdot \bold{\nabla} \rho_{k+1}^{\ast2}-\C_0 \strain(\bold{u}^{\ast2}_{k}) \cdot \bold{\nabla} \rho_{k}^{\ast2}\right)\|_{\bold{L}^1(\Omega)}\\
&\leq
\|(\rho_{k+1}^{\ast2})^{p-1}\C_0 \strain(\bold{u}^{\ast2}_{k+1}) \cdot \bold{\nabla} \left(\rho_{k+1}^{\ast2} -\rho_{k}^{\ast2}\right)\|_{\bold{L}^1(\Omega)}+ \|(\rho_{k+1}^{\ast2})^{p-1}\left(\C_0 \strain(\bold{u}^{\ast2}_{k+1})-\C_0 \strain(\bold{u}^{\ast2}_{k})\right) \cdot \bold{\nabla} \rho_{k}^{\ast2}\|_{\bold{L}^1(\Omega)}\\
&\leq c\left(\|\rho^{\ast2}_{k+1}-\rho^{\ast2}_{k}\|_{H^1(\Omega)}+\|\bold{u}_{k+1}^{\ast1}-\bold{u}_{k}^{\ast1}\|_{\bold{H}^1(\Omega)}\right)\rightarrow 0.
\end{aligned}
\]
Likewise,
\begin{align*}
&\|\left((\rho_{k+1}^{\ast2})^{p-1}-(\rho_{k}^{\ast2})^{p-1}\right)\left(\C_0 \strain(\bold{u}^{\ast2}_{k}) \cdot \bold{\nabla} \rho_{k}^{\ast2}-\C_0 \strain(\bold{u}^{\ast2}_{\infty})  \cdot \bold{\nabla} \rho_{\infty}^{\ast2}\right)\|_{\bold{L}^1(\Omega)}\\
\leq &\|\C_0 \strain(\bold{u}^{\ast2}_{k}) \cdot \bold{\nabla} \rho_{k}^{\ast2}-\C_0 \strain(\bold{u}^{\ast2}_{\infty}) \cdot \bold{\nabla} \rho_{\infty}^{\ast2}\|_{\bold{L}^1(\Omega)}\to 0.
\end{align*}
By appealing to Theorem \ref{thm:medmin} again, we deduce
\[
\left((\rho_{k+1}^{\ast2})^{p-1}-(\rho_{k}^{\ast2})^{p-1}\right)\C_0 \strain(\bold{u}^{\ast2}_{\infty}) \cdot \bold{\nabla} \rho_{\infty}^{\ast2}\to 0 \quad\text{a.e. in}~\Omega,
\]
\[ |\left((\rho_{k+1}^{\ast2})^{p-1}-(\rho_{k}^{\ast2})^{p-1}\right)\C_0 \strain(\bold{u}^{\ast2}_{\infty})  \cdot \bold{\nabla} \rho_{\infty}^{\ast2}|\leq|\C_0 \strain(\bold{u}^{\ast2}_{\infty}) \cdot \bold{\nabla} \rho_{\infty}^{\ast2}|\in L^1(\Omega).
\]
Then Lebesgue's dominated convergence theorem implies
\[
\|\left((\rho_{k+1}^{\ast2})^{p-1}-(\rho_{k}^{\ast2})^{p-1}\right)\C_0 \strain(\bold{u}^{\ast2}_{\infty}) \cdot \bold{\nabla} \rho_{\infty}^{\ast2}\|_{\bold{L}^1(\Omega)}\to 0.
\]
Therefore, we obtain
\begin{equation}\label{pf:err->0_10}
\| (\rho_{k+1}^{\ast2})^{p-1} \C_0 \strain(\bold{u}^{\ast2}_{k+1}) \bold{\nabla} \rho_{k+1}^{\ast2}-(\rho_k^{\ast2})^{p-1} \C_0 \strain(\bold{u}^{\ast2}_k) \cdot \bold{\nabla} \rho_k^{\ast2}\|_{\bold{L}^1(\Omega)}\to 0.
\end{equation}
Repeating the preceding argument yields
\begin{align}
&\quad\|(\rho_{k+1}^{\ast2})^{p}\C_0\strain(\bold{u}_{k+1}^{\ast2})-(\rho_{k}^{\ast2})^{p}\C_0\strain(\bold{u}_{k}^{\ast2})\|_{\bold{L}^2(\Omega)}\nonumber\\
&\leq \|(\rho_{k+1}^{\ast2})^{p}\left(\C_0\strain(\bold{u}_{k+1}^{\ast2})-\C_0\strain(\bold{u}_{k}^{\ast2})\right)\|_{\bold{L}^2(\Omega)}+
\|\left((\rho_{k+1}^{\ast2})^{p}-(\rho_{k}^{\ast2})^{p}\right)\C_0\strain(\bold{u}_{k}^{\ast2})\|_{\bold{L}^2(\Omega)}\nonumber\\
&\leq \|(\rho_{k+1}^{\ast2})^{p}\left(\C_0\strain(\bold{u}_{k+1}^{\ast2})-\C_0\strain(\bold{u}_{k}^{\ast2})\right)\|_{\bold{L}^2(\Omega)}+
\|\left((\rho_{k+1}^{\ast2})^{p}-(\rho_{k}^{\ast2})^{p}\right)\left(\C_0\strain(\bold{u}_{k}^{\ast2})-\C_0\strain(\bold{u}_{\infty}^{\ast2})\right)\|_{\bold{L}^2(\Omega)}\nonumber\\
&\quad + \|\left((\rho_{k+1}^{\ast2})^{p}-(\rho_{k}^{\ast2})^{p}\right)\C_0\strain(\bold{u}_{\infty}^{\ast2})\|_{\bold{L}^2(\Omega)}\to 0.\label{pf:err->0_11}
\end{align}
It follows from \eqref{pf:err->0_08}--\eqref{pf:err->0_11} that
\[
\eta_{k+1,2}^2(\bold{u}_{k+1}^{\ast2},\rho_{k+1}^{\ast2})
\leq \alpha_2\eta_{k,2}^2(\bold{u}_{k}^{\ast2},\rho_{k}^{\ast2}) + \xi_{k,2},
\]
where for sufficiently small $\delta>0$,
\[
\alpha_2=(1+\delta)(1-\theta_2(1- \tfrac{1}{\sqrt{2}}))\in(0,1)
\]
and
\[
\begin{aligned}
\xi_{k,2}:&=     c(1+\delta^{-1})(\|(\rho_{k+1}^{\ast2})^{p}\C_0\strain(\bold{u}_{k+1}^{\ast2})-(\rho_{k}^{\ast2})^{p}\C_0\strain(\bold{u}_{k}^{\ast2})\|_{\bold{L}^2(\Omega)}^2 \\
&\quad+\| (\rho_{k+1}^{\ast2})^{p-1} \C_0 \strain(\bold{u}^{\ast2}_{k+1}) \bold{\nabla} \rho_{k+1}^{\ast2}-(\rho_k^{\ast2})^{p-1} \C_0 \strain(\bold{u}^{\ast2}_k) \bold{\nabla} \rho_k^{\ast2}\|_{\bold{L}^1(\Omega)}^2)\to 0.
\end{aligned}
\]
Like before, we obtain
$$\lim_{k\to\infty}\eta_{k,2}(\rho_{k}^{\ast2},\bold{u}_{k}^{\ast2})=0.$$
Then by the condition \eqref{pf:err->0_02}, we have
$$\displaystyle\lim_{k\to\infty}\eta_{k,1}(\rho_{k}^{\ast2},\bold{u}_{k}^{\ast2})=0.$$

\noindent {\bf Step 3.} Since $\{(\rho_{k}^{\ast1},\bold{u}_{k}^{\ast1})\}_{k\geq0}$ and  $\{(\rho_{k}^{\ast2},\bold{u}_{k}^{\ast2})\}_{k\geq0}$ comprise $\{(\rho_{k}^{\ast},\bold{u}_{k}^{\ast})\}_{k\geq0}$, the vanishing limits obtained in Steps 1 and 2 for the two subsequences ensure the desired identity \eqref{err->0}.
\end{proof}

\begin{remark}
{The proof of Lemma \ref{lem:err->0}, to be precise, the derivation of the estimates \eqref{pf:err->0_04} and \eqref{pf:err->0_08}, requires that the exponent $p$ be an integer greater than $1$. This is due to the crucial use of the element inverse inequality for polynomials; see, e.g., \eqref{eqn:poly-p-inverse}. The proof technique does not apply directly to the case of a non-integer $p$.}
\end{remark}

The next result gives the asymptotic vanishing property (of the residual operators).
\begin{lemma}\label{lem:residual}
If $p>1$ is an integer, then for the convergent
subsequence $\{(\rho^{\ast}_{k_j},\bold{u}^{\ast}_{k_j})\}_{j\geq0}$ in Theorem \ref{thm:medmin}, there hold
\begin{equation}\label{residual01}
\liminf_{j\to\infty}\left(\widetilde{\beta}\mathcal{F}_\gamma'(\rho^\ast_{k_j})(\phi-\rho^\ast_{k_j})
- \int_{\Omega} p (\rho^\ast_{k_j})^{p-1} \mathcal{C}_0 \strain (\bold{u}^\ast_{k_j}) :
\strain (\bold{u}^\ast_{k_j}) (\phi -\rho^\ast_{k_j})  \dx \right) \geq 0, \quad
\forall\phi \in \Atilde,
\end{equation}
\begin{equation}\label{residual02}
\lim_{j\to\infty} \left(\int_\Omega (\rho_{k_j}^\ast)^p\C_0 \strain(\bold{u}^\ast_{k_j}):
\strain(\bold{v}) \dd x - \int_\Omega \bold{f} \cdot \bold{v} \dd x - \int_{\Gamma_N}
\bold{g}\cdot \bold{v} \dd s \right)= 0,\quad  \forall\bold{v} \in  \bold{V}^{D},
\end{equation}
with
$$\mathcal{F}_\gamma'(\rho^\ast_{k_j})(\phi-\rho^\ast_{k_j})=\gamma\int_\Omega
\bold{\nabla} \rho^\ast_{k_j} \cdot \bold{\nabla} (\phi - \rho^\ast_{k_j}) \dx
+ \dfrac{1}{\gamma}\int_{\Omega} W'(\rho^\ast_{k_j})(\phi - \rho^\ast_{k_j})
\dx.$$
\end{lemma}

To prove Lemma \ref{lem:residual}, we need a constraint-preserving interpolation operator for
$L^1(\Omega)$-functions introduced in \cite{JinXu:2019}. Let $ \mathcal{N}_k$ be the set of all
nodes of $\cT_k$ and $\{\phi_{x}\}_{x\in\mathcal{N}_k} $ be the nodal basis functions in
$V_k$. For each $x\in\mathcal{N}_k$, the support of $\phi_x$ is denoted by $\omega_x$, i.e., the union of
all elements in $\cT_k$ with $x$ being a vertex, and $\cT_k(\omega_x)$ is the triangulation of $\omega_x$
with respect to $\cT_k$. Then we define an interpolation operator $\Pi_k:L^1(\Om)\to V_k$ by
\begin{equation}\label{eqn:cn_int_def}
\Pi_k v := \sum_{x\in\mathcal{N}_k} \frac{1}{|\omega_x|}\int_{\omega_x}v \dx \phi_x.
\end{equation}
It follows from the definition \eqref{eqn:cn_int_def} that
\begin{equation}\label{eqn:cn_int_prop1}
\underline{\rho} \leq  \Pi_k v \leq 1,\quad \text{if}~ \underline{\rho} \leq  v \leq 1~\text{a.e. in}~\Omega.
\end{equation}
Further, by the standard trapezoidal rule,
\[
\begin{aligned}
\int_\Omega \Pi_k v \dx & = \sum_{x\in\mathcal{N}_k} \frac{1}{|\omega_x|}\int_{\omega_x} v \dx \int_{\omega_x} \phi_x \dx = \sum_{x\in\mathcal{N}_k} \frac{1}{|\omega_x|}\int_{\omega_x} v \dx \sum_{T\in \cT_k(\omega_x)}\int_{T} \phi_x \dx \\
&= \sum_{x\in\mathcal{N}_k} \frac{1}{|\omega_x|}\int_{\omega_x} v \dx \sum_{T\in \cT_k(\omega_x)} \frac{|T|}{3}= \frac{1}{3} \sum_{x\in\mathcal{N}_k} \sum_{T\in \cT_k(\omega_x)} \int_T v \dx = \int_\Omega v \dx.
\end{aligned}
\]
So $\Pi_k$ preserves the integral constraint in $\Atilde$, i.e.,
\begin{equation}\label{eqn:cn_int_prop2}
\int_\Omega \Pi_k v \dx = V^0,\quad\text{if}~\int_\Omega v \dx = V^0.
\end{equation}
Moreover, there hold the following stability and error estimates \cite[Lemma 5.3]{JinXu:2019}.

\begin{lemma}\label{lem:err-NZ}
For all $v\in H^1(\Omega)$ on all $T\in\cT_k$ and any $F\in\partial T \cap \mathcal{F}_k$, there hold
\begin{equation}\label{int_est1}
\|\Pi_k v\|_{L^r(T)} \leq c \|v\|_{L^2(D_T)}, \quad
\|\nabla\Pi_k v\|_{L^r(T)} \leq c \|\nabla v\|_{L^r(D_T)},
\end{equation}
\begin{equation}\label{int_est2}
\| v - \Pi_k v \|_{L^r(T)} \leq c h_T \| \nabla v \|_{L^2(D_T)}, \quad\| v - \Pi_k v\|_{L^2(F)} \leq c h_T^{1/2} \|\nabla v\|_{L^2(D_T)}.
\end{equation}
\end{lemma}

\begin{proof}[Proof of Lemma \ref{lem:residual}]
Like before, we relabel the convergent subsequence $\{(\rho^{\ast}_{k_j},\bold{u}^{\ast}_{k_j})\}_{j\geq0}$ as $\{(\rho^{\ast}_{k},\bold{u}^{\ast}_{k})\}_{k\geq0}$. The basic idea of the proof has been described when motivating the two error estimators in Section \ref{sec:adapt}. Below we expand the idea rigorously.
Since $(\rho_k^\ast,\bold{u}_k^\ast)$ is a minimizer to problem
\eqref{dismin}--\eqref{dismin-state} over the triangulation $\cT_k$,
it satisfies the following discrete variational inequality (see the first variational inequality of \eqref{disoptsys}):
\begin{equation}\label{pf:residual_01}
\widetilde{\beta}\mathcal{F}_\gamma'(\rho^\ast_{k})(\phi_k-\rho^\ast_{k})
- \int_{\Omega} p (\rho^\ast_{k})^{p-1} \mathcal{C}_0 \strain (\bold{u}^\ast_{k}) : \strain (\bold{u}^\ast_{k}) (\phi_k -\rho^\ast_{k})  \dx \geq 0,
\quad \forall\phi_k \in \Atilde_k.
\end{equation}
Due to the conservation properties \eqref{eqn:cn_int_prop1}--\eqref{eqn:cn_int_prop2}, there holds
$\Pi_k\phi \in \Atilde_k$ for any $\phi\in \Atilde $.
Letting $\phi_k=\Pi_k\phi$ in \eqref{pf:residual_01} gives that for any $\phi\in \Atilde$,
\begin{align}
&\quad\widetilde{\beta}\mathcal{F}_\gamma'(\rho^\ast_{k})(\phi-\rho^\ast_{k})
- \int_{\Omega} p (\rho^\ast_{k})^{p-1} \mathcal{C}_0 \strain (\bold{u}^\ast_{k}) : \strain
(\bold{u}^\ast_{k}) (\phi -\rho^\ast_{k})  \dx  \nonumber\\
&= \widetilde{\beta} \mathcal{F}_\gamma'(\rho^\ast_{k})(\phi-\Pi_k\phi)
- \int_{\Omega} p (\rho^\ast_{k})^{p-1} \mathcal{C}_0 \strain (\bold{u}^\ast_{k}) : \strain (\bold{u}^\ast_{k}) (\phi -\Pi_k\phi)  \dx \nonumber\\
&\quad + \widetilde{\beta} \mathcal{F}_\gamma'(\rho^\ast_{k})(\Pi_k\phi-\rho^\ast_{k})
- \int_{\Omega} p (\rho^\ast_{k})^{p-1} \mathcal{C}_0 \strain (\bold{u}^\ast_{k}) : \strain (\bold{u}^\ast_{k}) (\Pi_k\phi -\rho^\ast_{k})  \dx \nonumber\\
&\geq \widetilde{\beta} \mathcal{F}_\gamma'(\rho^\ast_{k})(\phi-\Pi_k\phi)
- \int_{\Omega} p (\rho^\ast_{k})^{p-1} \mathcal{C}_0
\strain (\bold{u}^\ast_{k}) : \strain (\bold{u}^\ast_{k}) (\phi -\Pi_k\phi)  \dx. \label{pf:residual1_aux}
\end{align}
By elementwise integration by parts, the Cauchy-Schwarz inequality, the definition of the
error indicator $\eta_{k,1}$ and error estimates for $\Pi_k$ from Lemma \ref{lem:err-NZ}, we get
\[
\begin{aligned}
&\quad\left| \widetilde{\beta} \mathcal{F}_\gamma'(\rho^{\ast}_{k})(\phi-\Pi_k\phi)
- \int_{\Omega} p (\rho^{\ast}_{k})^{p-1} \mathcal{C}_0 \strain (\bold{u}^{\ast}_{k}) : \strain (\bold{u}^{\ast}_{k}) (\phi -\Pi_k\phi)  \dx  \right|\\
&= \left|  \sum_{T \in \cT_k} \int_T R_{T,1}(\rho_{k}^{\ast}, \bold{u}_k^{\ast})(\phi-\Pi_k\phi)\dx + \sum_{F\in\mathcal{F}_k}\int_F J_{F,1}(\rho_k^{\ast})(\phi-\Pi_k\phi) \dd s \right|  \\
& \leq \sum_{T \in \cT_k} \left(\|R_{T,1}(\rho_{k}^{\ast}, \bold{u}_k^{\ast})\|_{L^2(T)}\|\phi-\Pi_k\phi\|_{L^{2}(T)}
+   \sum_{F\subset\partial T}\| J_{F,1}(\rho_k^{\ast})\|_{L^{2}(F)}\|\phi-\Pi_k\phi\|_{L^{2}(F)}\right) \\
& \leq c \sum_{T \in \cT_k} \left(h_T\|R_{T,1}(\rho_{k}^{\ast}, \bold{u}_k^{\ast})\|_{L^2(T)} + h_T^{1/2}\sum_{F\subset\partial T}\| J_{F,1}(\rho_k^{\ast})\|_{L^{2}(F)}\right) \|\bold{\nabla}\phi\|_{L^{2}(D_T)} \\
& \leq  c  \eta_{k,1}(\rho_{k}^{\ast}, \bold{u}_k^{\ast}) \|\bold{\nabla}\phi\|_{L^{2}(\Omega)}.
\end{aligned}
\]
Now the assertion \eqref{residual01} follows from \eqref{pf:residual1_aux} and Lemma \ref{lem:err->0}.
The proof of \eqref{residual02} is similar. On the space $\bold{V}_k^D$, we define a vectorial operator
$\bold{I}_k^{sz}$ with each component being the Scott-Zhang quasi-interpolation \cite{ScottZhang:1990}. Since
$(\rho_k^\ast,\bold{u}_k^\ast)$ solves \eqref{dismin-state}, by applying the operator $\bold{I}_k^{sz}$ to any
$\bold{v}\in \bold{V}^D$, elementwise integration by parts and using the definition of $\eta_{k,2}$, error
estimates for $\bold{I}_k^{sz}$ \cite{ScottZhang:1990} and the Cauchy-Schwarz inequality, we obtain that for any $\bold{v}\in \bold{V}^D $,
\begin{align*}
&\quad \left|\int_\Omega (\rho^\ast_k)^p\C_0 \strain(\bold{u}^\ast_{k}) : \strain(\bold{v}) \dd x - \int_\Omega \bold{f} \cdot \bold{v} \dx - \int_{\Gamma_N} \bold{g}\cdot \bold{v} \dd s\right|  \\
&=\left|\int_\Omega (\rho^\ast_k)^p\C_0 \strain(\bold{u}^\ast_{k}) : \strain(\bold{v}-\bold{I}_k^{sz}\bold{v}) \dd x - \int_\Omega \bold{f} \cdot (\bold{v}-\bold{I}_k^{sz}\bold{v}) \dx - \int_{\Gamma_N} \bold{g}\cdot (\bold{v}-\bold{I}_k^{sz}\bold{v}) \dd s\right| \\
& = \left| \sum_{T\in\cT_k}\int_{T} -\bold{R}_{T,2}(\rho_k^\ast,\bold{u}_k^\ast) \cdot (\bold{v}-\bold{I}_k^{sz}\bold{v})  \dx + \sum_{F\in\mathcal{F}_k^i} \int_F \bold{J}_{F,2} \cdot (\bold{v} - \bold{I}_k^{sz}\bold{v}) \dd s + \sum_{F\in\mathcal{F}_k^N} \int_F \bold{J}_{F,2} \cdot (\bold{v} - \bold{I}_k^{sz}\bold{v}) \dd s \right|\\
& \leq c \sum_{T\in\cT_k}\eta_{k,2}(\rho_k^\ast,\bold{u}_k^\ast,T)\|\bold{\nabla}\bold{v}\|_{\bold{L}^2(D_T)}\leq c\eta_{k,2}(\rho_k^\ast,\bold{u}_k^\ast)\|\bold{\nabla}\bold{v}\|_{\bold{L}^2(\Omega)}.
\end{align*}
Like above, Lemma \ref{lem:err->0} directly implies \eqref{residual02}.
\end{proof}

\begin{remark}
The proof of Lemma \ref{lem:residual} indicates that the left hand sides of \eqref{residual01} and \eqref{residual02} are actually bounded from above by $\eta_{k,1}$ and $\eta_{k,2}$, associated with the convergent subsequence, respectively. This key observation motivates the two \textit{a posteriori} error estimators in the module \texttt{ESTIMATE} of Algorithm \ref{alg_afem_topopt}.
\end{remark}

\begin{lemma}\label{lem:medcont}
If $p>1$ is an integer, then the minimizer $(\rho_\infty^\ast,\bold{u}_\infty^\ast)$ to problem \eqref{medmin}-\eqref{medmin-state} in Theorem \ref{thm:medmin} satisfies
\begin{equation}\label{medoptsys}
\left\{
\begin{array}{ll}
\displaystyle\widetilde{\beta} \mathcal{F}_\gamma'(\rho^\ast_\infty)(\phi-\rho^\ast_\infty) - \int_{\Omega} p (\rho^\ast_\infty)^{p-1} \mathcal{C}_0
\strain (\bold{u}^\ast_\infty) : \strain (\bold{u}^\ast_\infty) (\phi -\rho^\ast_\infty)  \dx  \geq 0, \quad \forall \phi \in \Atilde,    \\[0.2cm]
\displaystyle\int_\Omega (\rho_\infty^\ast)^p\C_0 \strain(\bold{u}^\ast_\infty) : \strain(\bold{v}) \dd x = \int_\Omega \bold{f} \cdot \bold{v} \dd x
+ \int_{\Gamma_N} \bold{g}\cdot \bold{v} \dd x,\quad  \forall\bold{v} \in  \bold{V}^{D}
\end{array}
\right.
\end{equation}
with
$$\mathcal{F}_\gamma'(\rho^\ast_\infty)(\phi-\rho^\ast_\infty)= \gamma\int_\Omega \bold{\nabla}
\rho^\ast_\infty \cdot \bold{\nabla} (\phi - \rho^\ast_\infty) \dx
+ \dfrac{1}{\gamma}\int_{\Omega} W'(\rho^\ast_\infty)(\phi - \rho^\ast_\infty) \dx .$$
\end{lemma}
\begin{proof}
We relabel the convergent subsequence $\{(\rho_{k_j}^\ast,\bold{u}^\ast_{k_j})\}_{j\geq 0}$ in Theorem \ref{thm:medmin}
as $\{(\rho_{k}^\ast,\bold{u}^\ast_{k})\}_{k\geq 0}$. In view of the $H^1(\Omega)$-convergence of $\{\rho_{k}^\ast\}_{k\geq 0}$
from Theorem \ref{thm:medmin}, we have
\begin{equation}\label{pf:medcont_01}
\int_\Omega \bold{\nabla} \rho^\ast_k \cdot \bold{\nabla} (\phi - \rho^\ast_k) \dx
\to \int_\Omega \bold{\nabla} \rho^\ast_\infty \cdot \bold{\nabla} (\phi - \rho^\ast_\infty) \dx, \quad \forall\phi \in \Atilde.
\end{equation}
Since $W'(s)\in C^1[\underline{\rho},1]$, by the pointwise convergence of $\{\rho_{k}^\ast\}_{k\geq0}$ and Lebesgue's
dominated convergence theorem \cite[Theorem 1.19]{EvansGariepy:2015},
\begin{equation}\label{pf:medcont_02}
\int_{\Omega} W'(\rho^\ast_k)(\phi - \rho^\ast_k) \dx
\to    \int_{\Omega} W'(\rho^\ast_\infty)(\phi - \rho^\ast_\infty) \dx,
\quad \forall\phi\in  \Atilde.
\end{equation}
Direct computation yields for any $\phi\in \widetilde A$,
\begin{align*}
&\quad \int_{\Omega} p (\rho^\ast_k)^{p-1} \mathcal{C}_0 \strain (\bold{u}^\ast_k) : \strain (\bold{u}^\ast_k) \phi \dx - \int_{\Omega}
p (\rho^\ast_\infty)^{p-1} \mathcal{C}_0 \strain (\bold{u}^\ast_\infty) : \strain (\bold{u}^\ast_\infty) \phi \dx \\
&= \int_{\Omega} p (\rho^\ast_k)^{p-1}(\mathcal{C}_0 \strain (\bold{u}^\ast_k) : \strain (\bold{u}^\ast_k)- \mathcal{C}_0
\strain (\bold{u}^\ast_\infty) : \strain (\bold{u}^\ast_\infty))\phi \dx\\
&\quad + \int_{\Omega} p \left((\rho^\ast_k)^{p-1}-(\rho^\ast_\infty)^{p-1}\right) \mathcal{C}_0 \strain (\bold{u}^\ast_\infty) : \strain (\bold{u}^\ast_\infty) \phi \dx.
\end{align*}
Due to $|\rho^\ast_k|,|\phi|\leq 1$, the $\boldsymbol{H}^1(\Omega)$-convergence of $\{\bold{u}_k^\ast\}_{k\geq0}$ to
$\bold{u}_\infty^\ast$  in Theorem \ref{thm:medmin} implies
\[
\left|\int_{\Omega} p (\rho^\ast_k)^{p-1}(\mathcal{C}_0 \strain (\bold{u}^\ast_k) : \strain (\bold{u}^\ast_k)- \mathcal{C}_0 \strain (\bold{u}^\ast_\infty) : \strain (\bold{u}^\ast_\infty))\phi \dx \right| \leq c\|\bold{u}^\ast_{k}-\bold{u}^\ast_\infty\|_{\bold{H}^1(\Omega)} \to 0.
\]
Since $|(\rho^\ast_k)^{p-1}-(\rho^\ast_\infty)^{p-1}|\leq 1$, the pointwise convergence of
$\{\rho_{k}^\ast\}_{k\geq0}$ and Lebesgue's dominated convergence theorem \cite[Theorem 1.19]{EvansGariepy:2015} imply
\[
\int_{\Omega} p \left((\rho^\ast_k)^{p-1}-(\rho^\ast_\infty)^{p-1}\right) \mathcal{C}_0 \strain
(\bold{u}^\ast_\infty) : \strain (\bold{u}^\ast_\infty) \phi \dx \to 0.
\]
Collecting the last three estimates leads to
\begin{equation}\label{pf:medcont_03}
\int_{\Omega} p (\rho^\ast_k)^{p-1} \mathcal{C}_0 \strain (\bold{u}^\ast_k) : \strain (\bold{u}^\ast_k)
\phi \dx \to \int_{\Omega} p (\rho^\ast_\infty)^{p-1} \mathcal{C}_0 \strain (\bold{u}^\ast_\infty) :
\strain (\bold{u}^\ast_\infty) \phi \dx,  \quad \forall\phi\in \Atilde.
\end{equation}
By invoking \eqref{dismin-state} with $\rho_\cT=\rho^\ast_k$ and \eqref{medmin-state} with $\rho_\infty=\rho^\ast_\infty$ and arguing as for \eqref{pf:medmin_05} (see the proof of Theorem \ref{thm:medmin}), we further have
\begin{equation}\label{pf:medcont_04}
\begin{aligned}
&\quad\int_\Omega (\rho_{k}^\ast)^p \C_0 \strain(\bold{u}_{k}^\ast):\strain(\bold{u}_{k}^\ast) \dx = \int_{\Omega} \bold{f} \cdot \bold{u}_{k}^\ast \dx
+ \int_{\Gamma_N} \bold{g} \cdot \bold{u}_{k}^\ast \dd s \\
& \to   \int_{\Omega} \bold{f} \cdot \bold{u}_{\infty}^\ast \dx
+ \int_{\Gamma_N} \bold{g} \cdot \bold{u}_{\infty}^\ast \dd s = \int_\Omega (\rho_{\infty}^\ast)^p \C_0 \strain(\bold{u}_{\infty}^\ast):\strain(\bold{u}_{\infty}^\ast) \dx.
\end{aligned}
\end{equation}
Then it follows from \eqref{pf:medcont_01}-\eqref{pf:medcont_04} and \eqref{residual01} in Lemma \ref{lem:residual} that
\[
\widetilde{\beta} \mathcal{F}_\gamma'(\rho^\ast_\infty)(\phi-\rho^\ast_\infty) - \int_{\Omega} p (\rho^\ast_\infty)^{p-1}
\mathcal{C}_0 \strain (\bold{u}^\ast_\infty) : \strain (\bold{u}^\ast_\infty) (\phi -\rho^\ast_\infty)
\dx  \geq 0, \quad \forall\phi \in \Atilde.
\]
For the variational problem in \eqref{medoptsys}, Theorem \ref{thm:medmin} and the dominated convergence imply
\[
\|\left((\rho_{k}^\ast)^p-(\rho_{\infty}^\ast)^p\right)\strain(\bold{v})\|_{L^2(\Omega)^{2\times2}}\to0,\quad
\|\C_0 (\strain(\bold{u}^\ast_k)-\strain(\bold{u}^\ast_\infty))\|_{L^2(\Omega)^{2\times2}}\to0, \quad \forall\bold{v}\in \bold{V}^D.
\]
Therefore,
\[
\int_\Omega (\rho_k^\ast)^p\C_0 \strain(\bold{u}^\ast_k) : \strain(\bold{v}) \dd x \to \int_\Omega (\rho_\infty^\ast)^p\C_0 \strain(\bold{u}^\ast_\infty) :
\strain(\bold{v}) \dd x, \quad \forall\bold{v}\in \bold{V}^D.
\]
This and \eqref{residual02} in Lemma \ref{lem:residual} imply the desired variational equation, which completes the proof of the lemma.
\end{proof}

Finally, Theorem \ref{thm:conv_afem} follows directly from Theorem \ref{thm:medmin} and Lemma \ref{lem:medcont}.

\begin{remark}\label{rem_conv}
In the convergence analysis, we have assumed that $p>1$ in the tensor $\C$ is an integer.
This assumption is only used in the proof of Lemma \ref{lem:err->0}. The rest of convergence
analysis does not rely on this assumption.
\end{remark}

\section{Conclusions}
In this work we have developed a novel adaptive phase-field approximation for compliance optimization, one fundamental problem in topology optimization. The approach is based on the phase-field representation of the medium (material versus void), and adaptive finite element discretization of the phase-field function and the displacement field. The estimators driving the adaptive refinement procedure consist of two parts: one for the objective functional (and the first-order necessary condition), and the other for the displacement field approximation. We provided a convergence analysis of the approximation in the sense that the sequence of the optimal states and designs contains a convergent subsequence to a solution of the first-order necessary optimality system. Furthermore, we presented numerical experiments which show clearly the accuracy and efficiency of the proposed adaptive phase-field approach. It is observed that the adaptive approach can yield optimal designs with comparable objective values but sharper and smoother interfaces at a lower computational cost, and hence it has significant potential for topology optimization.

\bibliographystyle{abbrv}
\bibliography{top}

\begin{thebibliography}{10}

\bibitem{AinsworthOden:2000}
M.~Ainsworth and J.~T. Oden.
\newblock {\em A {P}osteriori {E}rror {E}stimation in {F}inite {E}lement
  {A}nalysis}.
\newblock Pure and Applied Mathematics. Wiley-Interscience, John Wiley \& Sons,
  New York, 2000.

\bibitem{Aurada:2012}
M.~Aurada, S.~Ferraz-Leite, and D.~Praetorius.
\newblock Estimator reduction and convergence of adaptive {BEM}.
\newblock {\em Appl. Num. Math.}, 62:787--801, 2012.

\bibitem{Bendsoe:1989}
M.~P. Bends{\o}e.
\newblock Optimal shape design as a material distribution problem.
\newblock {\em Struct. Opt.}, 1(4):193--202, 1989.

\bibitem{BendsoeKikuchi:1988}
M.~P. Bends{\o}e and N.~Kikuchi.
\newblock Generating optimal topologies in structural design using a
  homogenization method.
\newblock {\em Comput. Methods Appl. Mech. Engrg.}, 71(2):197--224, 1988.

\bibitem{BendsoeSigmund:1999}
M.~P. Bends\o{e} and O.~Sigmund.
\newblock Material interpolation schemes in topology optimization.
\newblock {\em Arch. Appl. Mech.}, 69(9--10):635--645, 1999.

\bibitem{BendsoeSigmund:2003}
M.~P. Bends{\o}e and O.~Sigmund.
\newblock {\em Topology {O}ptimization: {T}heory, {M}ethods and
  {A}pplications}.
\newblock Springer-Verlag, Berlin, 2003.

\bibitem{BlankGarckeHecht:2016}
L.~Blank, H.~Garcke, C.~Hecht, and C.~Rupprecht.
\newblock Sharp interface limit for a phase field model in structural
  optimization.
\newblock {\em SIAM J. Control Optim.}, 54(3):1558--1584, 2016.

\bibitem{BourdinChambolle:2003}
B.~Bourdin and A.~Chambolle.
\newblock Design-dependent loads in topology optimization.
\newblock {\em ESAIM Control Optim. Calc. Var.}, 9:19--48, 2003.

\bibitem{BourdinChambolle:2006}
B.~Bourdin and A.~Chambolle.
\newblock The phase-field method in optimal design.
\newblock In M.~P. Bends{\o}e, N.~Olhoff, and O.~Sigmund, editors, {\em {IUTAM
  Symposium on Topological Design Optimization of Structures, Machines and
  Materials}}, pages 207--215. Springer, Berlin, 2006.

\bibitem{Braess:2007}
D.~Braess.
\newblock {\em {Finite Elements: Theory, Fast Solvers, and Applications in
  Elasticity Theory}}.
\newblock Cambridge University Press, Cambridge, third edition, 2007.

\bibitem{BrennerScott:2008}
S.~C. Brenner and L.~R. Scott.
\newblock {\em The {M}athematical {T}heory of {F}inite {E}lement {M}ethods}.
\newblock Springer, New York, third edition, 2008.

\bibitem{BruggiVerani:2011}
M.~Bruggi and M.~Verani.
\newblock A fully adaptive topology optimization algorithm with goal-oriented
  error control.
\newblock {\em Comput. Struct.}, 89(15--16):1481--1493, 2011.

\bibitem{BurgerStainko:2006}
M.~Burger and R.~Stainko.
\newblock Phase-field relaxation of topology optimization with local stress
  constraints.
\newblock {\em SIAM J. Control Optim.}, 45(4):1447--1466, 2006.

\bibitem{CahnHilliard:1958}
J.~Cahn and J.~Hilliard.
\newblock Free energy of a non-uniform system i--interfacial free energy.
\newblock {\em J. Chem. Phys.}, 28(2):258--267, 1958.

\bibitem{Ciarlet:2002}
P.~G. Ciarlet.
\newblock {\em The {F}inite {E}lement {M}ethod for {E}lliptic {P}roblems}.
\newblock SIAM, Philadelphia, PA, 2002.

\bibitem{Ciarlet:2013}
P.~G. Ciarlet.
\newblock {\em Linear and {N}onlinear {F}unctional {A}nalysis with
  {A}pplications}.
\newblock SIAM, Philadelphia, PA, 2013.

\bibitem{CostaAlves:2003}
J.~C.~A. {Costa Jr} and M.~K. Alves.
\newblock Layout optimization with $h$--adaptivity of structures.
\newblock {\em Int. J. Numer. Method Eng.}, 58(1):83--102, 2003.

\bibitem{SturlerPaulino:2008}
E.~{de Sturler}, G.~H. Paulino, and S.~Wang.
\newblock Topology optimization with adaptive mesh refinement.
\newblock In J.~F. Abel and J.~R. Cooke, editors, {\em {Proceedings of the
  Sixth International Conference on Computation of Shell and Spatial Structures
  IASS-IACM 2008: Spanning Nano to Mega}}, Ithaca, NY, USA, 2008.

\bibitem{TroyaTortorelli:2018}
M.~A.~S. {de Troya} and D.~A. Tortorelli.
\newblock Adaptive mesh refinement in stress-constrained topology optimization.
\newblock {\em Struct. Multidisc. Optim.}, 58(6):2369--2386, 2018.

\bibitem{DedeBordenHughes:2012}
L.~Ded\`e, M.~J. Borden, and T.~J.~R. Hughes.
\newblock Isogeometric analysis for topology optimization with a phase field
  model.
\newblock {\em Arch. Comput. Methods Eng.}, 19(3):427--465, 2012.

\bibitem{DiazSigmund:1995}
A.~D\'{i}az and O.~Sigmund.
\newblock Checkerboard patterns in layout optimization.
\newblock {\em Struct. Optim.}, 10(1):40--45, 1995.

\bibitem{EschenauerOlhoff:2001}
H.~A. Eschenauer and N.~Olhoff.
\newblock Topology optimization of continuum structures: A review.
\newblock {\em Appl. Mech. Rev.}, 54(4):331--390, 2001.

\bibitem{EvansGariepy:2015}
L.~C. Evans and R.~F. Gariepy.
\newblock {\em Measure Theory and Fine Properties of Functions}.
\newblock CRC Press, Boca Raton, FL, revised edition, 2015.

\bibitem{Fleury:1989}
C.~Fleury.
\newblock {CONLIN}: an efficient dual optimizer based on convex approximation
  concepts.
\newblock {\em Struct. Optim.}, 1(2):81--–89, 1989.

\bibitem{GLNS:2023}
H.~Garcke, K.~F. Lam, R.~Nürnberg, and A.~Signori.
\newblock Phase field topology optimisation for 4{D} printing.
\newblock {\em ESAIM Control Optim. Calc. Var.}, 29:24, 46pp, 2023.

\bibitem{GuoBabuska:1993}
B.-Q. Guo and I.~Babu\v{s}ka.
\newblock On the regularity of elasticity problems with piecewise analytic
  data.
\newblock {\em Adv. in Appl. Math.}, 14(3):307--347, 1993.

\bibitem{HQZ}
X.~Hu, M.~Qian, and S.~Zhu.
\newblock Accelerating a phase field method by linearization for eigenfrequency
  topology optimization.
\newblock {\em Struct. Multidisc. Optim.}, 66:242, 17 pp., 2023.

\bibitem{ItoJin:2015}
K.~Ito and B.~Jin.
\newblock {\em Inverse {P}roblems: Tikhonov {T}heory and {A}lgorithms}.
\newblock World Scientific Publishing Co. Pte. Ltd., Hackensack, NJ, 2015.

\bibitem{ItoKunisch:2008}
K.~Ito and K.~Kunisch.
\newblock {\em {Lagrange Multiplier Approach to Variational Problems and
  Applications}}.
\newblock SIAM, Philadelphia, PA, 2008.

\bibitem{JinXu:2019}
B.~Jin and Y.~Xu.
\newblock Adaptive reconstruction for electrical impedance tomography with a
  piecewise constant conductivity.
\newblock {\em Inverse Problems}, 36(1):014003, 28pp, 2020.

\bibitem{LambeCzekanski:2018}
A.~B. Lambe and A.~Czekanski.
\newblock Topology optimization using a continuous density field and adaptive
  mesh refinement.
\newblock {\em Int. J. Numer. Methods Eng.}, 113(3):357--373, 2018.

\bibitem{LazarovSigmund:2011}
B.~S. Lazarov and O.~Sigmund.
\newblock Filters in topology optimization based on {H}elmholtz-type
  differential equations.
\newblock {\em Int. J. Numer. Methods Eng.}, 86(6):765--781, 2011.

\bibitem{LiYang:2022}
F.~Li and J.~Yang.
\newblock A provably efficient monotonic-decreasing algorithm for shape
  optimization in stokes flows by phase-field approaches.
\newblock {\em Comput. Methods Appl. Mech. Engrg.}, 398:115195, 2022.

\bibitem{Modica:1987}
L.~Modica.
\newblock The gradient theory of phase transitions and the minimal interface
  criterion.
\newblock {\em Arch. Rational Mech. Anal.}, 98(2):123--142, 1987.

\bibitem{ModicaMortola:1977b}
L.~Modica and S.~Mortola.
\newblock Il limite nella {$\Gamma $}-convergenza di una famiglia di funzionali
  ellittici.
\newblock {\em Boll. Un. Mat. Ital. A (5)}, 14(3):526--529, 1977.

\bibitem{ModicaMortola:1977}
L.~Modica and S.~Mortola.
\newblock Un esempio di {$\Gamma ^{-}$}-convergenza.
\newblock {\em Boll. Un. Mat. Ital. B (5)}, 14(1):285--299, 1977.

\bibitem{NochettoSiebertVeeser:2009}
R.~H. Nochetto, K.~G. Siebert, and A.~Veeser.
\newblock Theory of adaptive finite element methods: an introduction.
\newblock In {\em {Multiscale, Nonlinear and Adaptive Approximation}}, pages
  409--542. Springer, Berlin, 2009.

\bibitem{NovotnySokołowski:2013}
A.~A. Novotny and J.~Sokołowski.
\newblock {\em Topological Derivatives in Shape Optimization}.
\newblock Springer, Heidelberg, 2013.

\bibitem{Plotnikov:2023}
P.~I. Plotnikov and J.~Sokolowski.
\newblock Geometric aspects of shape optimization.
\newblock {\em J. Geom. Anal.}, 33(7):206, 57pp, 2023.

\bibitem{QianHuZhu:2012}
M.~Qian, X.~Hu, and S.~Zhu.
\newblock A phase field method based on multi-level correction for eigenvalue
  topology optimization.
\newblock {\em Comput. Methods Appl. Mech. Engrg.}, 401:115646, 2022.

\bibitem{Razvany:2009}
G.~I.~N. Rozvany.
\newblock A critical review of established methods of structural topology
  optimization.
\newblock {\em Struct. Multidisc. Optim.}, 37(3):217--237, 2009.

\bibitem{ScottZhang:1990}
L.~R. Scott and S.~Zhang.
\newblock Finite element interpolation of nonsmooth functions satisfying
  boundary conditions.
\newblock {\em Math. Comp.}, 54(190):483--493, 1990.

\bibitem{SigmundPetersson:1998}
O.~Sigmund and J.~Petersson.
\newblock Numerical instabilities in topology optimization: A survey on
  procedures dealing with checkerboards, mesh-dependencies and local minima.
\newblock {\em Struct. Optim}, 16(1):68--75, 1998.

\bibitem{Stainko:2006}
R.~Stainko.
\newblock An adaptive multilevel approach to the minimal compliance problem in
  topology optimization.
\newblock {\em Commun. Numer. Methods Eng.}, 22(2):109--118, 2006.

\bibitem{StolpeSvanberg:2001}
M.~Stolpe and K.~Svanberg.
\newblock An alternative interpolation scheme for minimum compliance topology
  optimization.
\newblock {\em Struct. Multidisc. Optim.}, 22(2):116--126, 2001.

\bibitem{Svanberg:1987}
K.~Svanberg.
\newblock Method of moving asymptotes –- a new method for structural
  optimization.
\newblock {\em Int. J. Numer. Meth. Eng.}, 24(3):359--373, 1987.

\bibitem{TakezawaNishiwakiKitamura:2010}
A.~Takezawa, S.~Nishiwaki, and M.~Kitamura.
\newblock Shape and topology optimization based on the phase field method and
  sensitivity analysis.
\newblock {\em J. Comput. Phys.}, 229(7):2697--2718, 2010.

\bibitem{Tavakoli:2014}
R.~Tavakoli.
\newblock Multimaterial topology optimization by volume constrained
  {A}llen-{C}ahn system and regularized projected steepest descent method.
\newblock {\em Comput. Methods Appl. Mech. Engrg.}, 276:534--565, 2014.

\bibitem{Verfurth:2013}
R.~Verf\"{u}rth.
\newblock {\em A {P}osteriori {E}rror {E}stimation {T}echniques for {F}inite
  {E}lement {M}ethods}.
\newblock Oxford University Press, Oxford, 2013.

\bibitem{WallinRistinmaa:2013}
M.~Wallin and M.~Ristinmaa.
\newblock Howard's algorithm in a phase-field topology optimization approach.
\newblock {\em Internat. J. Numer. Methods Eng.}, 94(1):43--59, 2013.

\bibitem{WallinRistinmaa:2014}
M.~Wallin and M.~Ristinmaa.
\newblock Boundary effects in a phase-field approach to topology optimization.
\newblock {\em Comput. Methods Appl. Mech. Engrg.}, 278:145--159, 2014.

\bibitem{WangZhou:2004}
M.~Y. Wang and S.~Zhou.
\newblock Phase field: a variational method for structural topology
  optimization.
\newblock {\em Comput. Model. Eng. Sci.}, 6(6):227--246, 2004.

\bibitem{ZhouRozvany:1991}
M.~Zhou and G.~I.~N. Rozvany.
\newblock The {COC} algorithm, {Part II}: Topological, geometrical and
  generalized shape optimization.
\newblock {\em Comput. Methods Appl. Mech. Engrg.}, 89(1--3):309--336, 1991.

\bibitem{ZhouWang:2007}
S.~Zhou and M.~Y. Wang.
\newblock Multimaterial structural topology optimization with a generalized
  {Cahn--Hilliard} model of multiphase transition.
\newblock {\em Struct. Multidisc. Optim.}, 33(2):89--111, 2007.

\end{thebibliography}

\end{document}